\documentclass[1p,preprint]{elsarticle}
\usepackage{ecrc_arxiv1}
\volume{00} 
\firstpage{1} 

\makeatletter
\renewcommand\paragraph{\@startsection{paragraph}{4}{\z@}%
           {12\p@ \@plus 6\p@ \@minus 3\p@}%
           {\p@}%
           {\normalfont\normalsize\itshape}}

\makeatother

\usepackage{algorithmic}
\usepackage{amssymb}
\usepackage{amsthm}
\usepackage[title]{appendix}
\usepackage{array}
\usepackage{cancel}
\usepackage{dsfont}
\usepackage{multirow}

\usepackage{epsfig}
\usepackage{calc}
\usepackage{amssymb}
\usepackage{amstext}
\usepackage{amsmath}
\usepackage{hyperref}
\usepackage{lineno}
\usepackage{multicol}
\usepackage{orcidlink}
\usepackage{pslatex}
\usepackage{setspace}
\definecolor{DodgerBlue}{rgb}{0.12, 0.56, 1.0}
\definecolor{Cerise}{rgb}{1, 0.3, 1.0}

\usepackage{xcolor}

\newcommand*{\colorboxed}{}
\def\colorboxed#1#{%
  \colorboxedAux{#1}%
}
\newcommand*{\colorboxedAux}[3]{%
  \begingroup
    \colorlet{cb@saved}{.}%
    \color#1{#2}%
    \boxed{%
      \color{cb@saved}%
      #3%
    }%
  \endgroup
}

\newtheorem{prop}{Proposition}

\newtheorem{definition}{Definition}

\usepackage{colortbl}
\usepackage{svg}
\usepackage{empheq}
\usepackage{tikz,calc}
\usetikzlibrary{math}
\tikzset{
	pics/carc/.style args={#1:#2:#3}{
		code={
			\draw[pic actions] (#1:#3) arc(#1:#2:#3);
		}
	}
}

\usepackage{empheq}
\definecolor{DodgerBlue}{rgb}{0.12, 0.56, 1.0}
\definecolor{Cerise}{rgb}{1, 0.3, 1.0}
\definecolor{Green}{rgb}{0.3, 0.6, 0.45}

\definecolor{vc}{RGB}{155, 89, 182}
\definecolor{vertc}{RGB}{167, 224, 45}
\definecolor{bleu}{RGB}{61, 57, 196}
\usepackage{xcolor}
\usepackage{diagbox}

\newtheorem{remark}{Remark}

\begin{document}
\begin{frontmatter}

\title{Co-simulation domain decomposition algorithm for hybrid EMT-Dynamic Phasor modeling}
\author{Hélèna Shourick \fnref{auth1}}
\address{SuperGrid Institute, 23 rue Cyprian,  69200 Villeurbanne}
\fntext[auth1]{helena.shourick@supergrid-institute.com}

\author{Damien Tromeur-Dervout \fnref{auth2}}
\address{University of Lyon, UMR5208 U.Lyon1-CNRS,  Institut Camille Jordan, 15 Bd Latarjet 69622 Villeurbanne}
\fntext[auth2]{damien.tromeur-dervout@univ-lyon1.fr}

\author{Laurent Chédot\fnref{auth3}}
\address{SuperGrid Institute, 23 rue Cyprian,  69200 Villeurbanne}
\fntext[auth3]{Laurent.chedot@supergrid-institute.com}

\begin{abstract}
An iterative coupling algorithm based on a restricted additive Schwarz domain decomposition is investigated to co-simulate electrical circuits  with  hybrid electromagnetic (EMT) and transient stability (TS) modeled using dynamic phasors. This co-simulation algorithm does not introduce any delay between the data exchanged at the co-simulation step. The pure linear convergence property of the iterative method allows it to be accelerated towards the true solution by a non-intrusive Aitken's acceleration  of the convergence post-processing, even if the domain decomposition interface conditions make the iterative method divergent. This provides a method less sensitive to the splitting.  Numerical tests on a linear RLC circuit combining EMT and TS modeling are provided.
\end{abstract}

\runauth{H. Shourick et al.}

\begin{keyword}
Co-simulation  \sep Restricted additive Schwarz \sep Aitken's convergence acceleration
\MSC[2010] 65\sep  65B05 \sep 65L80 \sep 65M55 \sep 68U20
\end{keyword}

\end{frontmatter}



\section{Introduction}
The introduction of renewable energies into the power grid leads to the use of more components based on power electronics. These components imply faster dynamics. Power system safety simulations, which cannot be handled by traditional Transient Simulations (TS) conducted with dynamic phasors (DP), require Electro-Magnetic Transient (EMT) simulations. Nevertheless,  the advantage of TS programs is their computational speed which makes them suitable for handling large-scale networks, however, their modeling is not sufficiently detailed and can only catch slow dynamics. On the other hand, EMT simulators can capture fast dynamics, but are limited in computational speed; therefore, they are used to simulate only small portions of the network. For large power grids, it can be expected that the need for high-level detail requiring Electro-Magnetic Transient (EMT) modeling will be localized near disturbances, and other parts of the network will use TS modeling. 

Based on this assumption, the co-simulation approach is an attractive candidate to handle these hybrid power system simulations. Nevertheless, the EMT-TS co-simulation has to face several locks as already underlined in \cite{IEEEinterfaceTSEMT} among which we can mention:
\begin{itemize}
\item the data exchange between TS and EMT simulators, including choice of the interface variables (i.e. the partitioning of the network), the   
data conversion between waveform and dynamic phasors or phasors;
\item the interaction protocol between TS and EMT simulators, including the time step size difference between EMT and TS, the coupling algorithm (i.e. iterative or non iterative, the signal rebuilding Zero order hold or high order hold), the convergence of the resulting coupling algorithm.
\end{itemize}
Advances in these areas have been proposed in recent years. Among these is the partitioning based on traveling wave patterns of the transmission line which introduce natural decoupling into the nodal equations of an EMT simulator due to transmission line latency. Note that this latency limits the flexibility of choosing DP time steps \cite{Mudunkotuwa2018}. Le-Huy \& al \cite{Le-Huy2017}  developed a simple hybrid line model that accounts for wave propagation in both electromagnetic and TS simulations i.e both ends exchange historical current and the delay is considered  as a phase shift for the TS part. Another co-simulation that splits the transmission line is the one of Rupasinghe \& al \cite{Rupasinghe2020} in which a Base-Frequency Phasor Adaptive Simulation Transient solver is developed, which allows to derive frequency-dependant equivalent of network components using base-frequency dynamic phasors at the discrete level.
Then the offset frequency parameter can switch from EMT to DP at the fundamental frequency. Plumier \& al \cite{plumier2016WR}has proposed a co-simulation algorithm coupling EMT and phasor which dynamically updates by an iterative procedure the equivalent impedance of Norton and Thevenin equivalent models representing  the boundary conditions between each subsystem. Some acceleration of the coupling algorithm is achieved by a prediction scheme which is similar to a Richardson extrapolation of temporal quantities.
 
Shu \& al \cite{Shu2018} proposed a two-level Schur complement through which the Thevenin equivalent obtained for each EMT subsystem could fully consider the coupling among different EMT subsystems and the TS central system. Abhyankar and Flueck \cite{Abhyankar2012} proposed an implicitly coupled\\ TSEMT algorithm where the set of TS and EMT equations are solved by a Newton's method at each TS time step. An important issue in this approach is the computation of the instantaneous Thevenin equivalent voltage for EMT time steps that are not on the temporal boundary. Their experiment has shown that it is peferable to use the TS solution at the end TS time step rather  than a linear interpolation between the two TS interval steps, probably due to the Newton's algorithm. Rimorov \& al \cite{PCoupling} focused on the problems of co-simulation stability and precision in the presence of delays and proposed a generalized interface framework related to the search for a power-conjugate interface that combines current and voltage through an "impedence" parameter resulting from non-physical related boundary conditions.\\
\indent In co-simulation algorithms such as non-iterative Jacobi, zero-order hold iterative co-simulation and non-iterative algorithm improving variables smoothing, the delay of one co-simulation step  (i.e. TS time step delay) between the given inputs and the retrieved outputs of the TS and EMT systems can lead to instabilities. Some iterative techniques such as the fixed-point method \cite{Eguillon2019IFOSMONDI} or the Newton-like method \cite{plumier2016WR,Abhyankar2012} can, even with a high order smoothing constraints,  solve the so-called ”constraint function” corresponding to the interface of the systems \cite{Eguillon2021IFOSMONDIJFM}.

\indent In this paper, we consider a fixed-point coupling algorithm based on the Schwarz domain decomposition technique which can be related to the dynamic iteration method (DI) (i.e waveform relaxation of \cite{Lelarasmee1982}) in which we used a restricted additive Schwarz (RAS) splitting \cite{Shourick2022}. These DI methods can be convergent or divergent depending on the domain partitioning and boundary conditions. Nevertheless, we used the good property of purely linear convergence or divergence (i.e. the error operator of the method does not depend on the iteration number) to accelerate the iterative method towards the true solution  with the Aitken's acceleration of convergence technique even with a divergent method \cite{GTD_IJNMF02,TromeurEssaim2013}. The advantages of our approach are:
\begin{itemize}
\item the fundamental concept of the method is to post-process the sequence of interfaces solutions generated by the domain decomposition solver. It can use different boundary conditions for the acceleration as long as they are linear in the variables (i.e Dirichlet , Neumann, Robin, those of \cite{PCoupling},...);
\item as it is a post-process the method can be used non-intrusively in the local solver;
\item it is not necessary for domain partitioning to be cut on the transmission lines;
\item it can also support an overlap of TS and EMT parts. We then have some components with the two representations.
\end{itemize}

The outline of the paper is as follows: section \ref{shourick_contrib_Sec2} describes the EMT and TS (modeled with dynamic phasors) modeling of the electrical network. Section \ref{shourick_contrib_Sec3} presents  the co-simulation  algorithm consisting on the heterogeneous Restricted Additive Schwarz domain decomposition  with a special attention to the translation operators between the EMT and TS  RAS iterates. As the convergence of the RAS can depends of the electrical network components, section \ref{shourick_contrib_Sec4} establishes the heterogeneous RAS EMT-TS error operator and defines the acceleration of its convergence with the  Aitken's acceleration of the convergence technique. 
 Section \ref{shourick_contrib_Sec5} presents results obtained on a linear RLC circuit before concluding in section  \ref{shourick_contrib_Sec6}.

\section{Mathematical modeling of the electrical network EMT and Dynamic Phasor modelling \label{shourick_contrib_Sec2}}

The electrical network can be view as a graph connecting electrical components through their connecting pins. The vertices or nodes of this graph are the pins of electrical components and its  edges  are the link between vertices of two connected components. Some physical quantities are defined on this pins such as currents and voltages. The common principle of different mathematical modeling  is based on the application of Kirshoff's laws, which establish the mathematical relationships between the different physical quantities of an electrical network. Kirshoff's laws are as follows:
\begin{itemize}
\item Kirchhoff's current law: The sum of the currents entering a node is equal to the sum of the currents leaving this node.
\item Kirchhoff's voltage law: The sum of the voltages around any closed loop is zero.
\end{itemize}

The nodal analysis expresses the potential at each node using Kirshkoff's laws and the component properties of the branches connected to that node. It creates an admittance matrix linking the sum of the current entering each node and the voltages. However, the nodal formulation does not allow to directly represent certain devices as current-dependent circuit element \cite{Mana}. The  Modified Nodal Analysis, introduced by \cite{Mana}, widely used in network modeling since \cite{shourick_contrib_MANA} and improved especially for modeling electromagnetic transients\cite{MAHSEREDJIAN20071514}, allows one to overcome those difficulties. This way the complete network is written in the form $$A_{N_t} z_{N_t}=b_{N_t}$$ where the admittance matrix of the nodal analysis is included in the matrix $A_{N_t}$ where $A_{N_t}$ is the linearized matrix at time $t$, if there is nonlinear devices. The $z_{N_t}$ are the current and voltages unknowns and $b_{N_t}$ contains the knowns current and voltages at time $t$. Notes that in this formulation the time discretization has already been performed.

The more general mathematical formulation of the electrical network considers the building of the differential algebraic equations system induced by the Kirshoff's law and the electrical  network's components modeled using differential equations: 
\begin{equation}
F(t,x(t),\dot{x}(t),y(t))=0 \label{generalDAE}
\end{equation}
 where $ x(t)\in \mathbb {R} ^{n}$ are the differential unknowns, $ \dot{x}(t)\in \mathbb {R} ^{n}$ are the derivative of $x$ with respect to time and $y(t)\in \mathbb {R} ^{m}$ are the algebraic unknowns. Some tools based on Modelica language generate such DAE system.

One can sometime separate in the DAE system the purely algebraic equations from the others. This way, the general DAE system \eqref{generalDAE} can be rewritten as:
 \begin{equation}
 \left\{ 
 \begin{array}{lcl}
\dot{x}(t)&=& f(t,x(t),y(t)) \\
 0&=&g(t,x(t),y(t))
 \end{array} 
 \right. \label{generalDAEsplit}
\end{equation}

Eq. \eqref{generalDAEsplit} is a nonlinear system. Its  state space representation  consists on the linearizing around a time $t_n$ of functions $f$and $g$ producing the state matrices $A$,$B$,$C$,$D$. The resulting linear system writes for $t\in [t_n,t_{n+1}]$ 
  \begin{equation}
 \begin{array}{lcl}
\dot{v}(t)&=& {A}v(t)+ {B}u(t) \\
w(t)&=&{C}v(t)+{D}u(t)
 \end{array} \label{generalDAEStateSpace}
 \end{equation} 
where $v$ are the state variables, $u$ the inputs and $w$ the outputs. These variables can be voltages and/or currents. This representation Eq. \eqref{generalDAEStateSpace} is those used in the FMI standard where the connected Functional Mock-up Units (FMU) are black box differential systems producing outputs reacting to inputs coming from others FMUs.\\
 \medskip

The mathematical modeling of the electrical network depends  also of the nature of transient phenomena that must be caught. Some hypothesis or not on the shape of the unknowns can lead to different levels of mathematical modeling. 
  
If no hypothesis is made on the shape of unknowns, then the DAE system  formulated by system \eqref{generalDAE} has to be solved to address some EMT stability. The high dynamics present constraint strongly the time step of numerical DAE solver. If a strong hyppothesis on the shape of the unknowns is made then  the system Eq. \eqref{generalDAE} can be simplified.\\
\indent The phasor modeling considers the electrical unknowns as sinusoidal signal with a given constant pulsation $\omega_0=2\pi f_0$ (expressed in radians/second) and a phase $\theta$ and can be represented by a complex constant amplitude. This frequency $f_0=\frac{1}{T_0}$ is normally 50Hz or 60Hz and $t_0$ is the period (expressed in second).
Due to the Fourier transform, a periodic waveform $x(\tau)$ can be
written on the interval $\tau\in[t-T,t]$ where $T$ is the observation period considered,  as $x(\tau)=\sum_{k=-\infty}^{+\infty}X_{k}e^{jk\omega_{0}\tau}$ where the Fourier coefficients $X_k$ are the phasors amplitude. The time derivative then applies on the sinusoidal part of the phasor. It results a simplification of system \eqref{generalDAE} in a linear system. If steady-state transients are searched, the time derivatives in system \eqref{generalDAE} can be omitted.

A compromise between phasor simulation and EMT simulationb is the Dynamic Phasors (DP)  simulation. A dynamic phasor is a phasor  whose amplitude $X_k$ and phase angle $\theta$ are time-dependent values by considering that the  waveform is not strictly periodic (in an almost periodic state). Their definition and properties are the following:
\begin{eqnarray}
\langle x\rangle_{k}&\stackrel{def}{=}&X_{k}(t)=\frac{1}{T}\int_{t-T}^{^{t}}x(\tau)e^{-jkw_{0}\tau}d\tau  \label{DynamicPhasor}\\
\langle\dfrac{dx}{dt}\rangle_{k}&=&\frac{d\langle x\rangle_{k}}{dt}+jk\omega_{0}\langle x\rangle_{k} \nonumber\\
\langle xy\rangle_{k}&=&\sum_{i}\langle x\rangle_{k-i}\langle y\rangle_{i}\nonumber
\end{eqnarray}

 The DAE system \eqref{generalDAE} is transformed in another DAE system with putting the each unknown as a sum of dynamic phasors. The sum range depends on the harmonic kept.  As harmonics are compute separately there is a multiplication of variables. Despite this enlargement of the DAE system, as time varying Fourier coefficients are slower than the original values, it allows time step much bigger than EMT (2 to 30 times larger). Demiray \cite{Demiray2008SimulationOP} discusses on the feasibility of using Dynamic phasor to simulate large networks, dynamic phasor simulations type allows to catch dynamics up to 60Hz with a computation time reasonable for large networks.  These results must be tempered by those of Hassani \& al \cite{HASSANI} that did not find advantage of DP over a simulation in the time domain.

Table \ref{typeSimuTable} is a summary of the modeling of the three type of simulations (EMT, Dynamic Phasors, Phasors) for the basic components.
\begin{table}[h!]
\begin{center}
\footnotesize
\begin{tabular}{|c|c|c|c|}
\hline 
 & EMT & Dynamic Phasor & Phasor\tabularnewline
\hline
system & $F(t,x(t),\dot{x}(t),y(t))=0$ &  $\langle F  (t,\langle x \rangle(t),\dot{\langle x \rangle (t)} ,\langle y \rangle (t))=0$ & $\tilde{F}(t,\textbf{X}(t),\textbf{Y}(t))=0$ \tabularnewline
\hline 
variable & free shape & $\langle x \rangle(t)=\sum_{k=0}^{m} \bar{x}_k(t) e^{j\omega_k t + \theta_k(t)}$ & $\textbf{X}(t)=\sum_{k=0}^{m} \bar{\textbf{X}}_k e^{j\omega_k t + \theta_k}$ \tabularnewline
\hline
resistance & $u=Ri$ & $\tilde{u}_{k}=R\tilde{i}_{k}$ & $\tilde{u}=R\tilde{i}$\tabularnewline
\hline 
inductor & $u=L\frac{\ensuremath{di}}{dt}$ & $\langle u \rangle _{k}=L \langle \dot{ i_{k} \rangle}+L \langle i \rangle _{k}kj\omega_{o}$ & $\tilde{ u}=L\tilde{i}j\omega_{o}$\tabularnewline
\hline 
 capacitor & $i=C\frac{\ensuremath{du}}{dt}$ & $\langle i\rangle _{k}=C\dot{\langle u_{k} \rangle }+C\langle u \rangle_{k}kj\omega_{o}$ & $\tilde{ i}=C\tilde{u} j\omega_{o}$\tabularnewline
\hline 
\end{tabular}
\end{center}
\caption{\label{typeSimuTable} EMT, phasor and Dynamic phasor representation}
\end{table}

\section{Co-simulation algorithm \label{shourick_contrib_Sec3}}

Let us consider the DAE system that follows where, for all $t \in [0, T]$,   $x(t)\in \mathbb{R}^{n_1}$ are the differential unknowns  and $y(t) \in  \mathbb{R}^{n_2}$  the algebraic ones: 
\begin{equation}
\left\{
\begin{array}{ll}
  F(t,\dot{x},x,y)&=0 \\
  g(t,x,y)&=0\\
       x(0)&=x_0\\
       y(0)&=y_0
    \end{array}
\right.  \label{DAEmono1}
\end{equation}
The variables $x$ and $y$ can represent voltage or current following the electrical components involved in the electrical circuit.
 
Our goal is to solve this DAE system by splitting it into several parts (at least two), i.e. splitting the set of unknowns into several subsets by gathering the differential or algebraic equations associated with the unknowns belonging to the same subset. For example, if we consider two subsets the orginal DAE system will be split into two DAE systems to be solved on the time interval $[T_n^+, T_{n+1}^-]$:
\begin{eqnarray}
\left\{
\begin{array}{ll}
  F_1(t,\dot{x}_1,x_1,y_1,\tilde{x}_2,\tilde{y}_2)&=0 \\
  g_1(t,x_1,y_1,\tilde{x}_2,\tilde{y}_2)&=0\\
       x_1(T_n^+)&=x^n_1\\
       y_1(T_n^+)&=y^n_1
    \end{array} 
\right.   &&
\left\{
\begin{array}{ll}
  F_2(t,\dot{x}_2,x_2,y_2,\tilde{x}_1,\tilde{y}_1)&=0 \\
  g_2(t,x_2,y_2,\tilde{x}_1,\tilde{y}_1)&=0\\
       x_2(T_n^+)&=x^n_2\\
       y_2(T_n^+)&=y^n_2
    \end{array} 
\right.  
\end{eqnarray}
Where $\tilde{x}_2$ and $\tilde{y}_2$ in the DAE subsystem 1 (respectively $\tilde{x}_1$ and $\tilde{y}_1$ in the DAE subsystem 2) are representations of the solutions $x_2$ and $y_2$ in the DAE subsystem 2 (respectively ${x}_1$ and ${y}_1$ in the DAE subsystem 1). They can be considered as inputs for the current DAE subsystem and must be updated at some point de rendez-vous time $T_n^+$ during the time simulation.

If $\tilde{x}_2=x_2$ and $\tilde{y}_2=y_2$ (respectively $\tilde{x}_1=x_1$ and $\tilde{y}_1=y_1$ ), the two DAE subsystems are said to be strongly coupled and can not be solved separately. Co-simulation techniques consist of having approximations for $\tilde{x}$ and $\tilde{y}$ such as Zero Order Hold (ZOH), where the $\tilde{x}$ and $\tilde{y}$ are frozen at their values at the time of the previous rendez-vous point. Some other polynomial approximations for $\tilde{x}$ and $\tilde{y}$ 
such as linear approximations from the values at previous rendez-vous point (First Order Hold) or with polynomials with higher degree (Second Order Hold or Third order Hold) can be used. Some extrapolation techniques with delay such as the $C(p,q,j)$ scheme of \cite{cpqj} where the input values are extrapolated to order $j^{th}$ from the solution's values taken at $p$ regular rendez-vous points in the past also exist.

The major drawback of such approaches, is the difference of the value of the solution of one subsystem and its representation in the other subsystem at the next point of rendez-vous time and it also limits the size of the macro time step separating two de rendez-vous points times. One solution to avoid this delay at the next rendez-vous point is to consider iterative algorithms that will update the inputs $\tilde{x}$ and $\tilde{y}$ in order that they have the same value as their value in the subsystem that computes them. 

Schematically the iterative algorithm on the time interval $[T^+_n,T^-_{n+1}]$ is written as follows: starting from initial inputs $\tilde{x}^{(0)}$ , $\tilde{y}^{(0)}$, the algorithm iterates over these values until they no longer change:
\begin{eqnarray}
\left\{
\begin{array}{ll}
  F_1(t,\dot{x}_1,x_1,y_1,\tilde{x}_2^{(k)},\tilde{y}_2^{(k)})&=0 \\
  g_1(t,x_1,y_1,\tilde{x}_2,\tilde{y}_2)&=0\\
       x_1(T_n^+)&=x^n_1\\
       y_1(T_n^+)&=y^n_1
    \end{array} 
\right.  , &&
\left\{
\begin{array}{ll}
  F_2(t,\dot{x}_2,x_2,y_2,\tilde{x}_1^{(k)},\tilde{y}_1^{(k)})&=0 \\
  g_2(t,x_2,y_2,\tilde{x}_1,\tilde{y}_1)&=0\\
       x_2(T_n^+)&=x^n_2\\
       y_2(T_n^+)&=y^n_2
    \end{array} 
\right. , \\
H_1(\tilde{x}_1^{(k+1)},\tilde{y}_1^{(k+1)}, x^{n+1}_1,y^{n+1}_1)=0, && H_2(\tilde{x}_2^{(k+1)},\tilde{y}_2^{(k+1)}, x^{n+1}_2,y^{n+1}_2)=0.
\end{eqnarray}
Where functions $H_1$ and $H_2$ are constraint functions on the inputs in order to guarantee the same values  of the inputs as the values that they  represent in the other subsystem at the end of the macro step simulation. The way these constraint functions are satisfied can lead to Newton type algorithms such as IFOSMONDI-JFM \cite{Eguillon2021IFOSMONDIJFM} or IFOSMONDI fixed-point type algorithm \cite{Eguillon2019IFOSMONDI}. Each of them has advantages and drawbacks, the most important drawback for the fixed-point algorithm is its non-contracting property leading in some cases to a non-convergent algorithm.

We focus in this work on the special choice of a fixed-point algorithm to satisfy the inputs constraint functions that is the Schwarz type domain decomposition method \cite{Schwarz}.

\subsection{Schwarz method for heterogeneous  EMT-TS}

In order to build the Schwarz co-simulation, we first take the representation of the discrete state space  Eq.\eqref{generalDAEStateSpace}  back and rewrite it for each of the representations  EMT and TS in order to fix the notations.
We assume that we have the TS representation of a domain $W$ as well as the EMT representation of the same domain. In simulations, it will not always be possible to have both TS and EMT representations on common parts of the network. But when it is possible, we can define an overlap and can compare the solutions obtained with the two local solvers. We assume  that the time step of the TS is a multiple of time step EMT: $\Delta t_{ts} = m \Delta t_{emt}$.

We consider a linear electrical circuit in the following,  non linear electrical circuits could be linearized with a state space representation and treated as a linear circuit over the time step. 
Let us rewrite the linear DAE \eqref{DAEmono1}  in its state space representation:
\begin{eqnarray}
\left\{\begin{array}{rcl}\mathbb{I} \dot{x}(t) +A x(t)+B y(t) &=& G_1(t),  x(0)=x_0,   \\
 C x(t) + D y(t)  &=& G_2(t),\, t\in [0,T]. \end{array} \right. \label{EqLinearDAE}
\end{eqnarray}

Where $x(t)\in \mathbb{R}^{n_1}$ and $y(t) \in  \mathbb{R}^{n_2}$ for all $t \in [0, T]$, $D$  is a $n_2 \times n_2$ nonsingular matrix, $A$ and $\mathbb{I}$ are  $n_1 \times n_1$ matrix , $\mathbb{I}$  matrix can be the identity or matrix composed of 1s and 0s depending on whether the $x$ variables contain voltages or potentials. $B$ is an $n_1 \times n_2$ matrix,
$C$ is an $n_2 \times n_1$ matrix, $G_1(t) \in \mathbb{R}^{n_1}$ and $G_2(t) \in \mathbb{R}^{n_2}$ are known input functions, as the DAE system is representing an electrical network, $G_1(t)$ and $G_2(t)$ are sources vector. Finally, $x_0 \in \mathbb{R}^{n_1}$ is a consistent initial value. Let $n = n_1 + n_2$.

We define the matrix $\mathbb{A}=\left(\begin{array}{cc} A & B \\ C & D \end{array} \right)$ corresponding to the linear operator of the DAE and we define $z(t)=[x(t),y(t)]^T$, $G(t)=[G_1(t), G_2(t)]^T$ and $\mathbb{I}_{d}=\left(\begin{array}{cc} \mathbb{I}_{n_1} & 0_{n_1\times n_2} \\ 0_{n_2 \times n_1} & 0_{n_2\times n_2} \end{array} \right)$. 
 \medskip
 \\
Then we can rewrite Eq. \eqref{EqLinearDAE} as:
\begin{eqnarray}
\begin{array}{rcl} \mathbb{I}_{d}\dot{z}(t) + \mathbb{A} z(t) &=& G(t),\;  x(0)=x_0, \; t\in [0,T]. \end{array}  \label{EqLinearDAE2}
\end{eqnarray}

The matrix $\tilde{\mathbb{A}}= \left( \begin{array}{cc}
\mathbb{I}_{n_1}   +\Delta t A &  \Delta t B\\
C & D \end{array} \right)\in \mathbb{R}^{n\times n}$ has a non-zero pattern and is associated  to the graph $G = (W, F)$, where the set of vertices $W = \left\{1,\ldots, n\right\}$ represents the $n$ unknowns and the set of edges $F = \left\{(i, j) | (\tilde{\mathbb{A}}_{i,j}) \neq 0\right\}$ represents the pairs of vertices that are coupled by a non-zero element in $\mathbb{A}$. Then, we assume that a graph partitioning was applied and that resulted in $N$ non-overlapping subsets $W_i^0$ whose union is $W$. 
Let $W_i^p$ be the $p$-overlap partition of $W$, obtained by including all the vertices immediately neighboring the vertices of $W_i^{p-1}$. Let $W_{i,e}^p=W_i^{p+1} \backslash W_i^{p}$. Then let $R_i^p \in \mathbb{R}^{n_i \times n}$  (  $R_{i,e}^p \in \mathbb{R}^{n_{i,e} \times n}$ and $\tilde{R}_i^0 \in \mathbb{R}^{n_i \times n}$ respectively) be the operator which restricts $w\in \mathbb{R}^n$ to the components of $w$ belonging to $W_i^p$   ($W_{i,e}^p$  and $W_i^0$ respectively, and the operator $\tilde{R}_i^0 \in \mathbb{R}^{n_i \times n}$ puts $0$ to the unknowns belonging to $W_i^p\backslash W_i^0$).  
Then we defined the local operators $\mathbb{A}_{i} = R_i^p \mathbb{A} R_i^{pT}$ and  $\mathbb{E}_{ie} = R_i^p \mathbb{A} R_{ie}^{pT}$ .\\

The DAE system is for a domain $W_i^p$, integrated between $t^n$ to $t^{n+1}$: 
\begin{eqnarray}
\underbrace{\left( \begin{array}{cc}
\mathbb{I}_{i}   +\Delta t A_{i} &  \Delta t B_{i}\\
C_{i} & D_{i} \end{array} \right)}_{\mathbb{A}_{i}} \underbrace{\left( \begin{array}{c}
x^{n+1}_{i} \\ y^{n+1}_{i} \end{array}\right)}_{z_{i}^{n+1}} =\underbrace{\left( \begin{array}{cc}
\mathbb{I}_{{i}}   &  0\\
0 & 0 \end{array} \right)}_{\mathbb{I}_{{d,i}}}\underbrace{ \left( \begin{array}{c}
x^{n}_{{i}}\\ y^{n}_{{i}} \end{array}\right)}_{z_{{i}} ^{n}}+&&  \nonumber\\
\underbrace{\left( \begin{array}{cc}
\mathbb{I}_{{ie}}   &  0\\ 0 & 0 \end{array} \right)}_{\mathbb{I}_{{d,ie}}}\underbrace{ \left( \begin{array}{c} x^{n}_{{ie}}\\ y^{n}_{{ie}} \end{array}\right)}_{z_{{ie}}^{n}} - \underbrace{\left( \begin{array}{cc}
\Delta t E^{A}_{ie} & \Delta t E^{B}_{ie}\\
E^{C}_{ie} & E^{D}_{ie}
\end{array}\right)}_{\mathbb{E}_{ie}} \underbrace{\left(\begin{array}{c}
x^{n+1}_{{i,e}} \\ y^{n+1}_{{i,e}}   \end{array}\right)}_{z_{{i,e}} ^{n+1}}+ \underbrace{\left( \begin{array}{c}\Delta t G^{n+1}_{1i} \\G^{n+1}_{2i}\end{array}\right)}_{G_i^{n+1}}. \label{EqLinearDiscretDAEWithDependency}
\end{eqnarray}

The term $\mathbb{I}_{{d,ie}} z_{{ie}}^{n}$ is coming from the fact that differential terms on the interface unknowns can be involved due to the splitting. The terms at time $t^n$ and the source term $G_i^{n+1}$ can be gather in a term $b_{i}^{n+1}$ independent of the solution $z_i^{n+1}$. Then the DAE system for the domain $W_i^p$ integrated between $t^n$ and $t^{n+1}$:

\begin{eqnarray}
\mathbb{A}_{i} z^{n+1}_{i}&=&b_{i}^{n+1} - \mathbb{E}_{ie} z_{{i,e}} ^{n+1} \label{EqLinearDiscretDAEWithDependency2}
\end{eqnarray}

Restricted Additive Schwarz iterative method \cite{tromeur_contrib_RAS} to solve Eq. \eqref{EqLinearDiscretDAEWithDependency2} consists to take the $z_{{i,e}}^{n+1,(k)}$ from previous iterate $(k)$ on the other parts to compute the $(k+1)$ iterate $z_{{i}}^{n+1,(k+1)}$ on the partition $W_i$. Starting from initial values $z_i^{n+1,(0)}$ it iterates:

\begin{eqnarray}
\mathbb{A}_{i} z^{n+1,(k+1)}_{i}&=&b_{i}^{n+1} - \mathbb{E}_{ie} z_{{i,e}} ^{n+1,(k)} \label{SchwarzLocalHomogene}
\end{eqnarray}

Let us rewrite the systems on a $W_i$ subdomain assuming that the values at the artificial interfaces retrieved by an EMT subdomain (respectively TS) necessarily come from a TS subdomain (respectively EMT). As the solution representation on the partition differ, the difficulty is to translate the exchanged quantities $z_{{i,e}} ^{n+1,(k)}$ between the TS and EMT and vice versa. These translations involve  the iterate solutions on several EMT time step and some combining of iterate solutions on TS part to define the EMT boundary  conditions over $m$ EMT time steps. 

\subsubsection{TS side}

In the dynamic Phasor case, the equations must first be adapted to the shape of the dynamic phasor by considering the differentiation property of dynamic phasors  and by multiplying the number of equations by the number $K$ of kept harmonics and by solving the real and imaginary parts separately.

Adapting the DAE system  Eq. \eqref{EqLinearDiscretDAEWithDependency} to the Dynamic Phasor shape solution leads to integrate on the TS part from $T^N$ to $T^{N+1}$ the DAE system as follows:
\begin{eqnarray}
\underbrace{\left( \begin{array}{cc}
\mathbb{I}_{i_{ts}} +\Delta t_{{ts}}  A_{i_{ts}}  &  \Delta t_{{ts}} B_{i_{ts}} \\
C_{i_{ts}}  & D_{i_{ts}}  \end{array} \right)}_{\mathbb{A}_{i_{ts}} } \underbrace{\left( \begin{array}{c}
x^{N+1}_{i_{ts}}  \\ y^{N+1}_{i_{ts}}  \end{array}\right)}_{w^{N+1}_{i_{ts}} } = \underbrace{\left( \begin{array}{cc}
\mathbb{I}_{{i}_{ts}}    &  0\\
0 & 0 \end{array} \right)}_{\mathbb{I}_{{d,i}_{ts}} }\underbrace{\left( \begin{array}{c}
x^{N}_{{i}_{ts}}  \\ y^{N}_{{i}_{ts}}  \end{array}\right)}_{w^{N}_{{i}_{ts}} } && \nonumber\\
 - \underbrace{\left( \begin{array}{cc}
\mathbb{I}_{{i,e}_{ts}}    &  0\\
0 & 0 \end{array} \right)}_{\mathbb{I}_{{d,ie}_{ts}} }\mathbb{T}^{emt}_{ts}\underbrace{\left(\begin{array}{c} X^{N}_{{i,e}_{emt}} \\ Y^{N}_{{i,e}_{emt}}\end{array}\right)}_{Z_{{i,e}_{emt}}^{N}} 
- \underbrace{\left( \begin{array}{cc}
\Delta t_{{ts}} E^{A}_{i_{ts}}  & \Delta t_{{ts}} E^{B}_{i_{ts}} \\
E^{C}_{i_{ts}}  & E^{D}_{i_{ts}} 
\end{array}\right)\mathbb{T}^{emt}_{ts}}_{\mathbb{E}^{emt}_{ie_{ts}} } 
\underbrace{\left(\begin{array}{c} X^{N+1}_{{i,e}_{emt}} \\ Y^{N+1}_{{i,e}_{emt}}\end{array}\right)}_{Z_{{i,e}_{emt}}^{N+1}} 
+ G^{N+1}_{i_{ts}} . \label{TSsubdomain}
\end{eqnarray}
The term $\mathbb{I}_{{d,ie}_{ts}} \mathbb{T}^{emt}_{ts}Z_{{i,e}_{emt}}^{N}$ is due to the fact that differential terms on the interface unknowns may be involved due to the splitting. The operator $\mathbb{T}^{emt}_{ts}$ is the translation operator from EMT to TS acting on the boundary values. By noting $\tilde{m}=\frac{T_0}{\Delta t_{emt}}$ the number of EMT time steps performed during a period, we can define $Z_{i,e_{emt}} \in \mathbb{R}^{(\tilde{m}\times n_{ie})\times 1}$ the vector gathering the $z_{ie} \in \mathbb{R}^{n_{ie}}$ over the $\tilde{m}$ EMT time steps. Then we have  $\mathbb{T}^{emt}_{ts}: \mathbb{R}^{ (\tilde{m}\times n_{ie}) \times 1} \longmapsto \mathbb{R}^{ (2K\times n_{ie}) \times 1} $. 
Let us notice that $Z_{{i,e}_{emt}}^{N}$ is a much larger vector than $\mathbb{T}^{emt}_{ts}Z_{{i,e}_ { emt} } ^{N}$. Moreover the $Z_{{i,e}_{emt}}^{N}$ have already been transformed into a TS form at the previous time step. Rather than keeping a large vector and recomputing the translation, we will keep $\mathbb{T}^{emt}_{ts}Z_{{i,e}_{emt}}^{N}$ calculated during the previous time step. Therefore, one can gather all the terms of the preceding time step in the same vector:
$w_{i,ie}^{N}=w_{i}^{N}+\mathbb{T}^{emt}_{ts}Z_{{ie}_{emt}}^{N}$, we rewrite the equation \eqref{TSsubdomain} with the appropriate definition of $\mathbb{I}_{{i,ie}_{ts}}$ for the gathering:
\begin{eqnarray}
\underbrace{\left( \begin{array}{cc}
\mathbb{I}_{i_{ts}} +\Delta t_{{ts}}  A_{i_{ts}}  & \Delta t_{{ts}} B_{i_{ts}} \\
C_{i_{ts}}  & D_{i_{ts}}  \end{array} \right)}_{\mathbb{A}_{i_{ts}} } \underbrace{\left( \begin{array}{c}
x^{N+1}_{i_{ts}}  \\ y^{N+1}_{i_{ts}}  \end{array}\right)}_{w^{N+1}_{i_{ts}} } &=& \underbrace{\left( \begin{array}{cc}
\mathbb{I}_{{i,ie}_{ts}}    &  0\\
0 & 0 \end{array} \right)}_{\mathbb{I}_{{d,i,ie}_{ts}} }\underbrace{\left( \begin{array}{c}
x^{N}_{{i,ie}_{ts}}  \\ y^{N}_{{i,ie}_{ts}}  \end{array}\right)}_{w^{N}_{{i,ie}_{ts}} }  \nonumber\\
- \underbrace{\left( \begin{array}{cc}
\Delta t_{{ts}} E^{A}_{i_{ts}}  &  \Delta t_{{ts}} E^{B}_{i_{ts}} \\
E^{C}_{i_{ts}}  & E^{D}_{i_{ts}} 
\end{array}\right)\mathbb{T}^{emt}_{ts}}_{\mathbb{E}^{emt}_{ie_{ts}} } 
\underbrace{\left(\begin{array}{c} X^{N+1}_{{i,e}_{emt}} \\ Y^{N+1}_{{i,e}_{emt}}\end{array}\right)}_{Z_{{i,e}_{emt}}^{N+1}} 
+ G^{N+1}_{i_{ts}} . \label{TSsubdomainFinal}
\end{eqnarray} 

\begin{remark}
The conditions imposed at the boundaries of a TS domain $W_{i,ts}^p$ can be considered as analogous to the Dirichlet conditions, in fact the differentiated part of the boundary conditions having passed into $w^{N,\ (\infty) }_{{i,ie}_{ts}}$  are fixed on the iterations of Schwarz and have no impact on the conditions at the interfaces.
\end{remark}

Let us detail the $\mathbb{T}^{emt}_{ts}$ that we developed. The three most common types of techniques used to translate EMT information into TS are: the curve fitting techniques as the least-sqares curve fitting technique used by Plumier\cite{PlumierPhD},  the methods of changing referential frames such as the $\alpha,\beta$ method used by Konara \cite{KonaraPhD} and by Zamroni \cite{ZamroniPhD}, and like the direct-quadrature-zero transformation method (dq0)\cite{PlumierPhD}, the 
methods based on the Fast Fourier Transform (FFT), as performed by Kumara in\cite{KumaraPhD}. Each of these techniques has these advantages and disadvantages, for the curve fitting technique it is the least precise of the 3 techniques, moreover it is not necessarily a linear operator (although often it is linear because it is often a linear interpolation). The methods of changing referential frame has the merit of being instantaneous, on the other hand it assumes that the three-phase current is necessarily balanced, which already induces there a loss of information, moreover we cannot hope to recover the slightest information on what has happened during the EMT intermediate time steps. This way of doing things will therefore be the best if what happens in the EMT part does not have too strong an impact on the TS part. The FFT is the most precise of the 3 transformations, and is linear, its disadvantage being that it is necessary to wait for a period to be able to apply it, which strongly restricts the choice of macro time steps. \\
 We slightly modified  FFT-based method, in order to get rid of the period restriction. To translate from the EMT to the dynamic phasor, we apply a fast Fourier transform on a history of the size of a period, $T$,  containing the values taken by the part calculated by the EMT at the intermediate time steps.\\
To perform an FFT on a signal, the history must cover a period of this signal. We can therefore take a macro time step of size $T_0$, the period and fill the history with the $m=\tilde{m}$ intermediate time steps. Since we decided that the macro time steps would be the TS time steps, we take $\Delta t _{ts}=T$ the size of a Period.

\begin{figure}[h!]
\centering
\begin{minipage}{13cm}
\begin{minipage}{6.2cm}
\includegraphics[scale=0.4]{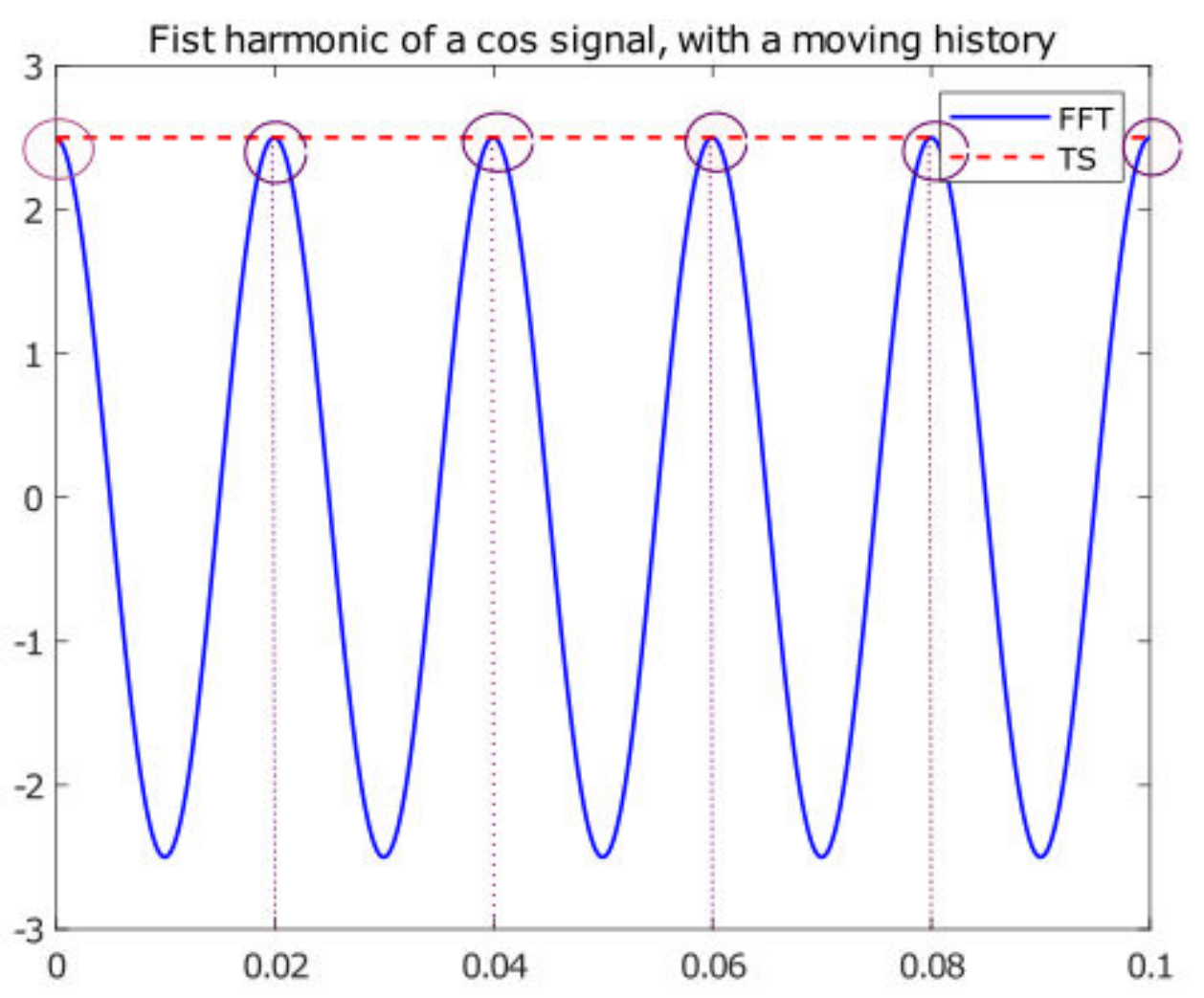}
\end{minipage}
\hfill
\begin{minipage}{6.2cm}
\includegraphics[scale=0.4]{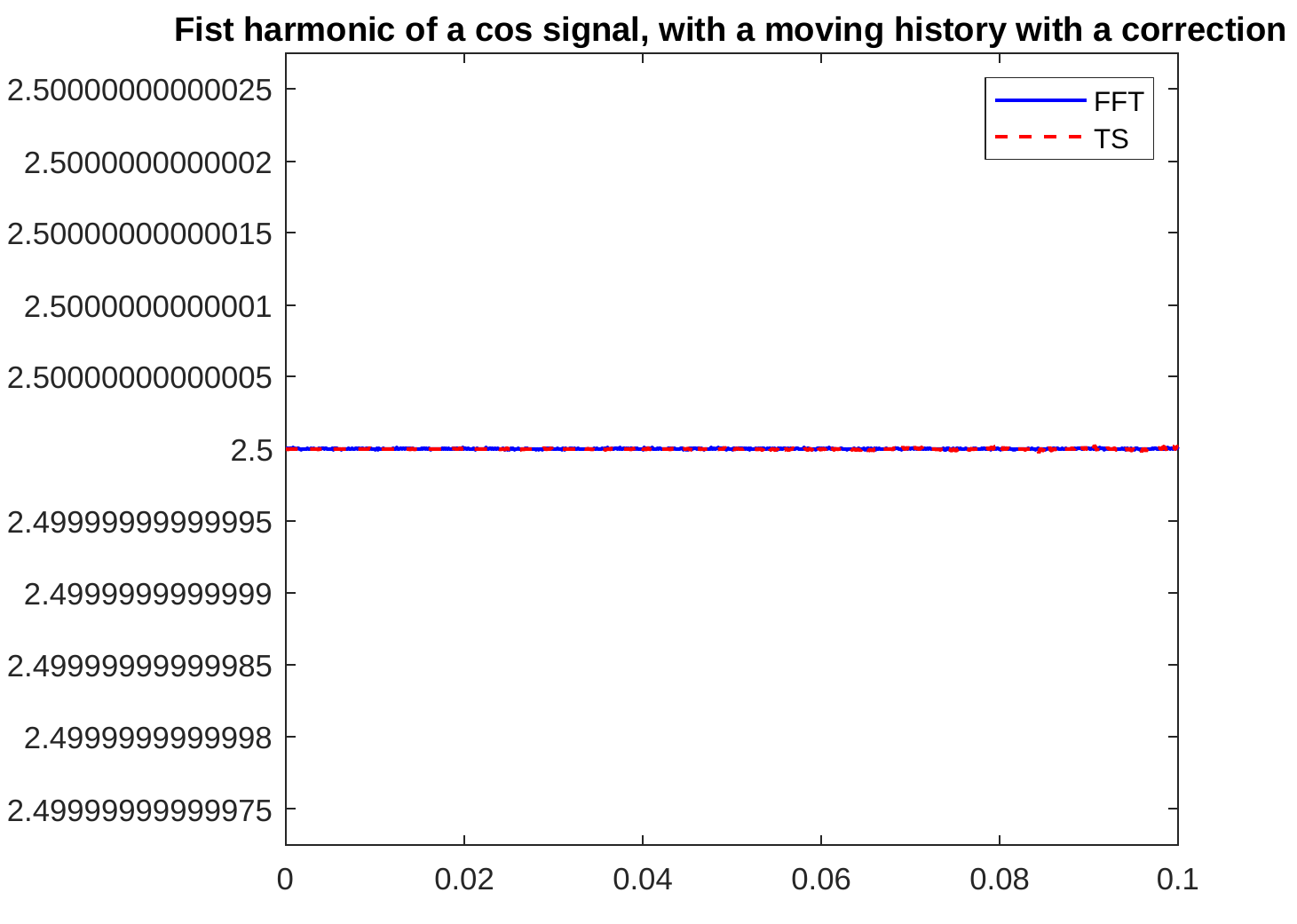}
\end{minipage}
\end{minipage}
\caption{$1^{st}$ Harmonic  of a moving history, computed with PhasorDynamic and an FFT. Left: without applying any correction, Right: with applying a correction \label{MovingFFT} }
\end{figure}

However, it is better to be able to take smaller TS time steps. We therefore take a $\Delta t _{ts}$ shorter than one period, however to perform the FFT, a history of the size of a period is always necessary, a “sliding” history will be produced. To do this, we will delete the $m$ (reminder $\Delta t _{ts} = m \Delta t _{emt} $) oldest values from the beginning of the history and put the $m$ new values at the end of the history. This operation is repeated at each instant of passage of information, i.e. after each TS time step (and therefore after $m$ time step EMT).\\
  \begin{definition}
   We recall that $m=\frac{\Delta t_{ts}}{\Delta t_{emt}}  $ and  $\tilde{m}=\frac{T}{\Delta t_{emt}}$, let's also denote $M=\frac{\tilde{m}}{m}$ the number of $\Delta t_{ts}$ during a period. 
  The history of the integration of $t^N$ to $t^{N+1}$, for the $k^{th}$ Schwarz iteration $Z_{{ie}_{emt}}^{N+1,(k)}$ is defined as follows:
  $$Z_{{ie}_{emt}}^{N+1,(k)}=[\underbrace{w_{{ie}_{emt}}^{m(N-M+1)+1,(\infty)},\hdots,w_{{ie}_{emt}}^{m(N-M+1)+m,(\infty)}, \hdots \hdots, z_{{ie}_{emt}}^{m\times N,(\infty)}}_{\text{Unactive boundary conditions}},\underbrace{w_{{ie}_{emt}}^{m\times N+1,(k)},\hdots,w_{{ie}_{emt}}^{m\times (N+1),(k)}}_{\mathbb{W}^p_{{ie}_{emt}} \text{Active Boundary conditions}}]^t$$
  \end{definition}
  \medskip 
However, an FFT on a moving story is also moving. See Figure \ref{MovingFFT} where we compare the harmonic $1^{st}$ of a cos signal calculated with a dynamic phasor and with an FFT (applied to the EMT signal), we can observe, on Figure left, that the FFT and the dynamic phasor meet only once per full period.
Indeed, for the Eq. \eqref{DynamicPhasor} the observation window will no longer be $[t-T,t]$ but $[t-T+\alpha  \Delta t_{ts} ,t +\alpha  \Delta t_{ts}]$, applying the appropriate change of variable in the Fourier transform formula, one can observe a $e^{jkw_{0}\alpha  \Delta t _{ts}}$ gap:
\begin{eqnarray}
{S_k}(t)= \frac{1}{T}\int_{t-T+\alpha  \Delta t_{ts} }^{t +\alpha  \Delta t_{ts}}s(\tau)e^{jkw_{0}\tau}d\tau =
 \frac{1}{T}\int_{t-T }^{t }s(X)e^{jkw_{0}(X+\alpha  \Delta t_{ts})}dX = e^{jkw_{0}\alpha  \Delta t_{ts}}\frac{1}{T}\int_{t-T }^{t }{s(X)e^{jkw_{0}X}}dX \nonumber
 \end{eqnarray}
We compensate this gap by applying $e^{-jkw_{0}\alpha  \Delta t _{ts}}$ to the FFT. Figure \ref{MovingFFT}(right)  shows the effect of the gap correction on the FFT that allows   the FFT of the EMT signal to correspond to the dynamic phasor.

\begin{definition}
The translation operator $\mathbb{T}^{emt}_{ts}$ consists of:
\begin{itemize}
\item Apply an FFT with a selection of the K modes corresponding to the K selected harmonics.
\item Next apply a compensation $e^{-jw_{0}\alpha  \Delta t _{ts}}$.
\item Finally separate the real and imaginary part of the obtained values.
\end{itemize}
\end{definition}

\begin{definition}
The $(k+1)$ iteration of the  Restrictive Additive Schwarz algorithm in the discrete case is written locally for the $W_{i,ts}^p$ partition from the TS type and with integrating between $T^N$ and $T^{N+1}$, if $ \mathbb{A}_{i_{ts}}$ invertible, as:
\begin{equation}
w^{N+1,(k+1)}_{i_{ts}} =\mathbb{A}_{i_{ts}}^{-1}(\mathbb{I}_{{i,ie}_{ts}}w_{i,ie}^{N,(\infty)}-\mathbb{E}^{emt}_{ie_{ts}}Z_{{i,e}_{emt}}^{N+1,(k)}+G^{N+1}_{i_{ts}})
\end{equation}\label{EqLinearDiscretDAERASHeterogeneTS}
\end{definition}

\subsubsection{EMT side}
For a $W_{i_{emt}}$  EMT partition, the discrete state space system Eq. \eqref{EqLinearDiscretDAEWithDependency} must be adapt as the data dependencies $z_{i,e}^{n+1}$ come from the TS partition. A translation operator from TS to EMT at time $t^q$: $\mathbb{T}^{ts}_{emt}(t^q):{ \mathbb{R}^{(2K\times n_{ie}}) \times 1} \longmapsto  \mathbb{R}^{ n_{ie} \times 1}$, where $K$ is the number of harmonics kept. As the time steps $\Delta t_{emt}$ are shorter than the macro time steps $\Delta t_{ts}$, the part due to the discretization of the interface values must therefore be calculated at each EMT time steps.

The integration from $t^n$ to $t^{n+1}=t^n+\Delta t_{emt}$ and $t^n,t^{n+1} \in [T^N,T^{N+1}]$ writes:
\begin{eqnarray}
\underbrace{\left( \begin{array}{cc}
\mathbb{I}_{i_{emt}}   +\Delta t_{emt} A_{i_{emt}} &  \Delta t_{emt} B_{i_{emt}}\\
C_{i_{emt}} & D_{i_{emt}} \end{array} \right)}_{\mathbb{A}_{i_{emt}}} \underbrace{\left( \begin{array}{c}
x^{n+1}_{i_{emt}} \\ y^{n+1}_{i_{emt}} \end{array}\right)}_{w_{i,emt}^{n+1}} =\underbrace{\left( \begin{array}{cc}
\mathbb{I}_{i_{emt}}   &  0\\
0 & 0 \end{array} \right)}_{\mathbb{I}_{d,i_{emt}}} \underbrace{\left( \begin{array}{c}
x^{n}_{i_{emt}}\\ y^{n}_{i_{emt}} \end{array}\right)}_{w_{i,emt}^{n}}+G^{n+1}_{i_{emt}} &&  \nonumber\\
 - \underbrace{\left( \begin{array}{cc}
\Delta t_{emt} E^{A}_{i_{emt}} &  \Delta t_{emt} E^{B}_{i_{emt}}\\
E^{C}_{i_{emt}} & E^{D}_{i_{emt}}
\end{array}\right)}_{E_{i,e_{emt}}^{ts}} \mathbb{T}^{ts}_{emt} (t^{n+1})\underbrace{\left(\begin{array}{c}
x^{N+1}_{{i,e}_{ts}}  \\ y^{N+1}_{{i,e}_{ts}}  \end{array}\right)}_{W_{{i,e}_{ts}} ^{N+1}}+ \underbrace{\left( \begin{array}{cc}
{\mathbb{I}}_{d,ie_{emt}} &  0\\
0&0
\end{array}\right)}_{\mathbb{I}_{ie_{emt}}}\mathbb{T}^{ts}_{emt} (t^{n})\ \underbrace{\left(\begin{array}{c}
x^{N+1}_{{i,e}_{ts}} \\ y^{N+1}_{{i,e}_{ts}} \end{array}\right)}_{W_{{i,e}_{ts}} ^{N+1}} . 
\end{eqnarray}
Where the operator $\mathbb{I}_{d,ie_{emt}}$ is $\mathbb{I}_{d,ie_{emt}}=R_{i_{emt}}^p \mathbb{I}_{d_{emt}}R_{{i,e}_{emt}}^{pT}$  and the term $ \mathbb{I}_{ie_{emt}}\mathbb{T}^{ts}_{emt} (t^{n})W_{{i,e}_{ts}} ^{N+1} $ is the part due to the fact that differential terms on the interface unknowns may be involved due to the splitting. \medskip \\
Let us consider the $m$  micro time steps realized by the EMT and by considering $t^n=T^N$, we can rewrite the behavior of the EMT on the whole time step of $T^{N}$ to $T ^{N+1}$ as:
{\small \begin{eqnarray}
\underbrace{\left( \begin{array}{cccc} 
 \mathbb{A}_{i_{emt}}&&&\\
 -\mathbb{I}_{i_{emt}}&  \mathbb{A}_{i_{emt}} &&\\
 &\ddots & \ddots  & \\
&& -\mathbb{I}_{i_{emt}} & \mathbb{A}_{i_{emt}}
\end{array}\right)}_{\mathbb{H}_{i_{emt}}}\underbrace{\left( \begin{array}{c}
w^{n+1}_{i_{emt}} \\ w^{n+2}_{i_{emt}}  \\ \vdots \\  w^{n+m}_{i_{emt}} 
 \end{array}\right)}_{\mathbb{W}_{i_{emt}}^{N+1} }=\underbrace{\left( \begin{array}{c}
 G^{n+1}_{i_{emt}} +\mathbb{I}_{i_{emt}} w^n_{i_{emt}}+ \mathbb{I}_{d,ie_{emt}} (x^n,y^n)^t\\ G^{n+2} \\ \vdots \\ _{i_{emt}} \\ G^{n+m}_{i_{emt}}
 \end{array}\right)}_{\mathbb{G}_{i_{emt}}^{N+1}}&& \nonumber \\
 +\underbrace{\left( \begin{array}{cccc} 
 E_{i,e_{emt}}^{ts} &&&\\
 -\mathbb{I}_{d,ie_{emt}}  &   E_{i,e_{emt}}^{ts}&\\
   & \ddots  &\ddots & \\
&&-\mathbb{I}_{d,ie_{emt}} &  E_{i,e_{emt}}^{ts}  
\end{array}\right)}_{\mathbb{E}^{ts}_{i,e_{emt}}}\underbrace{\left( \begin{array}{c}
  W_{{i,e}_{ts}}^{N+1}(t^{n+1})\\ W_{{i,e}_{ts}}^{N+1}(t^{n+2})  \\\vdots \\   W_{{i,e}_{ts}}^{N+1}(t^{n+m}) 
 \end{array}\right)}_{\mathbb{W}^{N+1}_{{i,e}_{ts}}} 
.&&\label{DAERASemt}
\end{eqnarray}}
With $W_{{i,e}_{ts}}^{N+1}(t^{n+1})=\mathbb{T}^{ts}_{emt} (t^{n+1})W_{{i,e}_{ts}} ^{N+1}$

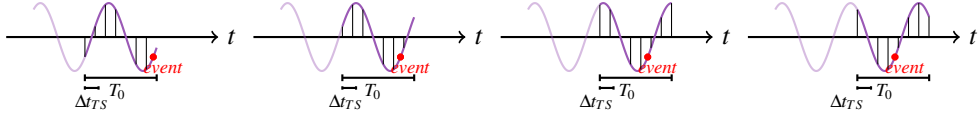
\begin{figure}[h!]

		\begin{tikzpicture}[scale=0.45]
	\tikzmath{ \fo =10; \A = 20.;\To = 1/\fo;\nb = 20; \Te = \To/\nb;};
	\draw[->, thick](-0.15*\A, 0)--(0.16*\A, 0) node[anchor = west]{$t$};
	\draw[samples=300, domain=-0.11:-0.035, color = vc,opacity=0.4, thick, scale = \A] plot ({\x},{sin(deg(\x*2*3.141*\fo)+ 90)/\A});
	\draw[samples=300, domain=-0.035:0.07, color = vc, thick, scale = \A] plot ({\x},{sin(deg(\x*2*3.141*\fo)+ 90)/\A});
	 \draw[thick] (-0.035*\A,-1.2) -- (0.07*\A,-1.2);
	 \draw[thick] (-0.035*\A,-1.3) -- (-0.035*\A,-1.1);
	 \draw[thick] (0.07*\A,-1.3) -- (0.07*\A,-1.1);
	  \draw (0.25,-1.1)  node[below]{\scriptsize $T_0$};
	  
	  	  \draw[thick] (-0.035*\A,-1.5) -- (-0.015*\A,-1.5);
	 \draw[thick] (-0.035*\A,-1.45) -- (-0.035*\A,-1.55);
	 \draw[thick] (-0.015*\A,-1.45) -- (-0.015*\A,-1.55);
	  \draw (-0.5,-1.4)  node[below]{\scriptsize $\Delta t _{TS}$};
\draw[samples=100, domain=0.065:0.065, color = red, thick, scale = \A, mark = *, mark size=0.1pt,only marks] plot ({\x},{sin(deg(\x*2*3.141*\fo)+ 90)/\A});
	 \draw (0.08*20,-0.5)  node[below,color=red]{\scriptsize $event$};
		\foreach \x in {-7,-4,...,13} {
		\tikzmath{ \y=\x*\Te; \ang = sin(deg(\y*2*3.141*\fo) + 90); };		
		\draw( {\y*\A} , {\ang} ) -- ( {\y*\A} ,0) ;
	}
	\end{tikzpicture}
		\begin{tikzpicture}[scale=0.45]
	\tikzmath{ \fo =10; \A = 20.;\To = 1/\fo;\nb = 20; \Te = \To/\nb;};
	\draw[->, thick](-0.15*\A, 0)--(0.16*\A, 0) node[anchor = west]{$t$};
	\draw[samples=300, domain=-0.11:-0.02, color = vc,opacity=0.4, thick, scale = \A] plot ({\x},{sin(deg(\x*2*3.141*\fo)+ 90)/\A});
	\draw[samples=300, domain=-0.02:0.085, color = vc, thick, scale = \A] plot ({\x},{sin(deg(\x*2*3.141*\fo)+ 90)/\A});
	 \draw[thick] (-0.02*\A,-1.2) -- (0.085*\A,-1.2);
	 \draw[thick] (-0.02*\A,-1.3) -- (-0.02*\A,-1.1);
	 \draw[thick] (0.085*\A,-1.3) -- (0.085*\A,-1.1);
	  \draw (0.5,-1.1)  node[below]{\scriptsize $T_0$};
	  	  	  \draw[thick] (-0.02*\A,-1.5) -- (0*\A,-1.5);
	 \draw[thick] (-0.02*\A,-1.45) -- (-0.02*\A,-1.55);
	 \draw[thick] (0*\A,-1.45) -- (0*\A,-1.55);
	  \draw (-0.3,-1.4)  node[below]{\scriptsize $\Delta t _{TS}$};
\draw[samples=100, domain=0.065:0.065, color = red, thick, scale = \A, mark = *, mark size=0.1pt,only marks] plot ({\x},{sin(deg(\x*2*3.141*\fo)+ 90)/\A});
	 \draw (0.08*20,-0.5)  node[below,color=red]{\scriptsize $event$};
		\foreach \x in {-4,-1,...,16} {
		\tikzmath{ \y=\x*\Te; \ang = sin(deg(\y*2*3.141*\fo) + 90); };		
		\draw( {\y*\A} , {\ang} ) -- ( {\y*\A} ,0) ;
	}
	\end{tikzpicture}
			\begin{tikzpicture}[scale=0.45]
	\tikzmath{ \fo =10; \A = 20.;\To = 1/\fo;\nb = 20; \Te = \To/\nb;};
	\draw[->, thick](-0.15*\A, 0)--(0.16*\A, 0) node[anchor = west]{$t$};
	\draw[samples=300, domain=-0.11:-0.005, color = vc,opacity=0.4, thick, scale = \A] plot ({\x},{sin(deg(\x*2*3.141*\fo)+ 90)/\A});
	\draw[samples=300, domain=-0.005:0.1, color = vc, thick, scale = \A] plot ({\x},{sin(deg(\x*2*3.141*\fo)+ 90)/\A});
	 \draw[thick] (-0.005*\A,-1.2) -- (0.1*\A,-1.2);
	 \draw[thick] (-0.005*\A,-1.3) -- (-0.005*\A,-1.1);
	 \draw[thick] (0.1*\A,-1.3) -- (0.1*\A,-1.1);
	  \draw (0.7,-1.1)  node[below]{\scriptsize $T_0$};
	  \draw[thick] (-0.005*\A,-1.5) -- (0.015*\A,-1.5);
	  	 \draw[thick] (-0.005*\A,-1.45) -- (-0.005*\A,-1.55);
	 \draw[thick] (0.015*\A,-1.45) -- (0.015*\A,-1.55);
	  \draw (-0.2,-1.4)  node[below]{\scriptsize $\Delta t _{TS}$};
\draw[samples=100, domain=0.065:0.065, color = red, thick, scale = \A, mark = *, mark size=0.1pt,only marks] plot ({\x},{sin(deg(\x*2*3.141*\fo)+ 90)/\A});
	 \draw (0.08*20,-0.5)  node[below,color=red]{\scriptsize $event$};
		\foreach \x in {-1,2,...,20} {
		\tikzmath{ \y=\x*\Te; \ang = sin(deg(\y*2*3.141*\fo) + 90); };		
		\draw( {\y*\A} , {\ang} ) -- ( {\y*\A} ,0) ;
	}
	\end{tikzpicture}
				\begin{tikzpicture}[scale=0.45]
	\tikzmath{ \fo =10; \A = 20.;\To = 1/\fo;\nb = 20; \Te = \To/\nb;};
	\draw[->, thick](-0.15*\A, 0)--(0.16*\A, 0) node[anchor = west]{$t$};
	\draw[samples=300, domain=-0.11:0.01, color = vc,opacity=0.4, thick, scale = \A] plot ({\x},{sin(deg(\x*2*3.141*\fo)+ 90)/\A});
	\draw[samples=300, domain=0.01:0.115, color = vc, thick, scale = \A] plot ({\x},{sin(deg(\x*2*3.141*\fo)+ 90)/\A});
	 \draw[thick] (0.01*\A,-1.2) -- (0.115*\A,-1.2);
	 \draw[thick] (0.01*\A,-1.3) -- (0.01*\A,-1.1);
	 \draw[thick] (0.115*\A,-1.3) -- (0.115*\A,-1.1);
	  \draw (1,-1.1)  node[below]{\scriptsize $T_0$};
	  	  	  	  \draw[thick] (0.01*\A,-1.5) -- (0.03*\A,-1.5);
	 \draw[thick] (0.01*\A,-1.45) -- (0.01*\A,-1.55);
	 \draw[thick] (0.03*\A,-1.45) -- (0.03*\A,-1.55);
	  \draw (-0.1,-1.4)  node[below]{\scriptsize $\Delta t _{TS}$};
\draw[samples=100, domain=0.065:0.065, color = red, thick, scale = \A, mark = *, mark size=0.1pt,only marks] plot ({\x},{sin(deg(\x*2*3.141*\fo)+ 90)/\A});
	 \draw (0.08*20,-0.5)  node[below,color=red]{\scriptsize $event$};
		\foreach \x in {2,5,...,23} {
		\tikzmath{ \y=\x*\Te; \ang = sin(deg(\y*2*3.141*\fo) + 90); };		
		\draw( {\y*\A} , {\ang} ) -- ( {\y*\A} ,0) ;
	}
	\end{tikzpicture}
\caption{Impact of an EMT event in the History \label{Hist}}
\end{figure}

\begin{figure}[h]
\centering
\includegraphics[scale=0.5]{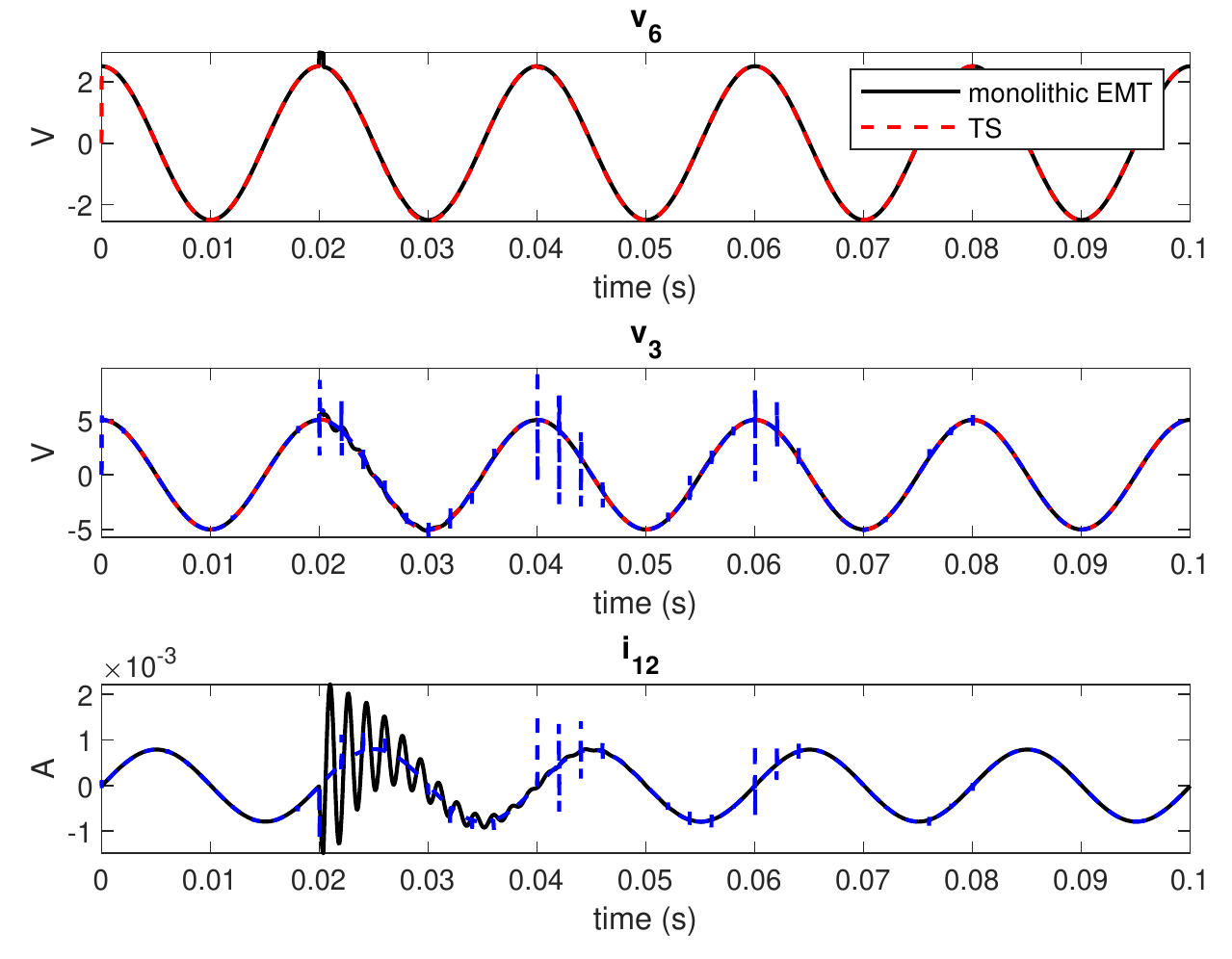}
\hfill
\includegraphics[scale=0.5]{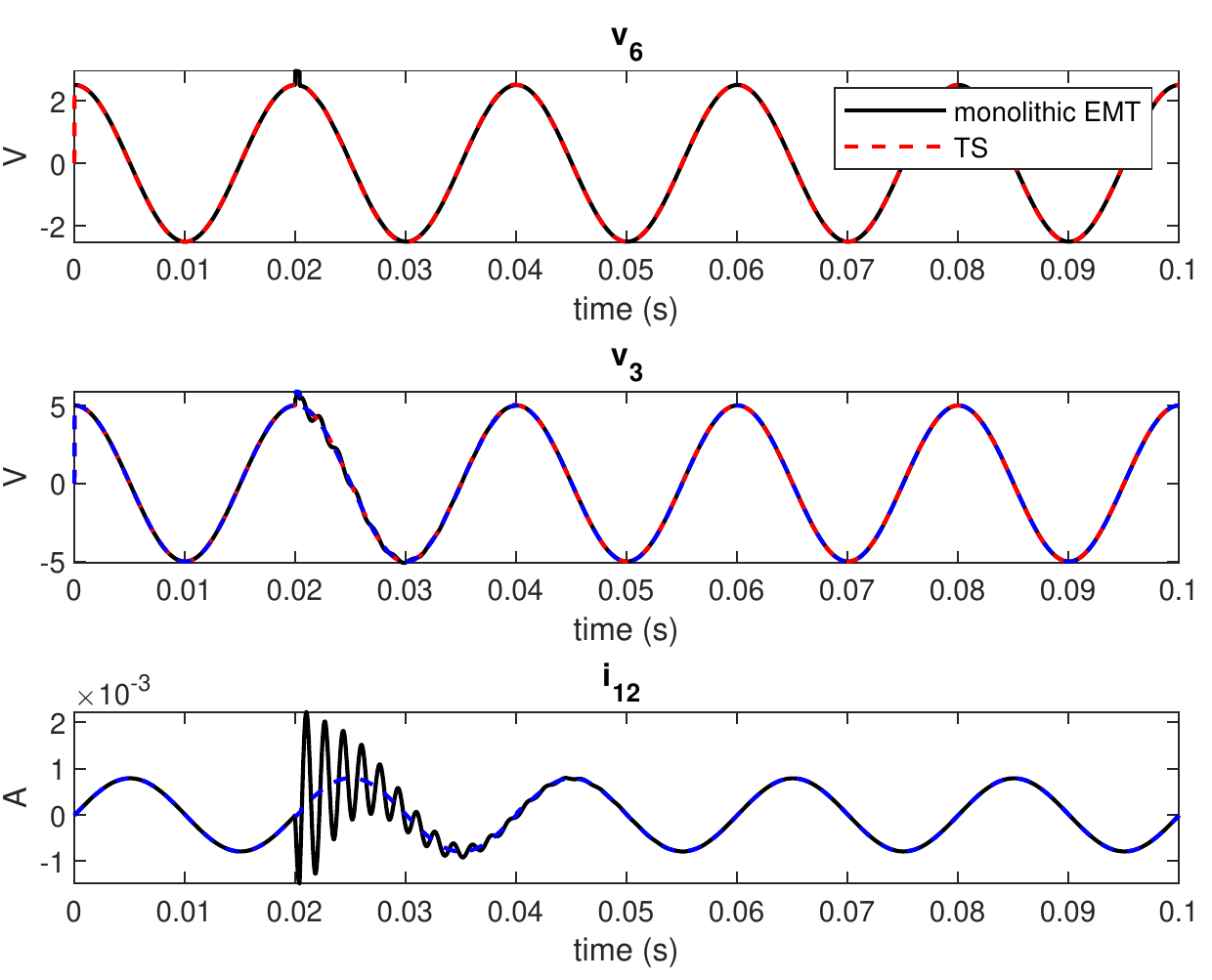}
\caption{ TS-EMT hybrid simulation  with a jump in the voltage source at $t=0.2s$ that last less than $\Delta t_{TS}$:   without (top) and with (bottom) linear interpolation of  $\mathbb{W}^{N+1}_{TS}$ \label{noise}}
\end{figure}

Let us detail the $\mathbb{T}^{ts}_{emt} (t^{n})$ translation operator .
To rebuild an EMT signal from the coefficients of the dynamic phasor, at each EMT time step we use the recombination of the $K \in I$ harmonics kept. However, when translating from EMT to dynamic phasor, a rolling history is used. With this moving history, an event has repercussions long after it ends. Indeed, as the figure (\ref{Hist}) shows, an EMT event may still be in the history long after the event has ended, and may therefore still have an effect on the translation of EMT to TS.
Due to the jump between the first two rows of $\mathbb{W}^{N+1}_{{i,e}_{ts}}$ in the Eq. \eqref{DAERASemt}, and the repercussions of an event remained in the history of the EMT, there are leaps in the solution on the side of the EMT. Indeed, these repercussions have almost no impact on the dynamic phasor part, but they have a significant impact on the EMT part. The fact that it is visible on the EMT part is due to the Gibbs phenomenon. Indeed, when there is a discontinuity jump, the Fourier sums (on which the translation of TS into EMT is based) reveal the Gibbs phenomenon which is expressed in the form of an overshoot at this jump then of swings. Gottlieb and Shu \cite{Gibbs} provide insight into this phenomenon as well as a way to obtain accurate function approximations, despite the presence of Gibbs' phenomenon, using filters. We made the choice to smooth the first jump (the one between $(x^n,y^n)^t$ and $W_{{i,e}_{ts}}^{N+1}(t^{n+1})$)  which is not a physical jump but a numerical jump, in order not to have the appearance of this phenomenon. \\
Figure \ref{noise}(left) shows these spurious oscillations  phenomena in EMT after an  event. Attenuations of these oscillations' amplitude  at each period  occur until the event is out of the history.
To smooth out these phenomena ,  the TS signal recombining in the EMT part is applied with a linear interpolation between two successive time steps of TS. $\alpha \in [1,m]$ is the number of EMT time steps over which the TS values are smoothed, where  $m$  is the number of EMT time steps during a TS time step, and $n $ is such that  $T^N=t^n$. Figure \ref{noise}(right) shows that the Gibbs phenomena disappears with the proposed smoothing.

Then, we can define the $~^{m+n}_{n+1}\mathbb{T}^{ts,\alpha}_{emt}$ translation operator:
\begin{definition}
We introduce  the operator $[~^{m+n}_{n+1}\mathbb{T}^{ts,\alpha,N}_{emt}]$ the translation operator from TS to EMT acting on the $m$ EMT time step starting from $t_{n+1}$ to $t_{n+m}$ as:
\begin{eqnarray}
[~^{m+n}_{n+1}\mathbb{T}^{ts,\alpha,N}_{emt}] &=& (\mathbb{T}^{ts,\alpha,N}_{emt} (t^{n+1}),\ldots, \mathbb{T}^{ts,\alpha,N}_{emt} (t^{n+m}))^T
\end{eqnarray}
where $\mathbb{T}^{ts,\alpha,N}_{emt} (t^{l})$ is defined as:
\begin{eqnarray}
\mathbb{T}^{ts,\alpha,N}_{emt} (t^{l})W_{{i,e}_{ts}} ^{N+1} =
\left\{
\begin{array}{ll}
  & \sum_{k \in I} (\frac{n+\alpha-l }{\alpha-1} w^{N}_{k;{i,e}_{ts}}+ \frac{n+1-l}{1-\alpha}w^{N+1}_{k;{i,e}_{ts}} ) e^{jk\omega_{0}t^{l}}.\text{ if }  l \in [n+1,n+\alpha] \\
& \sum_{k \in I} w^{N+1}_{k;{i,e}_{ts}}  e^{jk\omega_{0}t^{l}}.\text{ if }  l \in [n+\alpha,n+m] 
\end{array} \right.
\end{eqnarray}
\end{definition}

\begin{definition}
The $(k+1)$ iteration of the  Restrictive Additive Schwarz algorithm in the discrete case is written localy for the $W_{i,ts}^p$ partition from the EMT type and with integrating between $T^N$ and $T^{N+1}$, if $\mathbb{H}_{i_{emt}}$ inversible, as:
\begin{equation}
\mathbb{W}_{i_{emt}}^{N+1,(k+1)} =\mathbb{H}_{i_{emt}}^{-1}(\mathbb{E}^{ts}_{i,e_{emt}} [~^{m+n}_{n+1}\mathbb{T}^{ts,\alpha,N}_{emt}] W_{{i,e}_{ts}} ^{N+1,(k)}+\mathbb{G}_{i_{emt}}^{N+1})
\end{equation}\label{EqLinearDiscretDAERASHeterogeneEMT}
\end{definition}
\subsection{Initialization}

Correct initialization is important in our method, already because if we start from a solution that is too far from the true solution, more iterations will be needed to obtain the true solution (if we do not use the acceleration method of the Aitken convergence). Bad initial conditions can also be considered as an event, an undesirable physical disturbance.\\
The difference with the usual methods is that in this heterogeneous method, the initial conditions do not impact only the terms due to the discretization. Indeed, this correct initialization is particularly important in our heterogeneous TS-EMT Schwarz because of the translation operators used.
\begin{itemize}
\item For the EMT to TS translation operator: in order to be able to start the translation from the EMT to TS, a history of the size of an already pre-filled period is necessary. It is therefore necessary to initialize over a period in time step EMT. With the moving history, the initial conditions can be present in the history on several macro time steps, thus these initial conditions have a diffuse impact on several time steps and not only on the initial moment of the simulation.
\item For the operator of translation of TS towards EMT: A linear interpolation is used between two successive macro time steps, for the first time step the impact of the initial condition is thus very strong.
\end{itemize}
However, it is recalled that one of the motivations of this co-simulation is to gain in computational efficiency, which is why it would make no sense to have an initialization that is costly in terms of computation time. To meet these requirements, two choices are available to us: the first is to calculate the monolithic steady-state simulation over a full period with a time step equal to the EMT time step.
One recovers the real part of the solution in steady state at each step of time to fill the history. The second possibility is to calculate the monolithic steady-state simulation over a complete period with a time step equal to the time step TS and to apply a linear extrapolation for each time step EMT to fill the history. \\

For a TS domain $W_{i,ts} $ the EMT  history initialization condition is given by: \\
$ Z_{{0,ie}_{emt}}=[\Re(z_{ie_{steadystate}}(0)), \hdots, \Re(z_{ie_{steadystate}}(\tilde{m}\Delta t_{emt} ))]^t$. \\
To be coherent the TS domain $W_{i,ts}$ is initialize as  $w_{0,i_{ts}} = \mathbb{T}^{emt}_{ts} Z_{{0,i}_{emt}}$\\

For an EMT subdomain $W_{d,emt} $, only the first value of $\mathbb{W}_{d_{emt}}^{1}$ needs to be initialized, and this initialization is given by: $w_{0,d_{emt}}= \Re(z_{d_{steadystate}}(\tilde{m}\Delta t_{emt} )) $. We also initialize $\mathbb{W}_{{de}_{ts}} $:  $ W_{{0,de}_{ts}}(t^l)= \sum_{k \in I} w_{0,k;{d,e}_{ts}}  e^{jk\omega_{0}t^{l}}$ for $ l \in [0,m] $

\subsection{Heterogeneous EMT-TS RAS formulation}
Finally, we have a local writing of the Schwarz algorithm for the different EMT and TS subsystems. We also have the initialization of the problem. We can therefore write the RAS  for integrating the DAE system Eq. \eqref{DAEmono1} from $T^N$ to $T^{N+1}$:
\begin{definition}
The $(k+1)^{th}$ iteration of the RAS can be written for two subdomains  as follows:
\begin{eqnarray}
\left\{ \begin{array}{lcl}
w^{N+1,(k+1)}_{1_{ts}} &=&\mathbb{A}_{1_{ts}}^{-1}(\mathbb{I}_{{1;1e}_{ts}}w_{1,1e}^{N,(\infty)}-\mathbb{E}^{emt}_{1,e_{ts}}Z_{{1e}_{emt}}^{N+1,(k)}+G^{N+1}_{1_{ts}}),\\
\mathbb{W}_{2_{emt}}^{N+1,(k+1)}& =&\mathbb{H}_{2_{emt}}^{-1}(\mathbb{E}^{ts}_{2,e_{emt}}[~^{m+n}_{n+1}\mathbb{T}^{ts,\alpha,N}_{emt}] W_{{2,e}_{ts}} ^{N+1,(k)}+\mathbb{G}_{2_{emt}}^{N+1})\end{array}  \right.  \label{SchwarzHetero}
\end{eqnarray}
\end{definition}

\section{Heterogeneous EMT-TS RAS error operator and Acceleration of the convergence \label{shourick_contrib_Sec4}}
The system  \eqref{SchwarzHetero} is iterated at each macro-time step until the error between two successive iterations reaches a set tolerance. In this section we  discuss the impact of the heterogeneity of the co-simulation on the convergence of the method and on the acceleration of the convergence strategy.

\subsection{Error operator}
The error operator needs additional definition due to the translation operators. 

Let's calculate the error on the TS side to see the impact of rolling history on the error:
\begin{eqnarray}
w^{N+1,(k+1)}_{1_{ts}} - w^{N+1,(\infty)}_{1_{ts}}&=&\mathbb{A}_{1_{ts}}^{-1}(-\mathbb{E}^{emt}_{1,e_{ts}}(Z_{{1e}_{emt}}^{N+1,(k)}-Z_{{1e}_{emt}}^{N+1,(\infty)})),\nonumber \\
w^{N+1,(k+1)}_{1_{ts}} - w^{N+1,(\infty)}_{1_{ts}}&=&\mathbb{A}_{1_{ts}}^{-1}(-\mathbb{E}^{emt}_{1,e_{ts}}...\nonumber\\
&& [0, \hdots ,0,\underbrace{w_{{1e}_{emt}}^{m\times N+1,(k)}-w_{{1e}_{emt}}^{m\times N+1,(\infty)},\hdots,w_{{1e}_{emt}}^{m\times (N+1),(k)}-w_{{1e}_{emt}}^{m\times (N+1),(\infty)}}_{\mathbb{W}^{N+1,(k)}_{{1e}_{emt}} -\mathbb{W}^{N+1,(\infty)}_{{1e}_{emt}}}]^T),\nonumber
\end{eqnarray}

\begin{definition}
Let $\tilde{R}_{hist}^{0T}$ be the prolongation operator on the history on the EMT side that takes the $m$ active boundary conditions for the TS side and complete the unactive boundary conditions history by $0$.
\end{definition}
Then we have:
\begin{eqnarray}
Z_{{1e}_{emt}}^{N+1,(k)}-Z_{{1e}_{emt}}^{N+1,(\infty)}&=& \tilde{R}_{hist}^{0T} (\mathbb{W}^{N+1,(k)}_{{1e}_{emt}} -\mathbb{W}^{N+1,(\infty)}_{{1e}_{emt}})\nonumber \\
&=& \tilde{R}_{hist}^{0T} (I_m  \otimes R_{1,e_{ts}}^p) (\mathbb{W}^{N+1,(k)}_{{2}_{emt}} -\mathbb{W}^{N+1,(\infty)}_{{2}_{emt}}) \label{ErrorOp1}
\end{eqnarray}

When recombining harmonics, there is no change in the level of information, so the recombination itself will not change the error. On the other hand, smoothing can have an impact. Let's calculate the EMT error to see the impact of smoothing on the error:

For $ l \in [n+1,n+\alpha]$ we can write:
\begin{eqnarray}
W_{{2,e}_{ts}} ^{N+1,(k)}(t^l)-W_{{2,e}_{ts}} ^{N+1,(\infty)}(t^l)&=& \sum_{i \in I} (\frac{n+\alpha-l }{\alpha-1} w^{N}_{i;{2,e}_{ts}}+ \frac{n+1-l}{1-\alpha}w^{N+1,(k)}_{i;{2,e}_{ts}} ) e^{j i \omega_{0}t^{l}} \nonumber \\
&&-\sum_{i \in I}( \frac{n+\alpha-l }{\alpha-1} w^{N}_{i;{2,e}_{ts}}+ \frac{n+1-l}{1-\alpha}w^{N+1,(\infty)}_{i;{2,e}_{ts}} ) e^{ji\omega_{0}t^{l}} \nonumber \\
W_{{2,e}_{ts}} ^{N+1,(k)}(t^l)-W_{{2,e}_{ts}} ^{N+1,(\infty)}(t^l)&=&\sum_{i \in I} \frac{n+1-l}{1-\alpha}(w^{N+1,(k)}_{i;{2,e}_{ts}} - w^{N+1,(\infty)}_{i;{2,e}_{ts}}) e^{ji\omega_{0}t^{l}}.\nonumber
\end{eqnarray}  

Therefore we have:
 \begin{eqnarray}
\mathbb{W}_{2_{emt}}^{N+1,(k+1)}-\mathbb{W}_{2_{emt}}^{N+1,(\infty)}& =&\mathbb{H}_{2_{emt}}^{-1}\mathbb{E}^{ts}_{2,e_{emt}}  [~^{m+n}_{n+1}\mathbb{T}^{ts,\alpha,N}_{emt}] ({W}^{N+1,(k)}_{{2,e}_{ts}}-{W}^{N+1,(\infty)}_{{2,e}_{ts}}).\nonumber\\
& =&\mathbb{H}_{2_{emt}}^{-1}\mathbb{E}^{ts}_{2,e_{emt}}  [~^{m+n}_{n+1}\mathbb{T}^{ts,\alpha,N}_{emt}] R_{2,e}^p ({w}^{N+1,(k)}_{{1}_{ts}}-{w}^{N+1,(\infty)}_{{1}_{ts}}). \label{ErrorOp2}
\end{eqnarray}
So there is a smoothing of $ \frac{n+1-l}{1-\alpha} \in [0,1] $ on the first $\alpha$ values of each macro step, we expect a behavior of reduction of the radius of convergence, i.e. the method will converge less quickly or will diverge less quickly. It is expected that the translation operator from TS to EMT has no impact on the linearity of the error, indeed this operator being a combination of two linear operators is linear.

Finally, by using Eq. \eqref{ErrorOp1} and \eqref{ErrorOp2}, we can express the error operator $P$ of the heterogeneous RAS EMT-TS method with two partitions:
\begin{prop}
The error operator $P$ of the heterogeneous RAS EMT-TS method writes:
\begin{eqnarray}
P&=& \left(\begin{array}{cc} 
0 & \mathbb{A}_{1_{ts}}^{-1}(-\mathbb{E}^{emt}_{1,e_{ts}} \tilde{R}_{hist}^{0T} (I_m  \otimes R_{1,e_{ts}}^p) \\
\mathbb{H}_{2_{emt}}^{-1}\mathbb{E}^{ts}_{2,e_{emt}}  [~^{m+n}_{n+1}\mathbb{T}^{ts,\alpha,N}_{emt}] R_{2,e}^p & 0
\end{array} \right). \label{SchwarzHeteroOperator}
\end{eqnarray}
We can writes the expression of the error between two consecutive RAS iterates:

\begin{eqnarray}
\left( \begin{array}{l}
w^{N+1,(k+1)}_{1_{ts}} -w^{N+1,(k)}_{1_{ts}} \\ \mathbb{W}_{2_{emt}}^{N+1,(k+1)}-\mathbb{W}_{2_{emt}}^{N+1,(k)} \end{array} \right)   &=&  \underbrace{\left( \begin{array}{cc} 0 & P_{ts} \\ P_{emt} & 0 \end{array} \right)}_{P} \left( \begin{array}{l}
w^{N+1,(k)}_{1_{ts}} -w^{N+1,(k-1)}_{1_{ts}} \\ \mathbb{W}_{2_{emt}}^{N+1,(k)}-\mathbb{W}_{2_{emt}}^{N+1,(k-1)} \end{array} \right) \label{RASerror}
\end{eqnarray}
\end{prop}

Equation Eq. \eqref{RASerror} works on the global solution belonging to $W_1^p \times \underbrace{W_2^p \times \ldots \times W_2^p}_{m \textrm{~times}}$. \\
Nonetheless, we can reduce the size of the operator error by working only on the exchanged values,  due to the restriction operator property $R_i^p\tilde{\mathbb{A}}=R_i^p\tilde{\mathbb{A}}(R_i^{pT}R_i^p  +R_{i,e}^{pT} R_{i,e}^p)$.

 \begin{definition}[The RAS  global interface and its restriction operator]
 We define the RAS global interface $\Gamma=\left\{ W_{1,e}^p,\underbrace{W_{2,e}^p,\ldots,W_{2,e}^p}_{m \textrm{~times}} \right\}$ of size  $n_\Gamma$. We define the restriction operator $R_\Gamma$ associated to the RAS global interface as  $R_\Gamma =  \left( \begin{array}{cc} R_{1,e}^p& 0 \\ 0 & I_m \otimes R_{2,e}^p \end{array} \right) \in \mathbb{R}^{n_\Gamma \times (n1+m \, n_2)}$.
 \end{definition}
\begin{prop}

 The error operator $P_\Gamma$ acting on the interface value of the RAS iteration writes $P_\Gamma=R_\Gamma P R_\Gamma^T$ and we have:
 \begin{eqnarray}
\left( \begin{array}{l}
w^{N+1,(k+1)}_{1,e_{ts}} -w^{N+1,(k)}_{1,e_{ts}} \\ \mathbb{W}_{2,e_{emt}}^{N+1,(k+1)}-\mathbb{W}_{2,e_{emt}}^{N+1,(k)} \end{array} \right)   &=&  P_\Gamma \left( \begin{array}{l}
w^{N+1,(k)}_{1,e_{ts}} -w^{N+1,(k-1)}_{1,e_{ts}} \\ \mathbb{W}_{2,e_{emt}}^{N+1,(k)}-\mathbb{W}_{2,e_{emt}}^{N+1,(k-1)} \end{array} \right) \label{RASerrorGamma}
\end{eqnarray}
\end{prop}

\subsection{Aitken's acceleration of the RAS convergence}
The error operators $P$  Eq. \eqref{SchwarzHeteroOperator} and  $P_\Gamma$ do not depend of the RAS iterate $(k)$ for linear electical network. Then we have the following proposition to obtain the true solution on the partitions interfaces by using the Aitken's acceleration of the convergence technique and the Eq. \eqref{RASerrorGamma}:

\begin{prop}
 If $1$ is not a eigenvalue of the operator $P_\Gamma$ is then possible to accelerate the RAS to the true solution  with the Aitken's acceleration of the convergence technique as follows:
\begin{eqnarray}
\left( \begin{array}{l}
w^{N+1,(\infty)}_{1,e_{ts}}  \\ \mathbb{W}_{2,e_{emt}}^{N+1,(\infty)} \end{array} \right)   &=&  (I_{n_\Gamma}-P_\Gamma)^{-1} \left(\left( \begin{array}{l}
w^{N+1,(k)}_{1,e_{ts}} \\ \mathbb{W}_{2,e_{emt}}^{N+1,(k)}\end{array} \right)-P_\Gamma \left( \begin{array}{l}
w^{N+1,(k-1)}_{1,e_{ts}} \\ \mathbb{W}_{2,e_{emt}}^{N+1,(k-1)}\end{array} \right)\right), \, k\geq 1 \label{Aitken}
\end{eqnarray}
\end{prop}

{\bf \em This RAS acceleration technique is the essential property of the proposed co-simulation algorithm, as we will see in the numerical results, the convergence or divergence of the RAS depends strongly on the power system components, the time steps chosen for TS and EMT, etc.  The only condition for the co-simulation algorithm to work is that the RAS method does not stagnate (i.e. there is not $1$ as an eigenvalue of $P_\Gamma$) or that the divergence of the RAS is not too strong to face numerical conditioning problems.}

The Aitken's acceleration technique can use the error operator in the whole domain and accelerate the whole solution at once, or compute the error operator only on artificial boundaries accelerate only the interface values and use these  converged values to obtain the whole solution with local  resolutions. The error operator can be computed  algebraically but  can also be computed numerically. To construct the error operator corresponding to $n$ values, it is necessary to perform $n+1$ Schwarz iterations. For this reason, most of the time, the error operator will be built only to speed up the values of the artificial boundaries.

\subsubsection{Numerical Computation of the error Operator}
 As the translation operators complexify the analytical computation of the error operator $P_\Gamma$,  its  numerical computation is preferred. The way to compute the error operator $P_\Gamma$ numerically is as follows:
\begin{prop}
By defining the error between $(k+1)$ and $(k)$ iterates on the global interface $\Gamma$ as $e^{k}=((w^{N+1,(k+1)}_{1,e_{ts}} -w^{N+1,(k)}_{1,e_{ts}})^T, ( \mathbb{W}_{2,e_{emt}}^{N+1,(k+1)}-\mathbb{W}_{2,e_{emt}}^{N+1,(k)})^T)^T$, 
The operator $P_\Gamma \in \mathbb{R}^{n_{\Gamma} \times n_{\Gamma}}$ can be computed algebraically from Eq. \eqref{RASerrorGamma} after  $n_{\Gamma}+1$ iterations, if the matrix $[e^{n_{\Gamma}},\ldots, e^1]$ is non-singular as:  
\begin{eqnarray}
P_\Gamma&=&[e^{n_{\Gamma}+1},\ldots, e^2][e^{n_{\Gamma}},\ldots, e^1]^{-1}. \label{P_algebraic}
\end{eqnarray}
\end{prop}

This error operator should only be calculated for the first-time step. Unless the topology changes or there is a non-linearity in which case the operator must be recomputed.

\subsubsection{Impact of heterogeneity in the choice of acceleration strategy}
We seek to accelerate convergence at the interfaces between two subdomains, a TS $W_{1_{ts}}$and an EMT $W_{2_{emt}}$. We will therefore first accelerate the convergence of the interfaces:  $ Z_{{1e}_{emt}}$ and $\mathbb{W}_{{2,e}_{ts}}$ . \\
 Let's calculate the size of the interface values, note $n_{i,e}$ the size of interface $i,e$: \\
First let's calculate for $  Z_{{1e}_{emt}}$  we just need to speed up the values taken at the last macro time step, the "Active Boundary Conditions": $\mathbb{Z}_{{1e}_{emt}}$  so the sample size to be accelerated is $n_{1,e}\times m$ (with $\Delta t_{ts}= m \Delta t_{emt}$).\\

Then let's calculate the size of the interface $\mathbb{W}_{{2,e}_{ts}}$, $2K \times n_{2,e}$, (with K the number of harmonics that we chose to keep).\\

 So the whole interface between these two sub-domains is of size $n_{1,e}\times m + 2K \times n_{2,e}$, so to speed up this interface it will be necessary to realize $n_{1, e}\times m + 2K \times n_{2,e}+1$ Schwarz iteration, which can be huge. Indeed, $m$ is often around a hundred, K is rarely greater than 2. We cannot afford to do so many iterations of Schwarz. Two solutions are available to us: 
 \begin{itemize}
\item Accelerate the values coming from the EMT after translation, in this way, the size of the whole interface to be accelerated will be $ 2K \times n_{1,e} + 2K \times n_{2,e}$.
\item One can choose to accelerate only the interface values coming from  the TS side, i.e. to compute the error operator $P$ linked only to the $\mathbb{W}_{{2e}_{ts}}$ values. Like so, an error operator  of size $ 2K \times n_{2,e}$ will be computed by performing $ 2K \times n_{2,e}+1$ iterations. Once the true values of $\mathbb{W}_{{2e}_{ts}}$ are obtained, they are used to compute the local resolution of the EMT part $\mathbb{W}_{2_{emt}}$. These EMT values are in turn used for the local resolution of the TS part $w^{N+1,(k+1)}_{1_{ts}}$ .
\end{itemize}
We'll prefer the second option because it further limits the number of Schwarz iterations needed.

 \section{Heterogeneous EMT-TS Numerical Results \label{shourick_contrib_Sec5}}
 
\begin{figure}[h]
\begin{minipage}{13cm}
\begin{minipage}{4cm}
\begin{tikzpicture}[scale=0.6]
 \draw[black!70] (-0.7,3) node[above]{$W$};
 \draw (-4,2) -- (-3,2);
 \draw (-2,2) -- (-1.4,2);
 \draw(0,2) -- (1.5,2);
 \draw (2,2) -- (3,2);
 \draw (3,2) -- (3,-2); 
 \draw (-4,-2) -- (-2.5,-2);
 \draw (-2,-2) -- (-1,-2);
 \draw (0,-2) -- (1,-2);
 \draw (2.4,-2) -- (3,-2);
 \draw (-4,2) -- (-4,-2);

%

  \draw[DodgerBlue] (-4,2) node[above]{\scriptsize $2$};
   \draw[DodgerBlue] (-1.6,2) node[above]{\scriptsize $3$};
   \draw[DodgerBlue] (0.6,2) node[above]{\scriptsize $4$};
   \draw[DodgerBlue] (3,2) node[above]{\scriptsize $5$};
   \draw[DodgerBlue] (-1.6,-2) node[below]{\scriptsize $7$};
   \draw[DodgerBlue] (0.6,-2) node[below]{\scriptsize $6$};
   \draw[DodgerBlue] (-4,-2) node[left]{\scriptsize $1$};
  
 \draw[thick] (1.5,2.5) -- (1.5,1.5);
 \draw[thick] (2,2.5) -- (2,1.5);

  \draw (1.7,2.4) node[above]{\footnotesize  $C_1$};
   
 \draw[thick] (-2,-2.5) -- (-2,-1.5);
 \draw[thick] (-2.5,-2.5) -- (-2.5,-1.5);
 
  \draw (-2.3,-1.6) node[above]{\footnotesize $C_2$};

  \draw[ thick] (-4,-2) -- (-4,-2.5);
 \draw[thick] (-4.5,-2.5) -- (-3.6,-2.5);
  \draw[thick] (-4.5,-2.5) -- (-4.3,-2.7);
 \draw[thick] (-4.2,-2.5) -- (-4,-2.7);
  \draw[thick] (-3.9,-2.5) -- (-3.7,-2.7);
   \draw[thick] (-3.6,-2.5) -- (-3.4,-2.7);
   
   
    \draw[thick] (-1.4,2) -- (-1.2,2.3);
    \draw[thick] (-1.2,2.3) -- (-1,1.7);
     \draw[thick] (-1,1.7) -- (-0.8,2.3);
 \draw[thick] (-0.8,2.3) -- (-0.6,1.7);
  \draw[thick] (-0.6,1.7) -- (-0.4,2.3);
  \draw[thick] (-0.4,2.3) -- (-0.2,1.7);
  \draw[thick] (-0.2,1.7) -- (0,2);
  
   \draw (-0.7,2.4) node[above]{\footnotesize $R_1$};
  
  \draw[thick] (1,-2) -- (1.2,-1.7);
  \draw[thick] (1.2,-1.7) -- (1.4,-2.3);
 \draw[thick] (1.4,-2.3) -- (1.6,-1.7);
 \draw[thick] (1.6,-1.7) -- (1.8,-2.3);
  \draw[thick] (1.8,-2.3) -- (2,-1.7);
  \draw[thick] (2,-1.7) -- (2.2,-2.3);
  \draw[thick] (2.2,-2.3) -- (2.4,-2);
   \draw (1.7,-1.6) node[above]{\footnotesize $R_2$};
  
  \draw[thick] (-4,0) circle(0.5);
  \draw[thick] (-4,0.5)-- (-4,-0.5);
  \draw[->,thick] (-3.45,-0.5)-- (-3.45,0.5);
  \draw (-3.45,0) node[right]{\footnotesize E cos $\omega t = \beta$};

   
   \begin{scope}[shift={(-3.5,2)},rotate=90]
{
 \draw[black!5!yellow!10!red!8!,opacity=0.3] (-0.2,0) rectangle (0.4,-1.5);
 \foreach \r in {0,...,2}
 {
  \draw[thick,scale=1/3,shift={(0,-\r)}]
	(0,0) .. controls ++(2,0) and ++(1,0) ..
	++(0,-1.5) .. controls ++(-1,0) and ++(-0.5,0) ..
	++(0,0.5);
 }
 \draw[thick,scale=1/3,shift={(0,-3)}] (0,0) .. controls ++(2,0) and ++(1,0) .. ++(0,-1.5);
}

\end{scope}
\draw (-2.8,2.4) node[above]{\footnotesize  $L_1$};

   \begin{scope}[shift={(-1.2,-2)},rotate=90]
{
 \draw[black!7!DodgerBlue!11,opacity=0.80] (-0.2,0) rectangle (0.4,-1.5);
 \foreach \r in {0,...,2}
 {
  \draw[thick,scale=1/3,shift={(0,-\r)}]
	(0,0) .. controls ++(2,0) and ++(1,0) ..
	++(0,-1.5) .. controls ++(-1,0) and ++(-0.5,0) ..
	++(0,0.5);
 }
 \draw[thick,scale=1/3,shift={(0,-3)}] (0,0) .. controls ++(2,0) and ++(1,0) .. ++(0,-1.5);
}
\end{scope}
 \draw (-0.5,-1.6) node[above]{\footnotesize $L_2$};

\coordinate (h) at (0.3,2.1);;
\coordinate (i) at (3.5,2.1);
\coordinate (j) at (3.5,-2.2);
\coordinate (k) at (-1.5,-2.2);
\coordinate (hk) at (-0.4,0.15);

\coordinate (m) at (-4.4,2.2);
\coordinate (mn) at (-2,2.9);
\coordinate (n) at (0.2,2.2);
\coordinate (no) at (-0.46,-0.32);
\coordinate (o) at (-1.5,-2.1);
\coordinate (p) at (-4.4,-2.15);
\coordinate (pm) at (-4.85,0);

\draw[dotted,draw=red!60!yellow!35!black,opacity=0.55]  (m) .. controls +(1,0.5) and +(-1,0.08) .. (mn)
               .. controls +(1,0.08) and +(-0.8,0.5) .. (n)
               .. controls +(0.5,-0.8) and +(-0.19,2) .. (no)
               .. controls +(-0.15,-1.5) and +(0.55,0.72) .. (o)
               .. controls +(-1,-0.2) and +(1,-0.2) .. (p)
                .. controls +(-0.5,0.8) and +(0.1,-0.8) .. (pm)
               .. controls +(0.1,0.8) and +(-0.5,-0.8) .. (m); 

  \draw[dotted,draw=black,opacity=0.5](h) .. controls +(0.8,0.6) and +(-0.8,0.6) .. (i)
               .. controls +(0.8,-0.9) and +(0.7,0.9) .. (j)
               .. controls +(-0.9,-0.8) and +(0.9,-0.7) .. (k)
                .. controls +(0.87,0.8) and +(0,-1.3) .. (hk)
                .. controls +(0,1.3) and +(-0.1,-1) .. (h); 
\end{tikzpicture}
\begin{eqnarray}
v_1&=&0   \nonumber\\
v_2-v_1-E-Z_s i_{12} &=& 0  \nonumber \\
v_3-v_2-L_1 \dfrac{d{i}_{23}}{dt} &=& 0  \nonumber \\
v_4-v_3 -R_1 i_{34} &=& 0 \nonumber
\end{eqnarray}
\end{minipage}
\hfill
\begin{minipage}{4cm}
\begin{eqnarray}
C_1 (\dfrac{d{v}_5}{dt}-\dfrac{d{v}_4}{dt})-i_{45} &=&0  \nonumber\\
v_6-v_5-R_2 i_{56} &=& 0  \nonumber\\
v_7-v_6-L_2 \dfrac{d{i}_{67}}{dt} &=& 0  \nonumber \\
C_2 (\dfrac{d{v}_1}{dt}-\dfrac{d{v}_7}{dt})-i_{71} &=& 0  \nonumber  \\
i_{12}-i_{23}&=& 0 \nonumber  \\
i_{23}-i_{34}&=& 0  \nonumber \\
i_{34}-i_{45}&=& 0  \nonumber \\
i_{45}-i_{56}&=& 0  \nonumber \\
i_{56}-i_{67}&=& 0   \nonumber\\
i_{67}-i_{71}&=& 0   \nonumber
\end{eqnarray}
\end{minipage}
\end{minipage}
\caption{Linear RLC circuit and its associated  EMT modeling  DAE system with\\ $x=\left\{ v_1, i_{23}, v4, v5, i_{67}, v_7\right\}$ and $y=\left\{ v_2, i_{12},  v_3, i_{34}, i_{45}, i_{56}, v_6, i_{71}\right\}$. $L1=L2=0.7$, \\ $C1=C2=1.10^{-6}$, $R1=R2=77$, $Zs=1.10^{-6}$, $\omega_0=2\pi\,50$. \label{FigRLC}}
\end{figure}
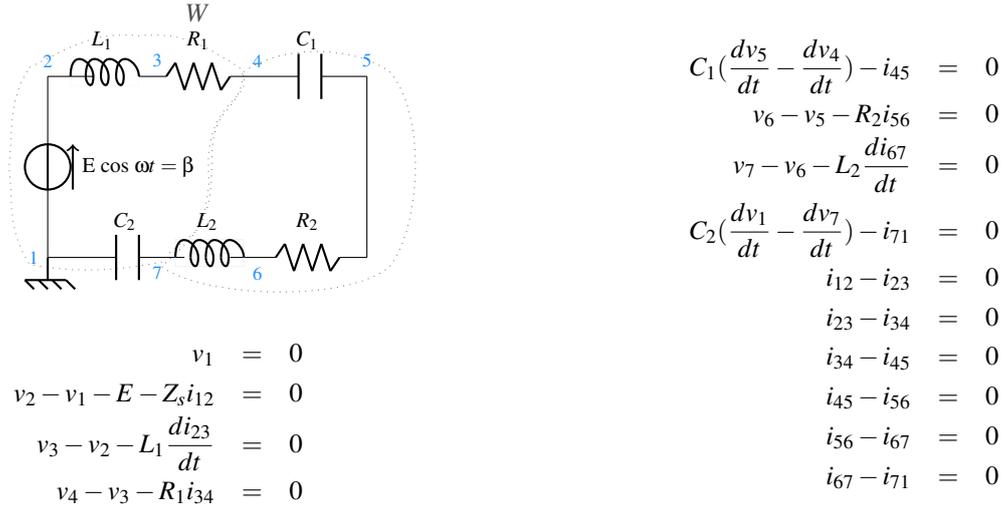

We consider a linear RLC circuit of Figure \ref{FigRLC}. This is a single-loop circuit, so the different sub-domains will necessarily be strongly coupled. It contains the basic components used in the electrical field to model most phenomena. The DAE associated with the small circuit is written component by component. We do not make any simplification usually done when there are several resistors or inductors in series, to keep more equations. Instead of considering the single loop current, we will consider that there is one current per component, which will allow us to observe more phenomena like error propagation. In order to keep the right number of equations and variables, and to be able to put two inductors in series,  we add a degree of freedom (that is, we remove one of the current equality equations).

\begin{figure}[h]
\centering
\begin{minipage}{15cm}

\begin{minipage}{6.4cm}
\begin{tikzpicture}[scale=0.7]

 \draw[black!70] (-0.7,3) node[above]{$W = W_1 \cup W_2$};
 \draw[black!70!yellow!65!red!65!] (2.55,2.7) node[above]{$W_2$};
 \draw[black!60!DodgerBlue!70] (-3.65,2.7) node[above]{$W_1$};
 \draw (-4,2) -- (-3,2);
 \draw (-2,2) -- (-1.4,2);
 \draw(0,2) -- (1.5,2);
 \draw (2,2) -- (3,2);
 \draw (3,2) -- (3,-2); 
 \draw (-4,-2) -- (-2.5,-2);
 \draw (-2,-2) -- (-1,-2);
 \draw (0,-2) -- (1,-2);
 \draw (2.4,-2) -- (3,-2);
 
  \draw (-4,2) -- (-4,-2);

   \fill[Green] (-4,2) circle(0.05);
   \fill[Green] (-1.6,2) circle(0.05);
   \fill[Green] (0.6,2) circle(0.05);
   \fill[Green] (3,2) circle(0.05);
   \fill[Green] (-1.6,-2) circle(0.05);
   \fill[Green] (0.6,-2) circle(0.05);
   \fill[Green] (-4,-2) circle(0.05);

  \draw[Green] (-4,2) node[above]{\scriptsize $2$};
   \draw[Green] (-1.6,2) node[above]{\scriptsize $3$};
   \draw[Green] (0.6,2) node[above]{\scriptsize $4$};
   \draw[Green] (3,2) node[above]{\scriptsize $5$};
   \draw[Green] (-1.6,-2) node[below]{\scriptsize $7$};
   \draw[Green] (0.6,-2) node[below]{\scriptsize $6$};
   \draw[Green] (-4,-2) node[left]{\scriptsize $1$};
  
 \draw[thick] (1.5,2.5) -- (1.5,1.5);
 \draw[thick] (2,2.5) -- (2,1.5);

  \draw (1.7,2.4) node[above]{\footnotesize  $C_1$};
   
 \draw[thick] (-2,-2.5) -- (-2,-1.5);
 \draw[thick] (-2.5,-2.5) -- (-2.5,-1.5);
 
  \draw (-2.3,-1.6) node[above]{\footnotesize $C_2$};

  \draw[ thick] (-4,-2) -- (-4,-2.5);
 \draw[thick] (-4.5,-2.5) -- (-3.6,-2.5);
  \draw[thick] (-4.5,-2.5) -- (-4.3,-2.7);
 \draw[thick] (-4.2,-2.5) -- (-4,-2.7);
  \draw[thick] (-3.9,-2.5) -- (-3.7,-2.7);
   \draw[thick] (-3.6,-2.5) -- (-3.4,-2.7);
   
   
    \draw[thick] (-1.4,2) -- (-1.2,2.3);
    \draw[thick] (-1.2,2.3) -- (-1,1.7);
     \draw[thick] (-1,1.7) -- (-0.8,2.3);
 \draw[thick] (-0.8,2.3) -- (-0.6,1.7);
  \draw[thick] (-0.6,1.7) -- (-0.4,2.3);
  \draw[thick] (-0.4,2.3) -- (-0.2,1.7);
  \draw[thick] (-0.2,1.7) -- (0,2);
  
   \draw (-0.7,2.4) node[above]{\footnotesize $R_1$};
  
  \draw[thick] (1,-2) -- (1.2,-1.7);
  \draw[thick] (1.2,-1.7) -- (1.4,-2.3);
 \draw[thick] (1.4,-2.3) -- (1.6,-1.7);
 \draw[thick] (1.6,-1.7) -- (1.8,-2.3);
  \draw[thick] (1.8,-2.3) -- (2,-1.7);
  \draw[thick] (2,-1.7) -- (2.2,-2.3);
  \draw[thick] (2.2,-2.3) -- (2.4,-2);
   \draw (1.7,-1.6) node[above]{\footnotesize $R_2$};
  
  \draw[thick] (-4,0) circle(0.5);
  \draw[thick] (-4,0.5)-- (-4,-0.5);
  \draw[->,thick] (-3.45,-0.5)-- (-3.45,0.5);
  \draw (-3.45,0) node[right]{\footnotesize E cos $\omega t = \beta$};

   
   \begin{scope}[shift={(-3.5,2)},rotate=90]
{
 \draw[black!6.8!DodgerBlue!11,opacity=0.80,fill=black!7!DodgerBlue!11,opacity=0.80] (-0.2,0) rectangle (0.4,-1.5);
 \foreach \r in {0,...,2}
 {
  \draw[thick,scale=1/3,shift={(0,-\r)}]
	(0,0) .. controls ++(2,0) and ++(1,0) ..
	++(0,-1.5) .. controls ++(-1,0) and ++(-0.5,0) ..
	++(0,0.5);
 }
 \draw[thick,scale=1/3,shift={(0,-3)}] (0,0) .. controls ++(2,0) and ++(1,0) .. ++(0,-1.5);
}

\end{scope}
\draw (-2.8,2.4) node[above]{\footnotesize  $L_1$};

   \begin{scope}[shift={(-1.2,-2)},rotate=90]
{
 draw[black!5!yellow!10!red!8!,opacity=0.3,fill=black!5!yellow!10!red!8!,opacity=0.7 (-0.2,0) rectangle (0.4,-1.5);
 \foreach \r in {0,...,2}
 {
  \draw[thick,scale=1/3,shift={(0,-\r)}]
	(0,0) .. controls ++(2,0) and ++(1,0) ..
	++(0,-1.5) .. controls ++(-1,0) and ++(-0.5,0) ..
	++(0,0.5);
 }
 \draw[thick,scale=1/3,shift={(0,-3)}] (0,0) .. controls ++(2,0) and ++(1,0) .. ++(0,-1.5);
}
\end{scope}
 \draw (-0.5,-1.6) node[above]{\footnotesize $L_2$};

\draw[blue,thick,->] (0.1,2.2) -- (0.5,2.2);
\draw[red,thick,<-] (-1.4,-2.1) -- (-1,-2.1);
\draw[blue](-0.1,2.1) node[above]{\footnotesize ${\bf i_{34},v_3}$};
\draw[red](-0.8,-1.9) node[below]{\footnotesize ${\bf i_{67},v_6}$};

\coordinate (h) at (0.3,2.1);;
\coordinate (i) at (3.5,2.1);
\coordinate (j) at (3.5,-2.2);
\coordinate (k) at (-1.5,-2.2);
\coordinate (hk) at (-0.4,0.15);

\coordinate (m) at (-4.4,2.2);
\coordinate (mn) at (-2,2.9);
\coordinate (n) at (0.2,2.2);
\coordinate (no) at (-0.46,-0.32);
\coordinate (o) at (-1.5,-2.1);
\coordinate (p) at (-4.4,-2.15);
\coordinate (pm) at (-4.85,0);

\draw[dotted,draw=black,fill=black!5!DodgerBlue!11,opacity=0.5]  (m) .. controls +(1,0.5) and +(-1,0.08) .. (mn)
               .. controls +(1,0.08) and +(-0.8,0.5) .. (n)
               .. controls +(0.5,-0.8) and +(-0.19,2) .. (no)
               .. controls +(-0.15,-1.5) and +(0.55,0.72) .. (o)
               .. controls +(-1,-0.2) and +(1,-0.2) .. (p)
                .. controls +(-0.5,0.8) and +(0.1,-0.8) .. (pm)
               .. controls +(0.1,0.8) and +(-0.5,-0.8) .. (m); 

  \draw[dotted,draw=red!60!yellow!35!black,fill=black!5!yellow!10!red!15!,opacity=0.55](h) .. controls +(0.8,0.6) and +(-0.8,0.6) .. (i)
               .. controls +(0.8,-0.9) and +(0.7,0.9) .. (j)
               .. controls +(-0.9,-0.8) and +(0.9,-0.7) .. (k)
                .. controls +(0.87,0.8) and +(0,-1.3) .. (hk)
                .. controls +(0,1.3) and +(-0.1,-1) .. (h); 

\end{tikzpicture}
\end{minipage}
\hfill
\begin{minipage}{6.4cm}
\includegraphics[scale=0.29]{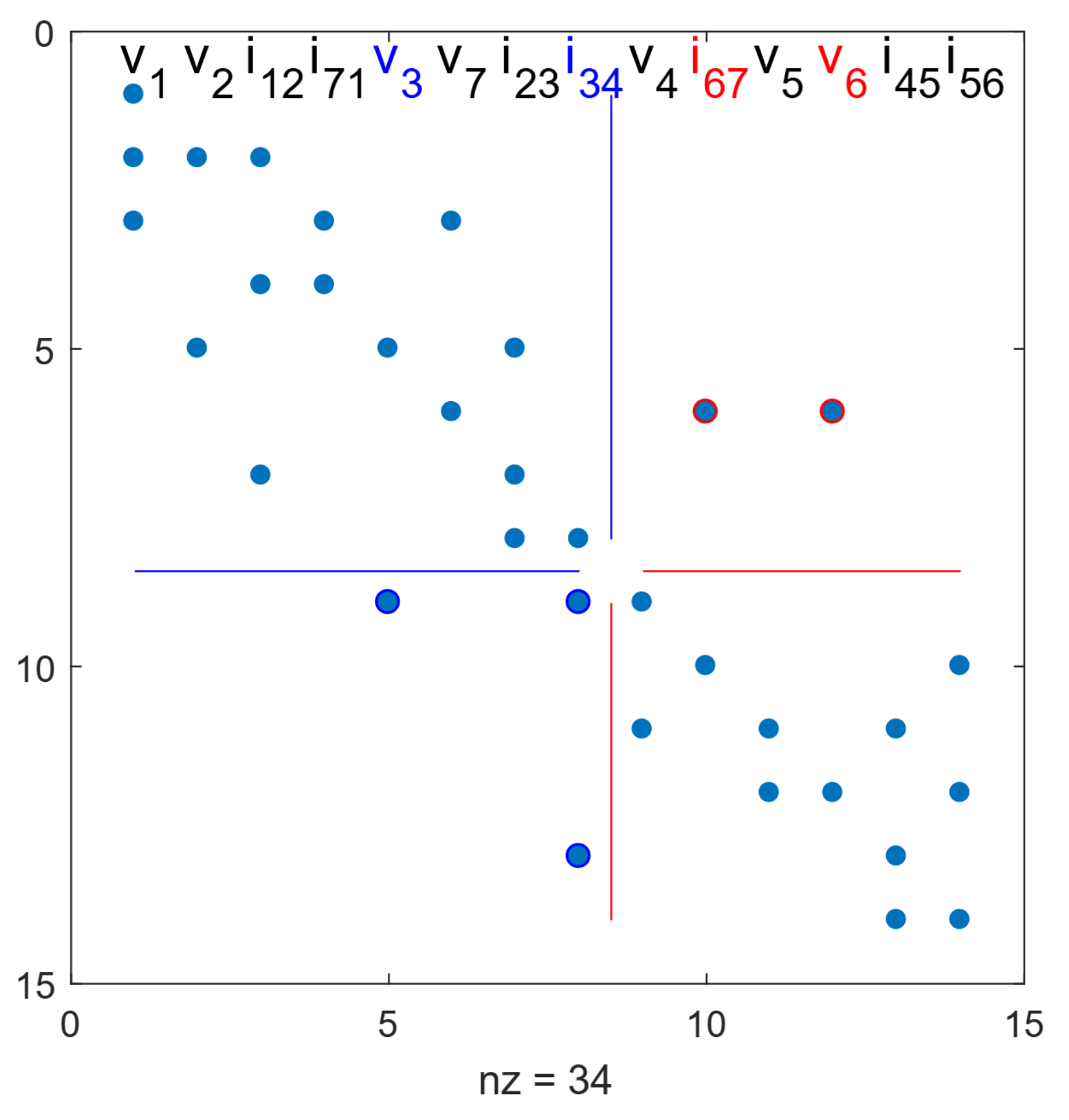}
\end{minipage}
\end{minipage}

\caption{Graph partitioning of the RLC circuit in two subdomains and the associated matrix partioning without overlap (top) and with overlap of 1 (bottom). EMT case \label{RLCCut}}
\end{figure}
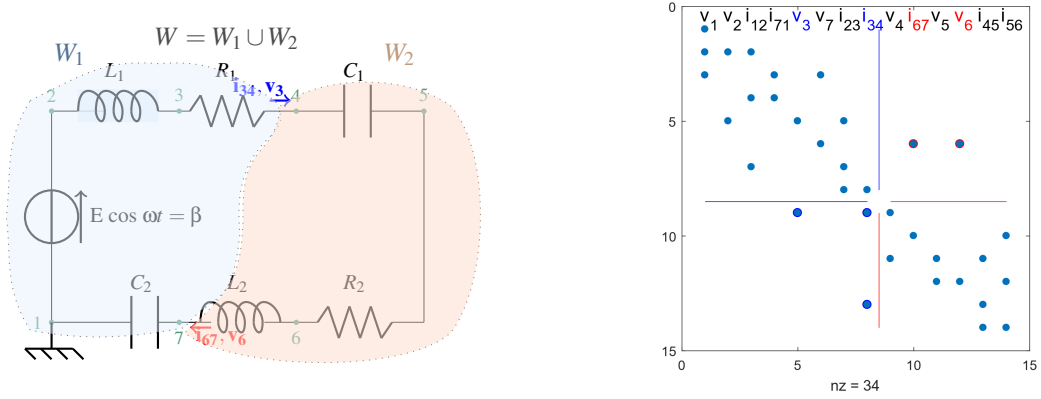

Figure \ref{RLCCut} is an example of this cutting for a small RLC circuit. The small linear system associated with the RLC circuit is partitioned into two subdomains using non-overlapping graph partitioning (Figure \ref{RLCCut} on top) and with an overlap of 6 (Figure \ref{RLCCut} bottom). This use of graph partitioning aims to have an equivalent computational load for each subdomain in a homogeneous case (the two subdomains are modeled in the same EMT or TS way). In this decomposition, convergence considerations are not taken into account. In this example the domain $W$ is cut into two with $W_1^0=\{v_1,v_2,i_{12},i_{71},v_3,v_7,i_{23},i_{34} \}$ and $W_2^0=\{v_4,i_{67},v_5,v_6,i_{56}\}$ , $W_1^1=\{v_1,v_2,i_{12},i_{71},v_3,v_7,i_{23},i_{34} ,i_{71},v{4},i_{67}\}$  and $W_2^1=\{v_3,v_7,i_{34},v_4,i_{67},v_5,v_6,i_{56}\}$ 

 In the EMT case, each subdomain needs two values from the other to solve its equations. In the TS case, as we choose to solve harmonics 0  and 1 and to solve real and imaginary part apart (the imaginary par of the harmonic 0 is always 0), each subdomain needs six values from the other to solve its own equations.

First,  we study  numerically the convergence of the method. As the RAS method diverged in the homogeneous case we can expect the same kind of results. Secondly, we study  the results obtained  by the heterogeneous co-simulation in order to see the gain of the EMT-TS modeling over the TS modeling. 
 \subsection{Heterogeneous EMT-TS RAS convergence results}
 
 \subsubsection{Effect of circuit topology on the convergence}
 In order to see only the effect of  the circuit topology (i.e the effect of the circuit components values) on the convergence, we use the classic FFT to translate from EMT to TS  keeping $\Delta t_{ts}$ to be a period. We also do not  smooth the TS to EMT translation either, but simply recombine the harmonics (i.e. $\alpha = 0$). The impact on the convergence of the method of these modifications to the translation operators will be given in a second time.\\
We know that the value of the components influences the convergence of the method (\cite{Pade_Tischendorf_WR_CVGCriterion_NUMAlgo_2019}), and more particularly the inductance and the capacitance have a very strong impact.

\begin{figure}[h!]
\centering
\includegraphics[scale=0.4]{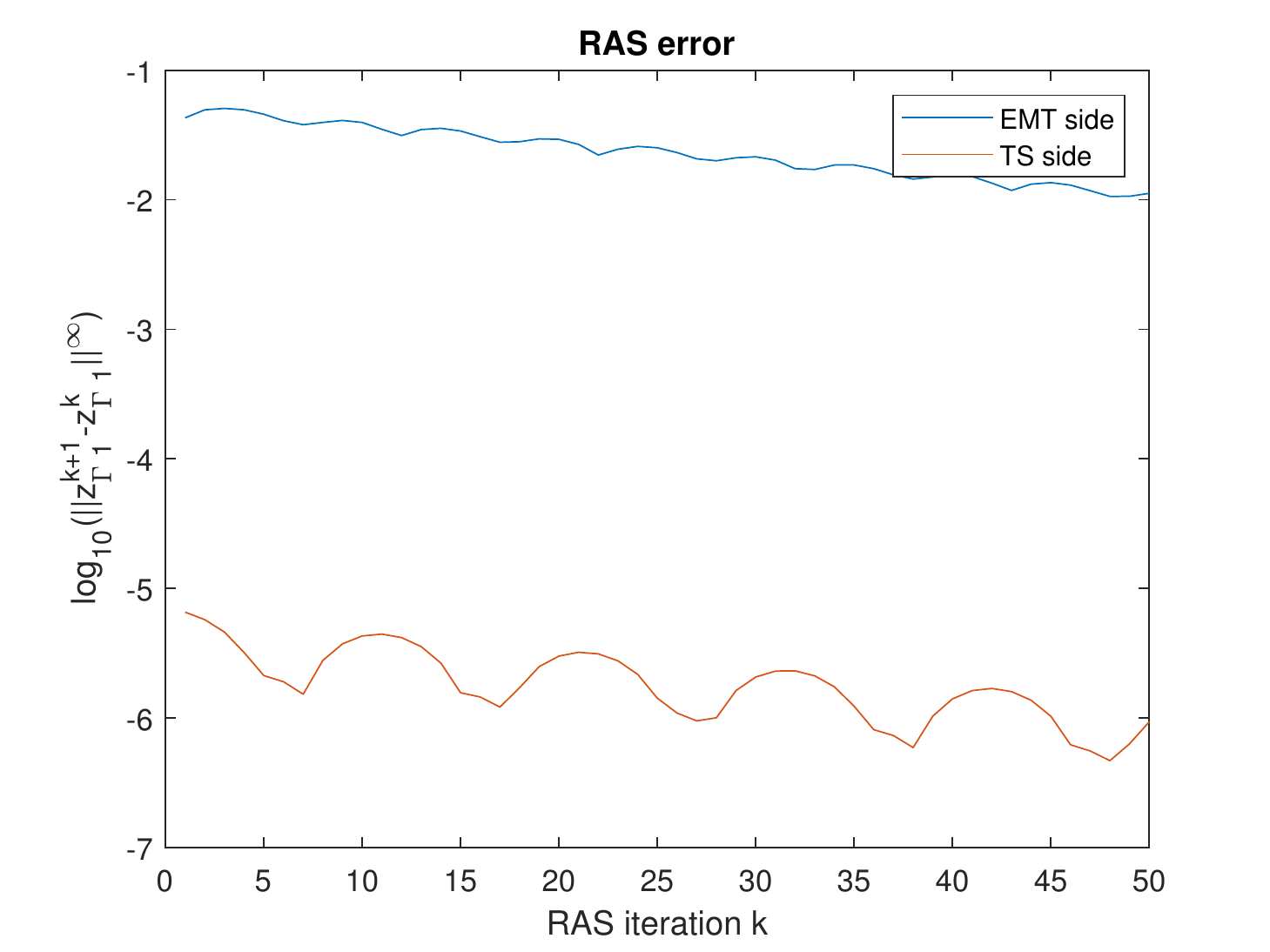}
\includegraphics[scale=0.4]{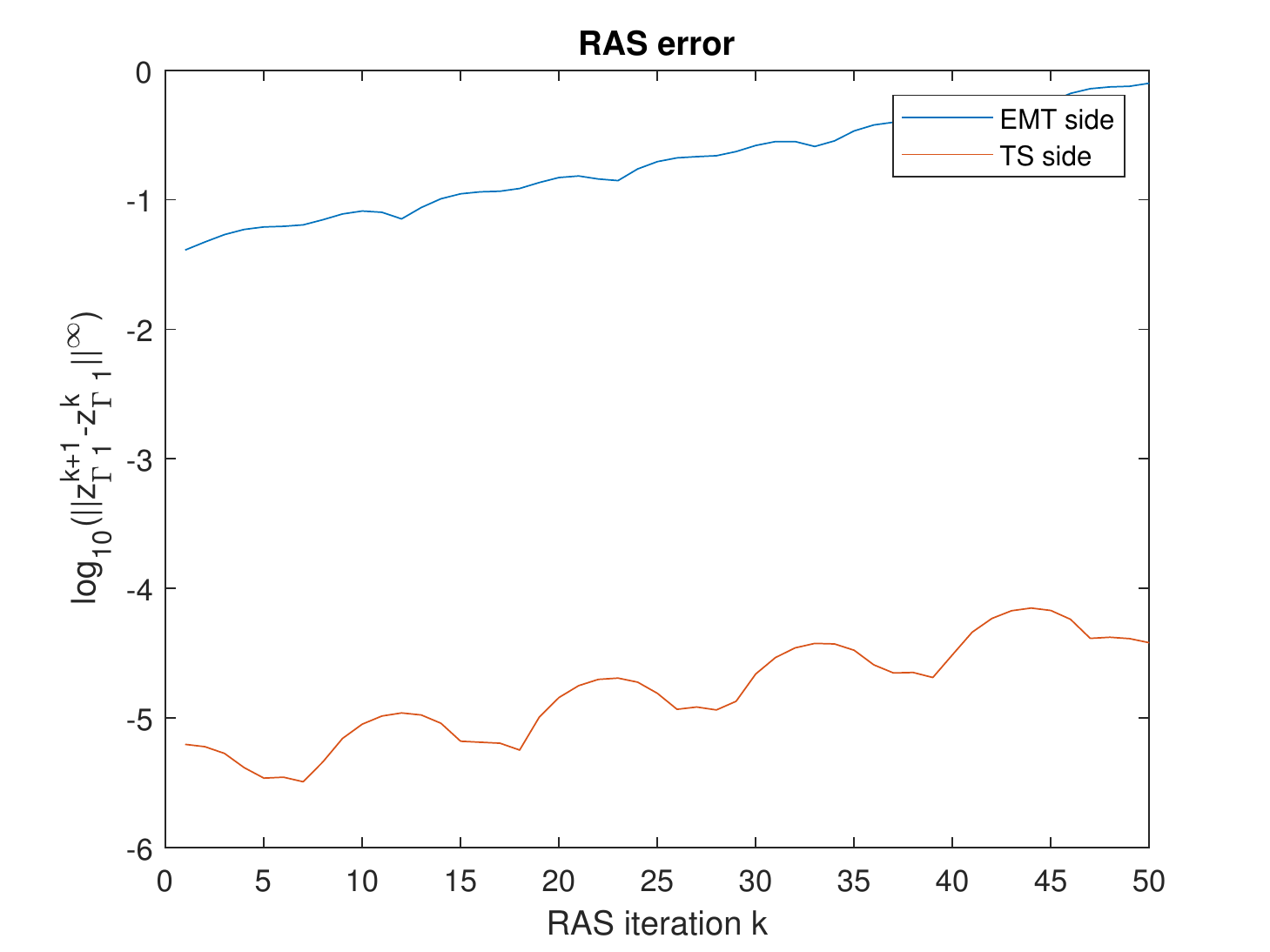}
\caption{ Convergence of the heterogeneous EMT-TS RAS error between two consecutive iterates with voltage source in the TS subsystem and with $C_1=C_2=1.10^{-6}$, $R_1=R_2=7$,$L_1=0.07$ and we have $L_2=1$  on the left and  $L_2=0.07$ on the right.}
\label{TopologyErrInduct}
\end{figure}

\begin{figure}[h!]
\centering
\includegraphics[scale=0.4]{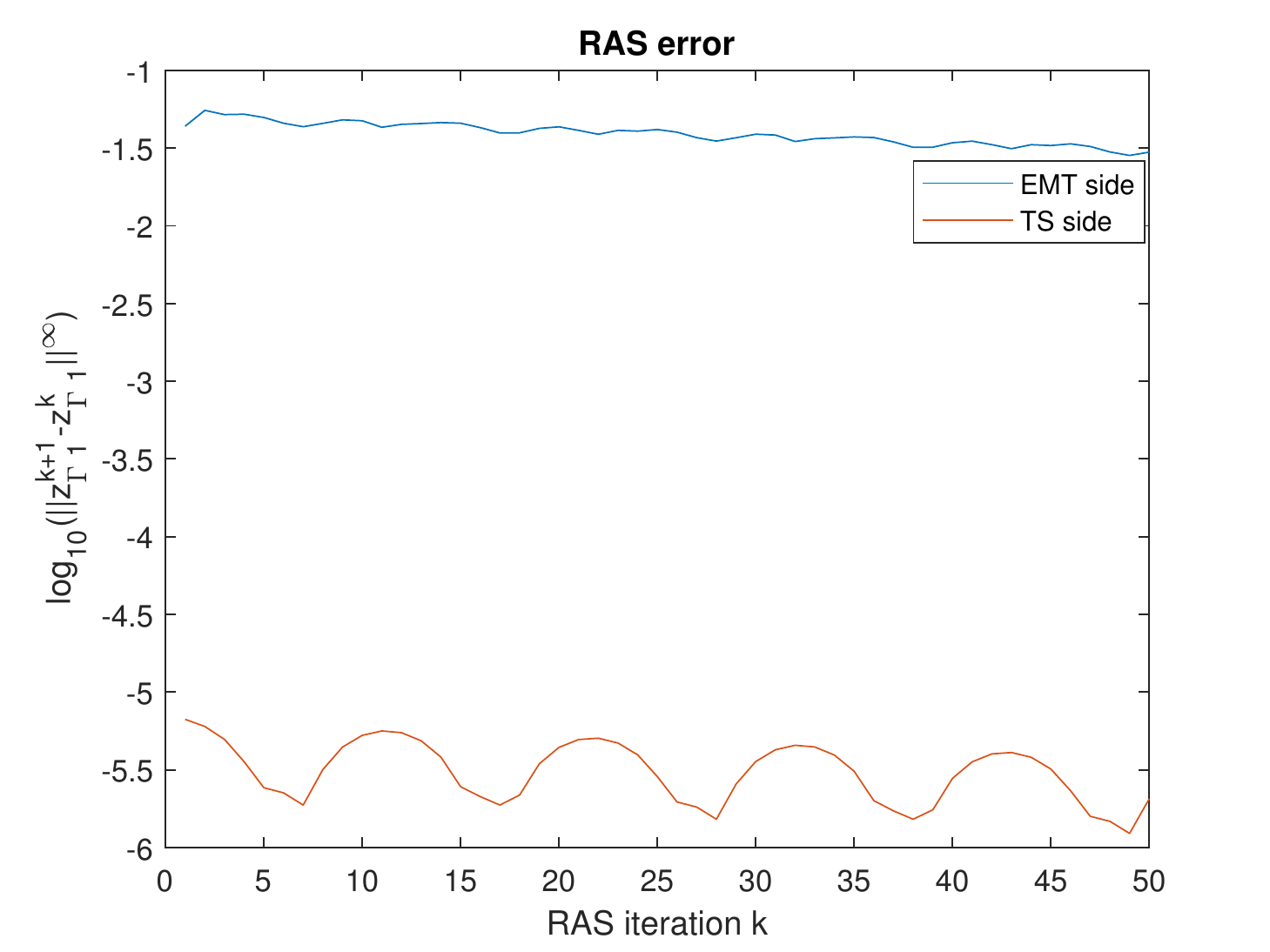}
\includegraphics[scale=0.4]{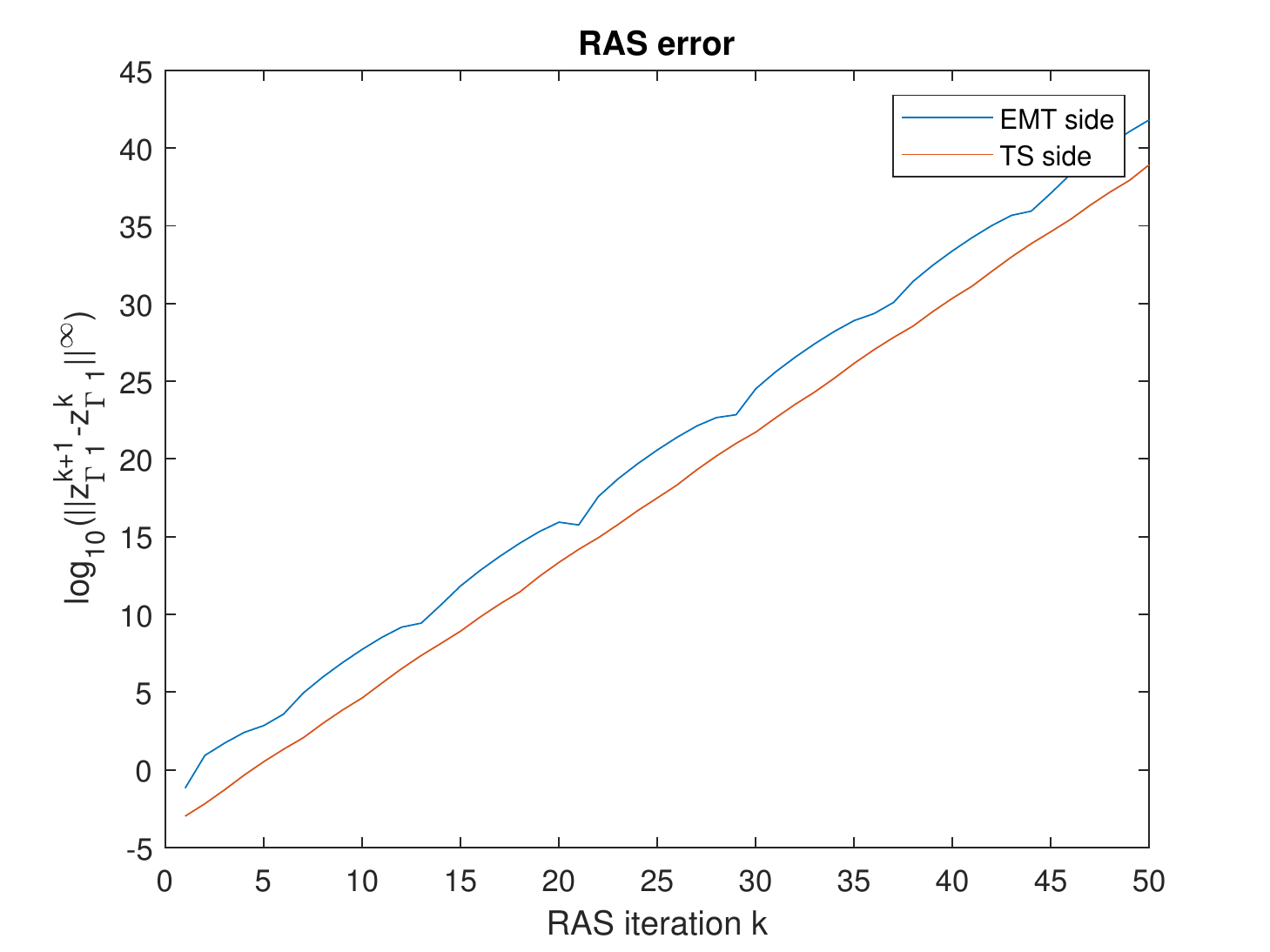}
\caption{Convergence of the heterogeneous EMT-TS RAS error between two consecutive iterates with voltage source in the TS subsystem and with, $R_1=R_2=7$,$L_1=L_2=0.4$ and we have $C_1=C_2=1.10^{-6}$ on the left and  $C_1=C_2=1.10^{-4}$ on the right.}
\label{TopologyErrCapa}
\end{figure}

Figure \ref{TopologyErrInduct} shows the error between two heterogeneous EMT-TS RAS iterations by changing the value of the inductance $L_2$ in the circuit. It exhibits a convergent method for $L_2=1$ while the method diverges for $L_2=0.07$ the other components being set.
Figure \ref{TopologyErrCapa} shows the error between two heterogeneous EMT-TS RAS iterations by changing the values of the capacitances $C_1$ and $C_2$ in the circuit. It exhibits a convergent method for $C_1=C_2=1.10^{-6}$ while the method strongly diverges for $C_1=C_2=1.10^{-4}$ the other components being set.
The choice of the part that is simulated with the TS or with the EMT also has an impact and more particularly the place where the voltage source is located has an impact in the convergence of the method.

\begin{figure}[h!]
\centering
\includegraphics[scale=0.4]{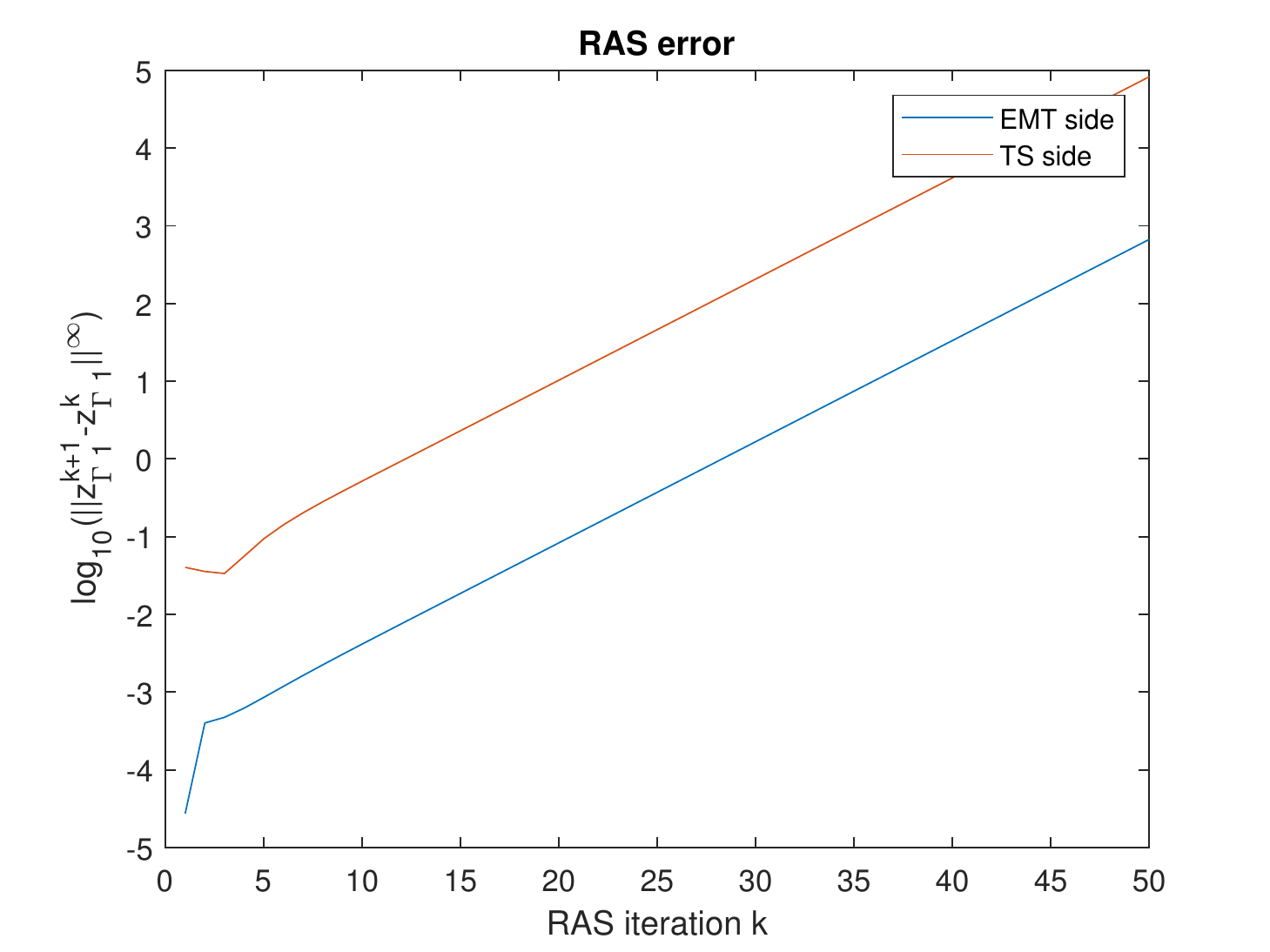}
\includegraphics[scale=0.4]{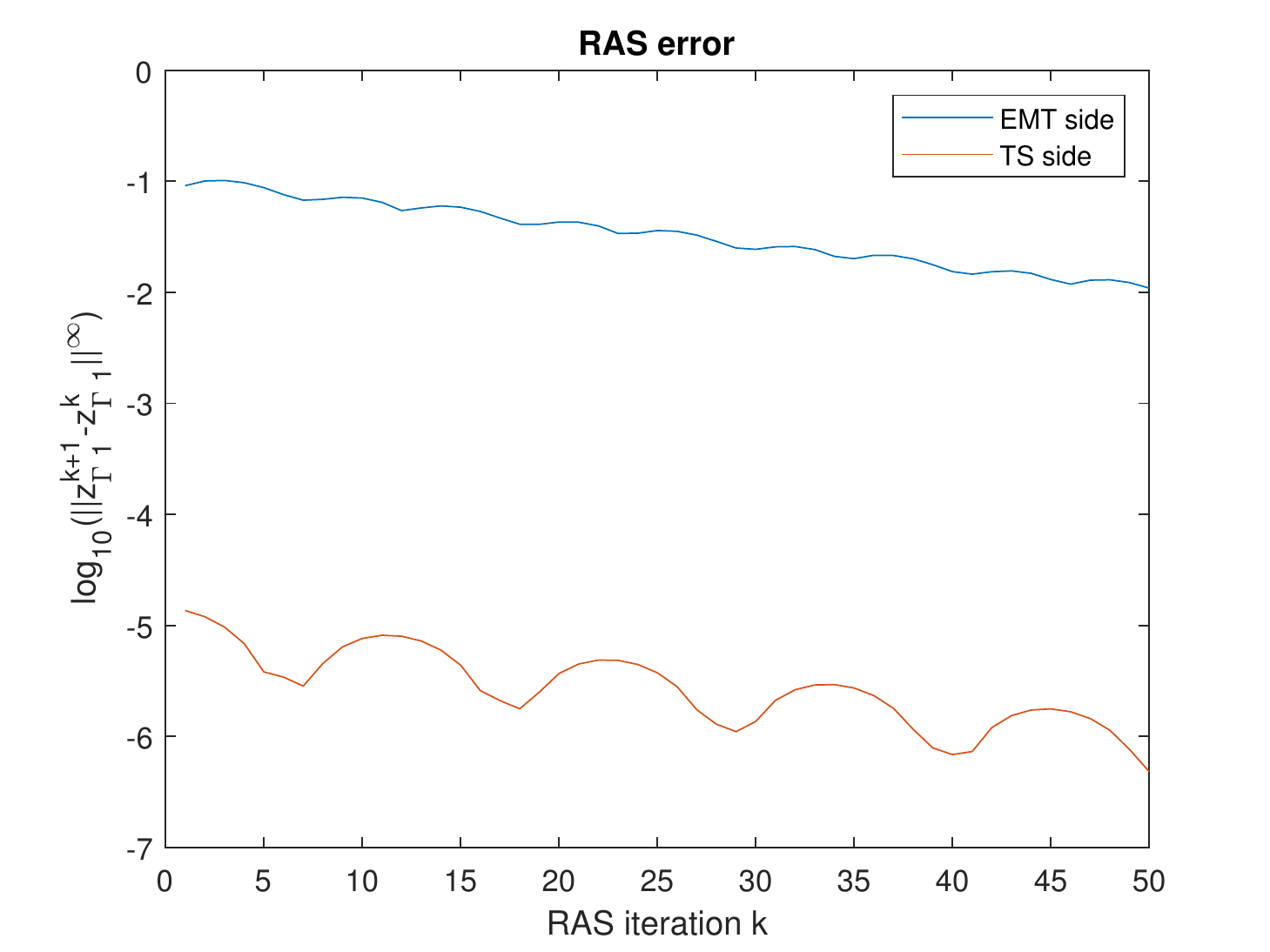}
\caption{Convergence of the heterogeneous EMT-TS RAS error between two consecutive iterates with, $R_1=R_2=7$,$L_1=L_2=0.5$ ,$C_1=C_2=1.10^{-6}$. On the left the source is in the EMT subsystem and on the right the source is in the TS subsystem.}
\label{SourceErr}
\end{figure}

Figure \ref{SourceErr}, gives the error between two consecutive  heterogeneous EMT-TS RAS iterates following that the voltage source is modeled with EMT (left) or is  modeled in TS (right). It exhibits a divergent method when the source is in the EMT part while the method is convergent if it is located in the TS part.
These results clearly demonstrate the need for the heterogeneous EMT-TS RAS method to have the Aitken's acceleration technique as in the homogeneous RAS case to be independent of the circuit topology.

 \subsubsection{Effect of the EMT Time steps $\Delta t_{emt}$ on the convergence}
We have seen with Figure \ref{SourceErr} that the convergence depends on the subsystem in which the source is located. However, the impact of the localization of the voltage source is coupled to the size of the EMT time steps. Indeed, a modification of the size of the EMT time steps does not have the same impact depending on the location of the voltage source. 
\begin{figure}[h!]
\centering
\includegraphics[scale=0.4]{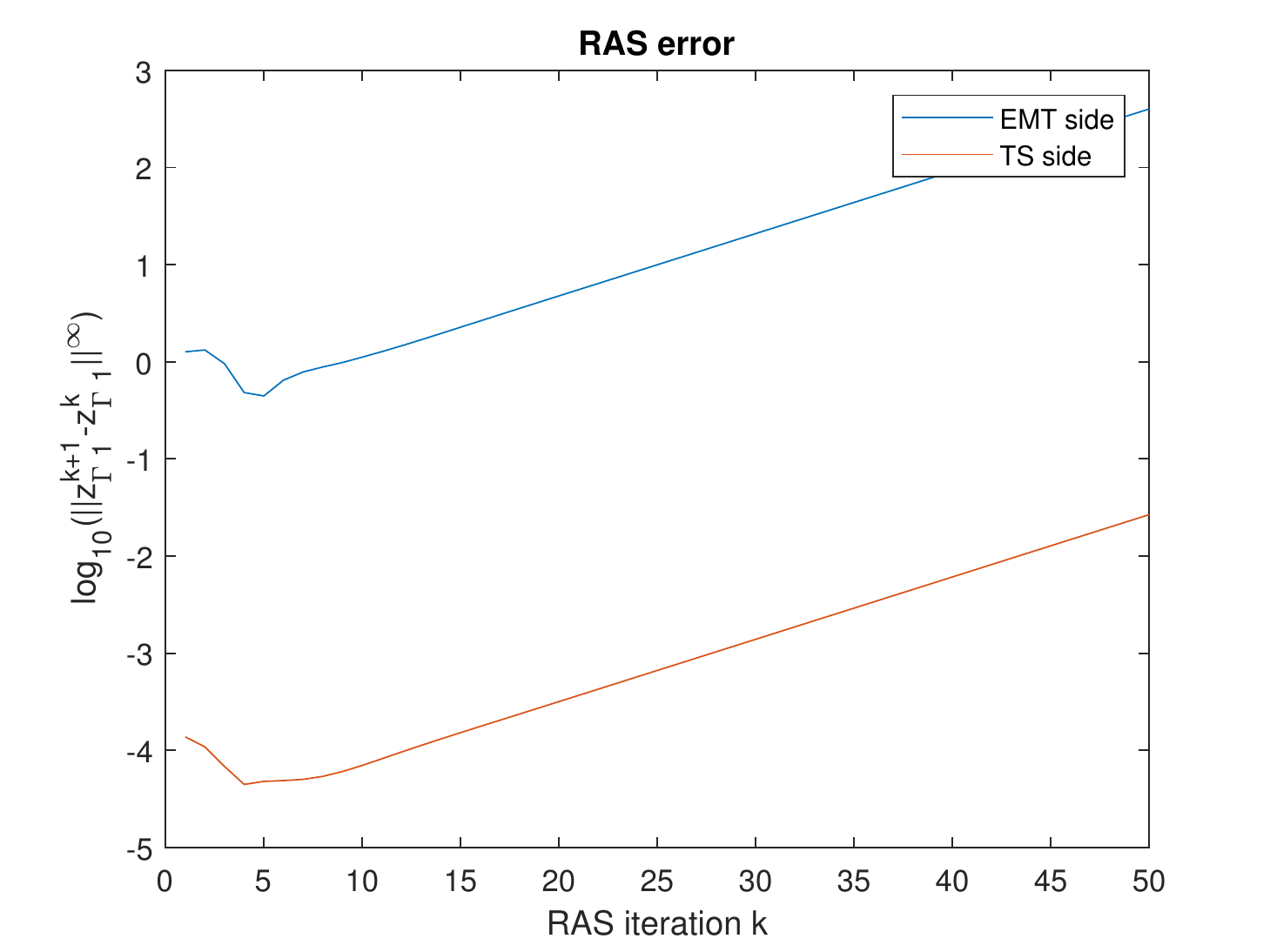}
\includegraphics[scale=0.4]{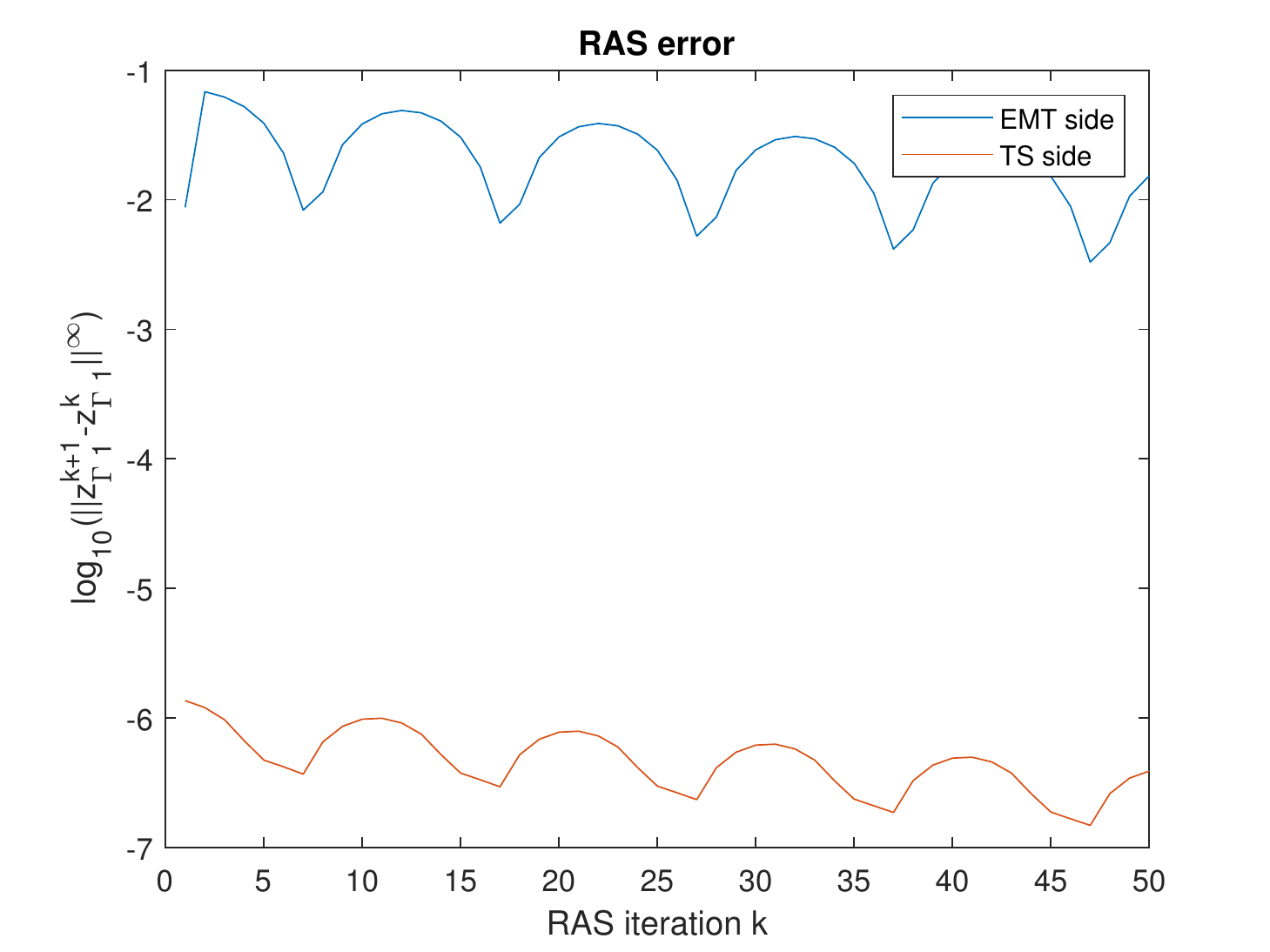}
\caption{ Convergence of the heterogeneous EMT-TS RAS error between two consecutive iterates with voltage source in the TS subsystem and  with, $R_1=R_2=7$,$L_1=L_2=0.5$ ,$C_1=C_2=1.10^{-6}$ on the left $\Delta t_{emt}=2.10^{-3}$ and on the right $\Delta t_{emt}=2.10^{-5}$.\label{PasDeTempsErrTS}}
\end{figure}
Figure \ref{PasDeTempsErrTS} gives the error between two consecutive heterogeneous RAS iterates with respect of the value of $\Delta t_{emt}$ for the circuit problem when the source is modeled in  TS. It exhibits that the method diverges with $\Delta t_{emt}=2.10^{-3}$ (left) while the method converges for $\Delta t_{emt}=2.10^{-5}$ (right).
\begin{figure}[h!]
\centering
\includegraphics[scale=0.4]{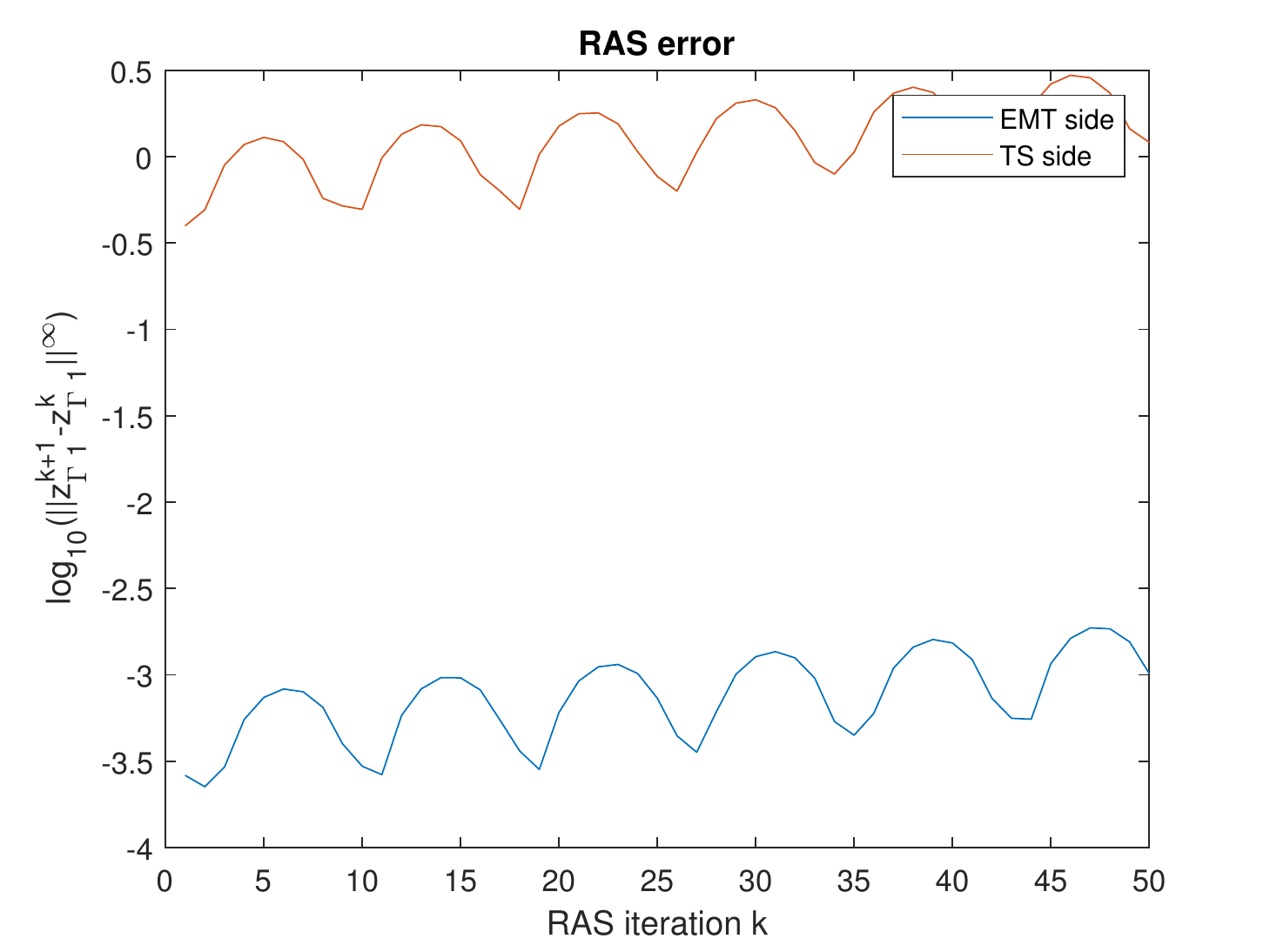}
\includegraphics[scale=0.4]{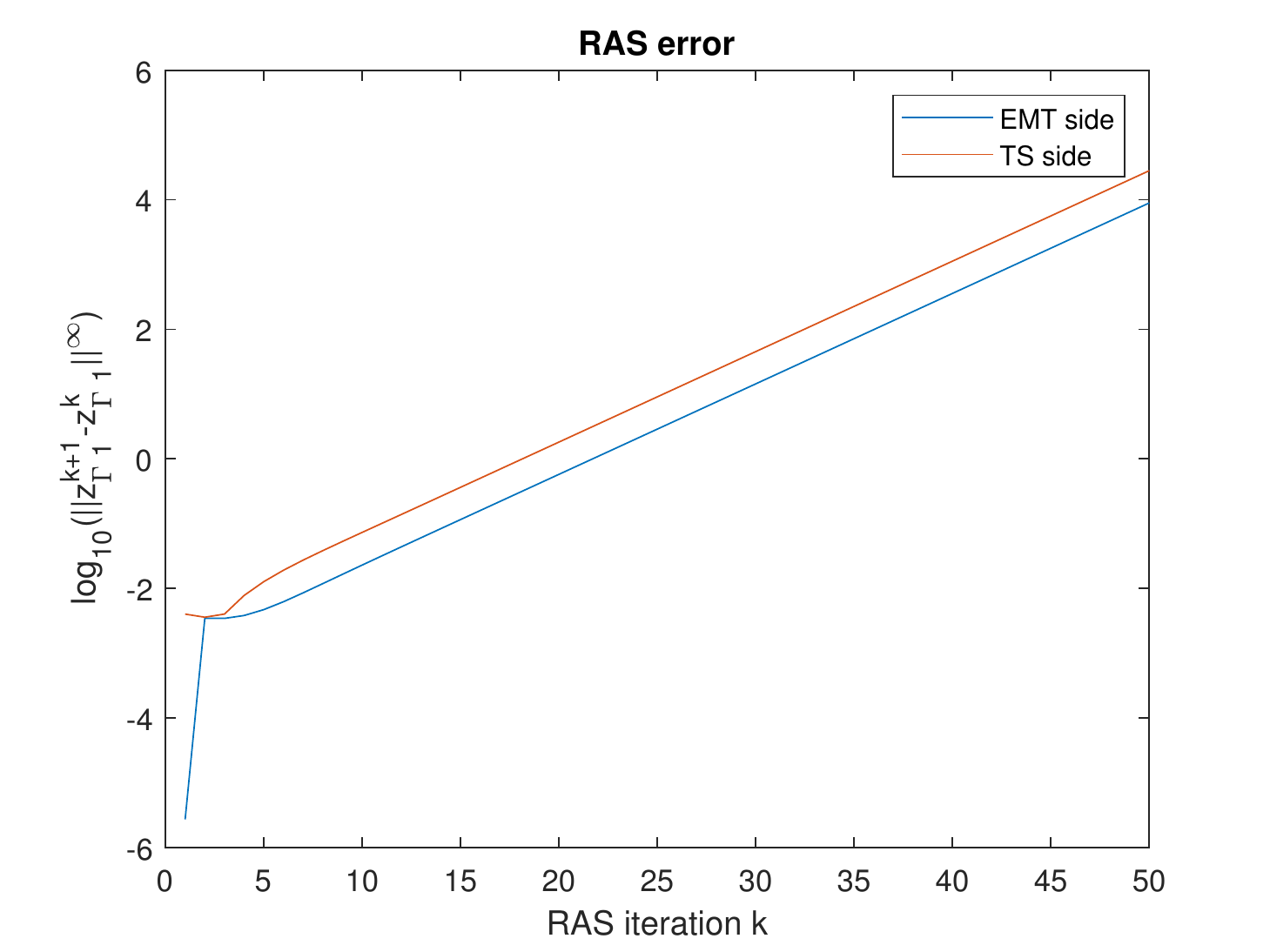}
\caption{Convergence of the heterogeneous EMT-TS RAS error between two consecutive iterates with voltage source in the EMT subsystem and  with, $R_1=R_2=7$,$L_1=L_2=0.5$ ,$C_1=C_2=1.10^{-6}$ on the left $\Delta t_{emt}=2.10^{-3}$ and on the right $\Delta t_{emt}=2.10^{-5}$.}
\label{PasDeTempsErrEMT}
\end{figure}
Figure \ref{PasDeTempsErrEMT} gives the error between two consecutive heterogeneous RAS iterates with respect of the value of $\Delta t_{emt}$ for the circuit problem when the source is modeled in EMT. It exhibits that the method diverges for both cases $\Delta t_{emt}=2.10^{-3}$ (left) and  $\Delta t_{emt}=2.10^{-5}$ (right). Indeed the more $\Delta t_{emt}$ decreases the more the error increases.

These results clearly demonstrate the need for the heterogeneous EMT-TS RAS method to have the Aitken's acceleration technique as in the homogeneous RAS case to be independent of the time step size choice.

 \subsubsection{Effect of the translation operators on the convergence}

Let us study the influence on the convergence of the percentage $\alpha$  in  the smoothing during the translation from TS to EMT.  We recall that it represents the number of EMT time steps on which the interpolation of TS interface values are performed.
\begin{figure}[h!]
\centering
\begin{minipage}{12cm}
\begin{minipage}{5.2cm}
\includegraphics[scale=0.4]{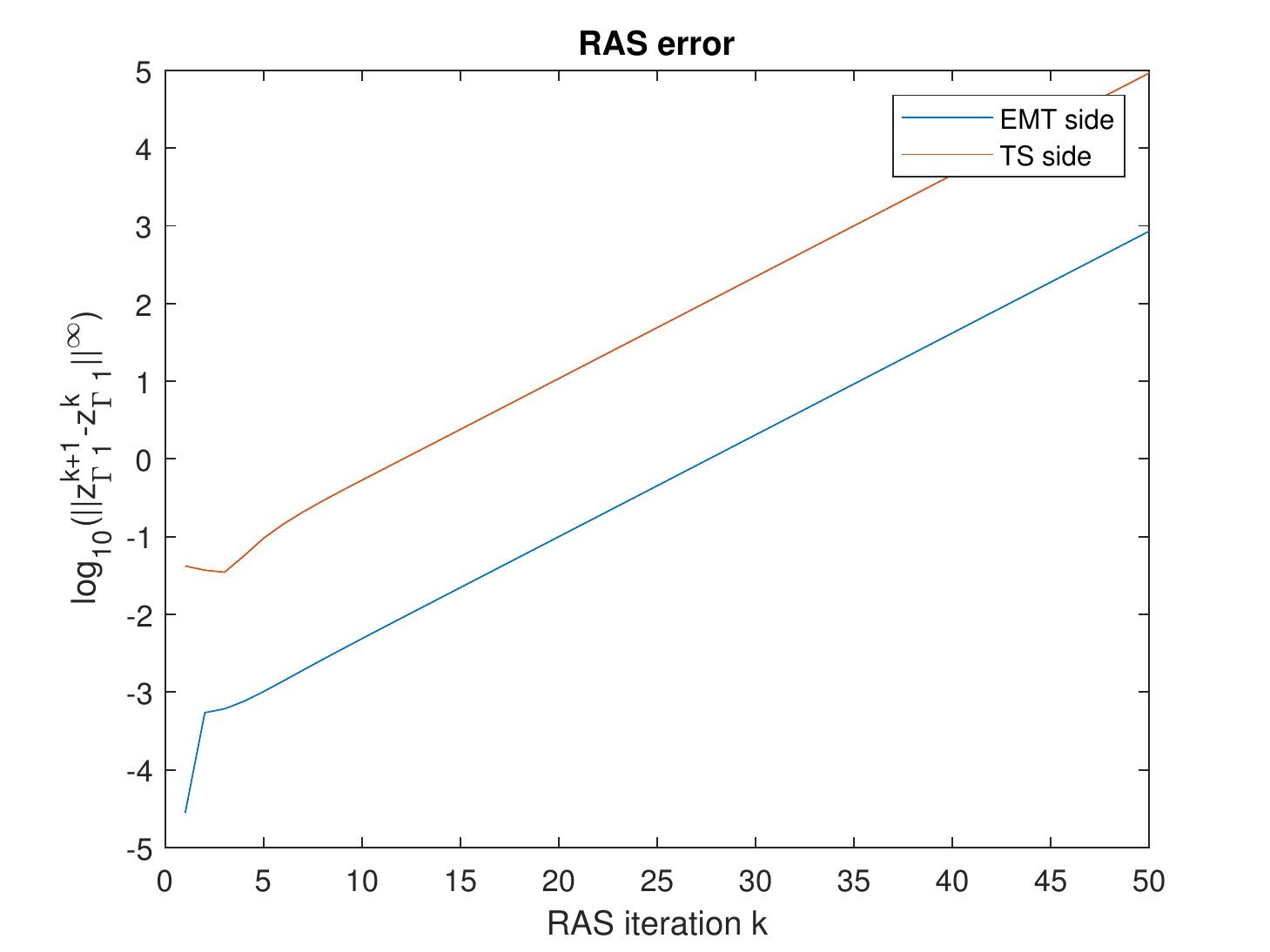}
\end{minipage}
\hfill
\begin{minipage}{5.2cm}
\includegraphics[scale=0.4]{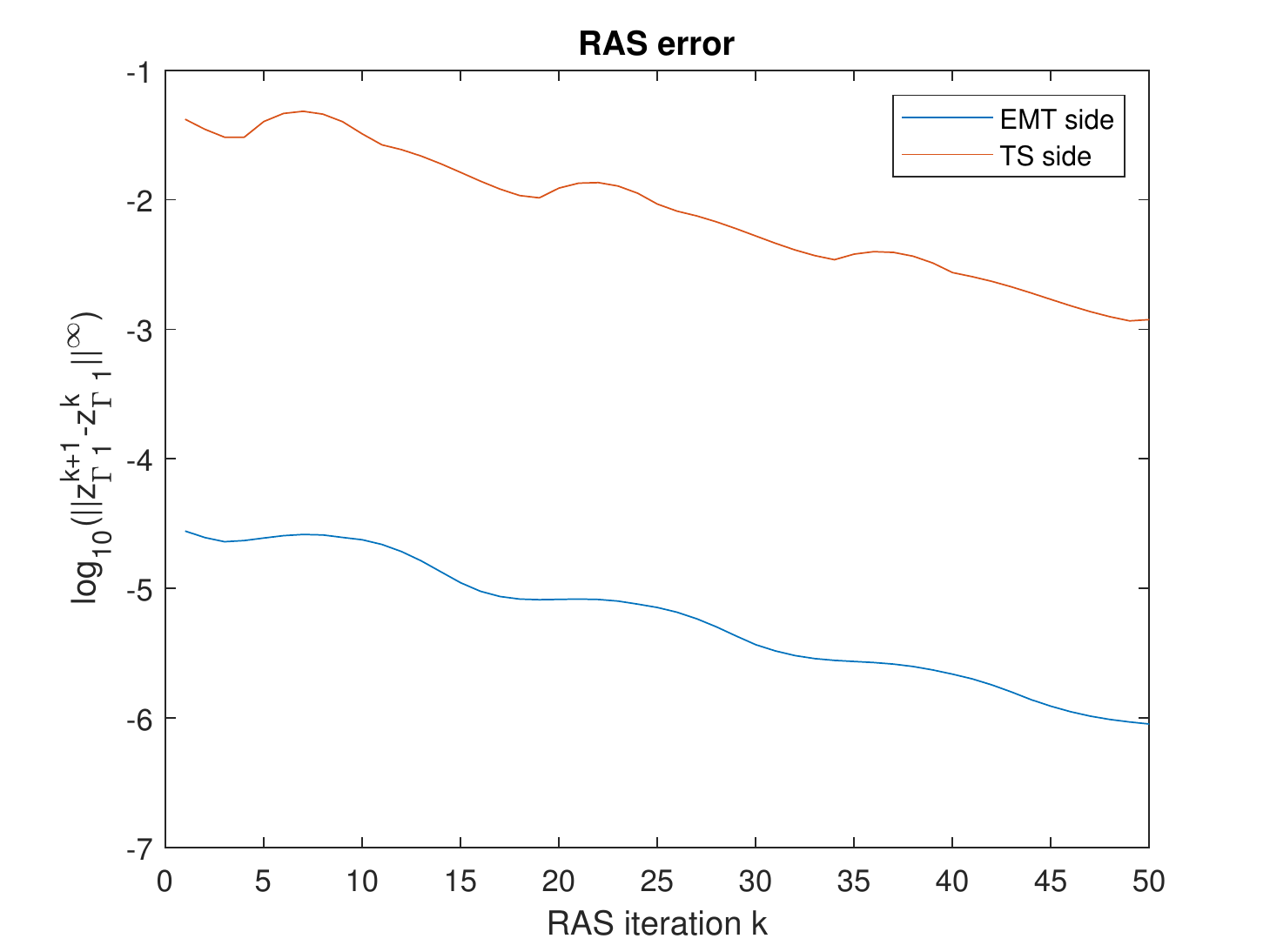}
\end{minipage}
\end{minipage} \\
\begin{minipage}{12cm}
\begin{minipage}{5.2cm}
\includegraphics[scale=0.4]{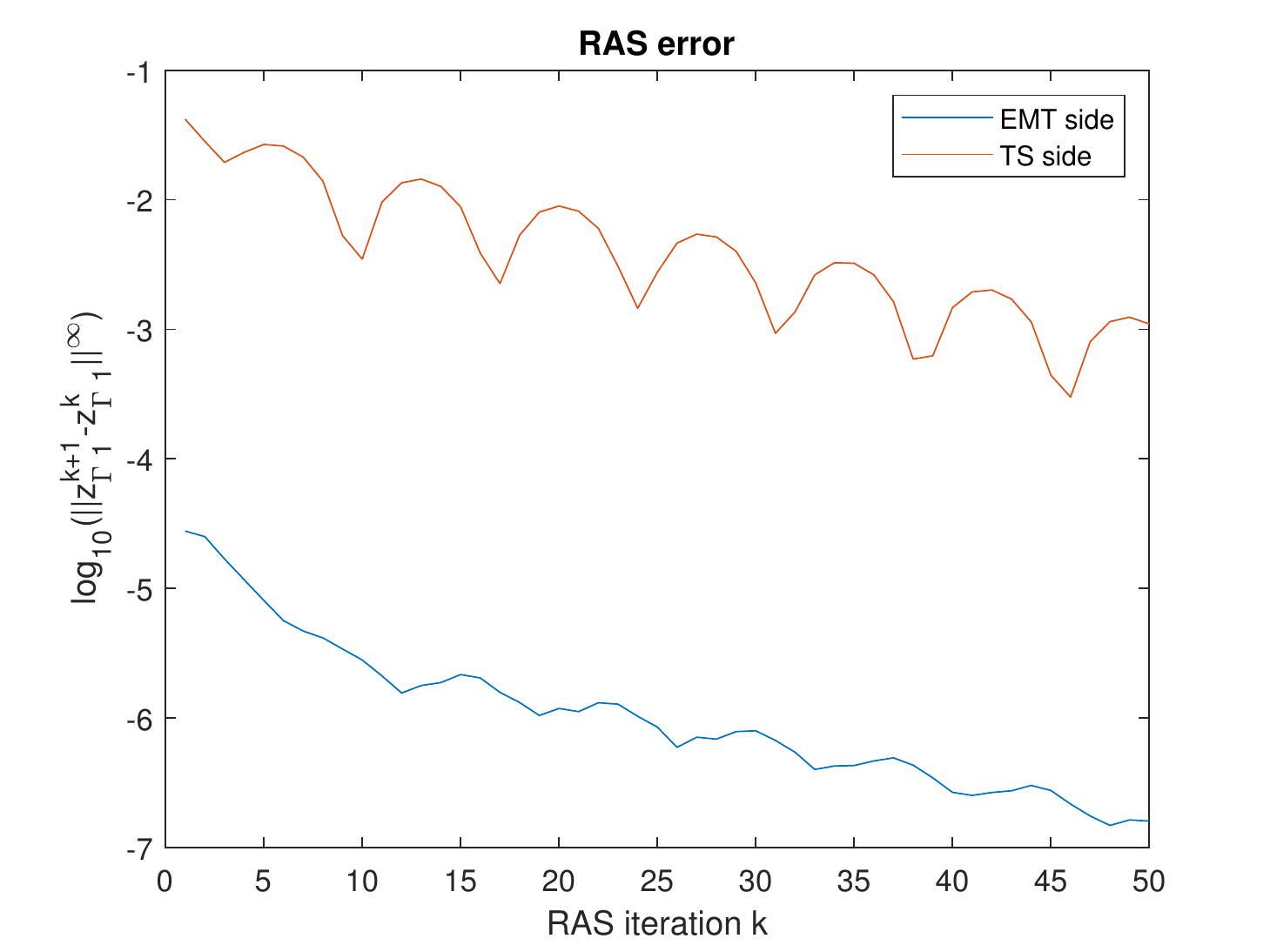}
\end{minipage}
\hfill
\begin{minipage}{5.2cm}
\includegraphics[scale=0.4]{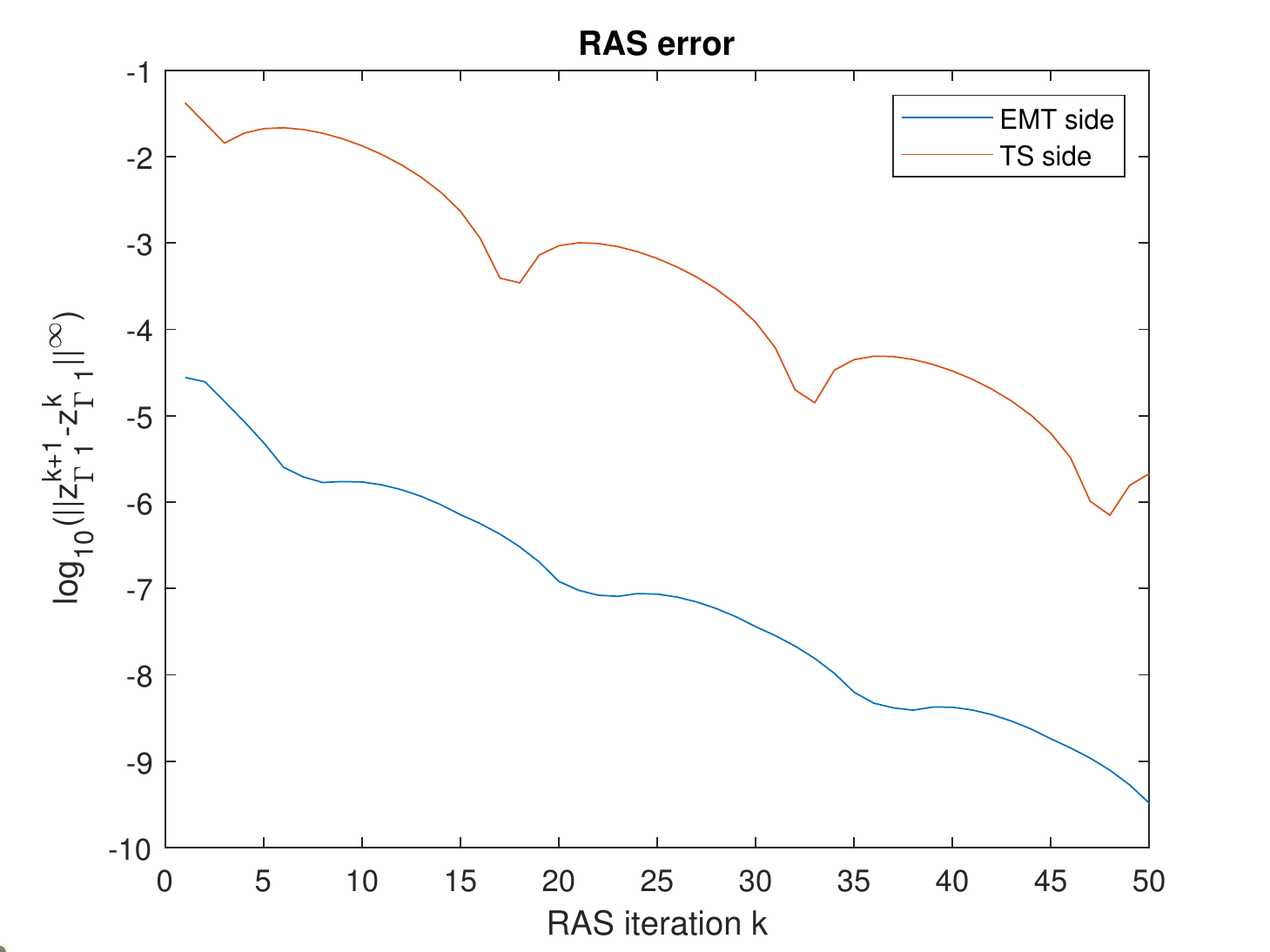}
\end{minipage}
\end{minipage} 
\caption{Convergence of the heterogeneous EMT-TS RAS error between two consecutive iterates with respect to the $\alpha$ parameter in the smoothing and  with voltage source in the EMT subsystem. $R_1=R_2=7$,$L_1=0.3$, $L_2=0.7$ ,$C_1=C_2=1.10^{-6}$. With top left $\alpha=0\%$, top right $\alpha=25\%$ , bottom left  $\alpha=50\%$ and bottom right  $\alpha=75\%$}
\label{ErrLissageEMT}
\end{figure}
Figure \ref{ErrLissageEMT} shows the error between two consecutive heterogeneous EMT-TS RAS iterates for  $\alpha=\left\{0\%, 25\%, 50\%, 75\%\right\}$. It exhibits that the $\alpha$ as a strong effect on the convergence as the method diverges when $\alpha=0\%$ (no smoothing) while it converges in other cases. It also shows that the larger is $\alpha$, the better is the convergence. 
  
\begin{figure}[h!]
\centering
\begin{minipage}{12cm}
\begin{minipage}{5.2cm}
\includegraphics[scale=0.4]{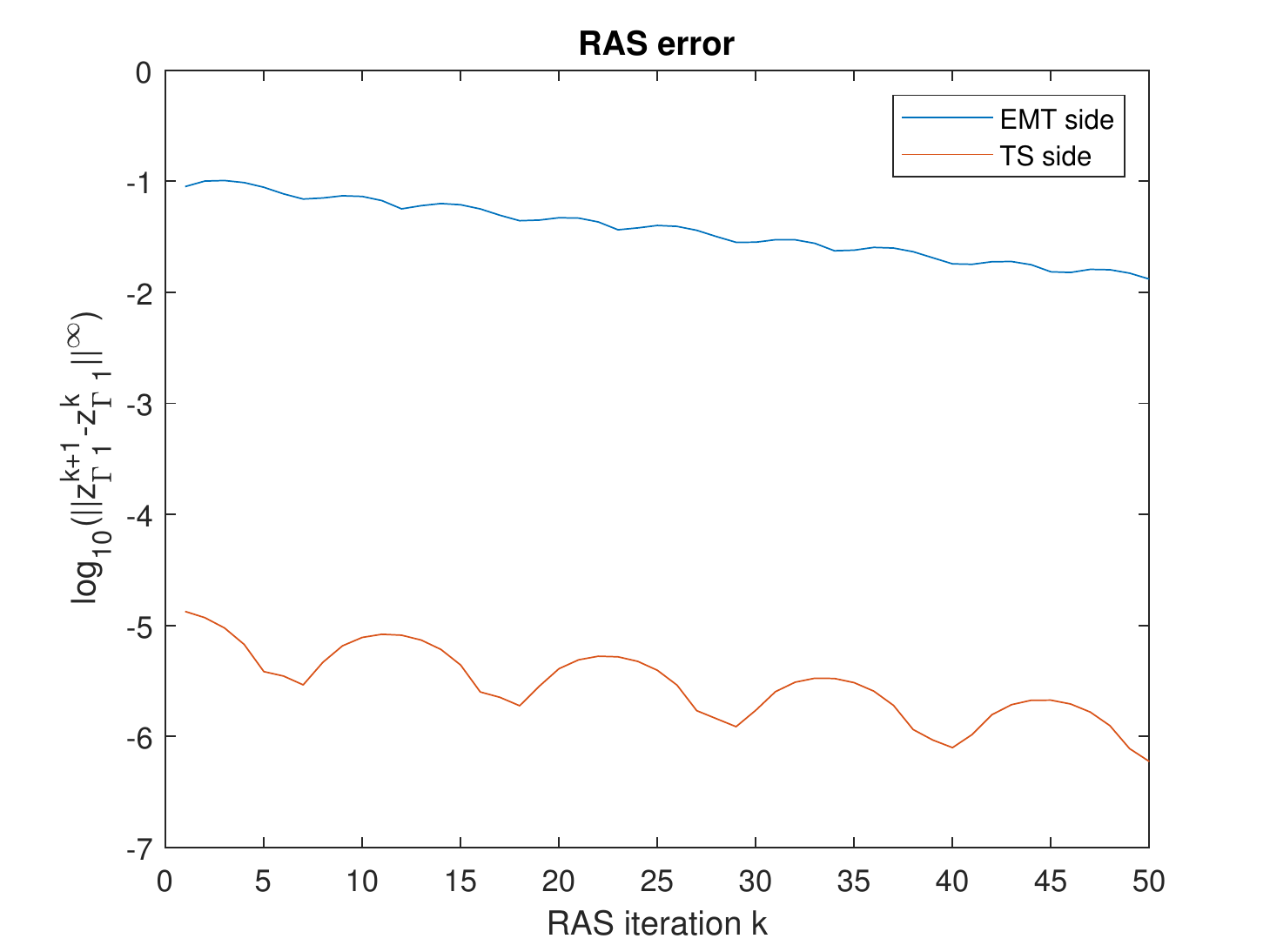}
\end{minipage}
\hfill
\begin{minipage}{5.2cm}
\includegraphics[scale=0.4]{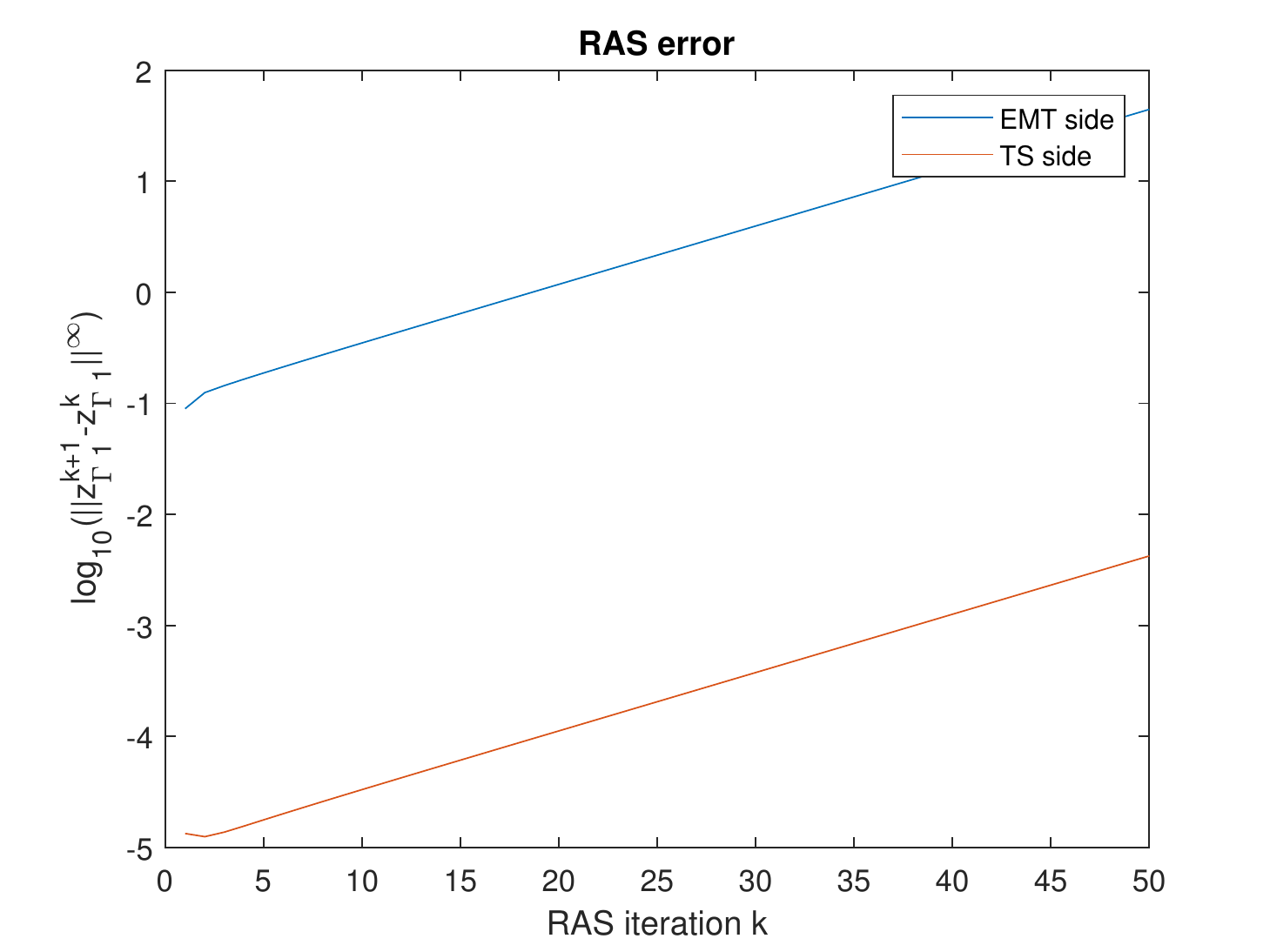}
\end{minipage}
\end{minipage} \\
\begin{minipage}{12cm}
\begin{minipage}{5.2cm}
\includegraphics[scale=0.4]{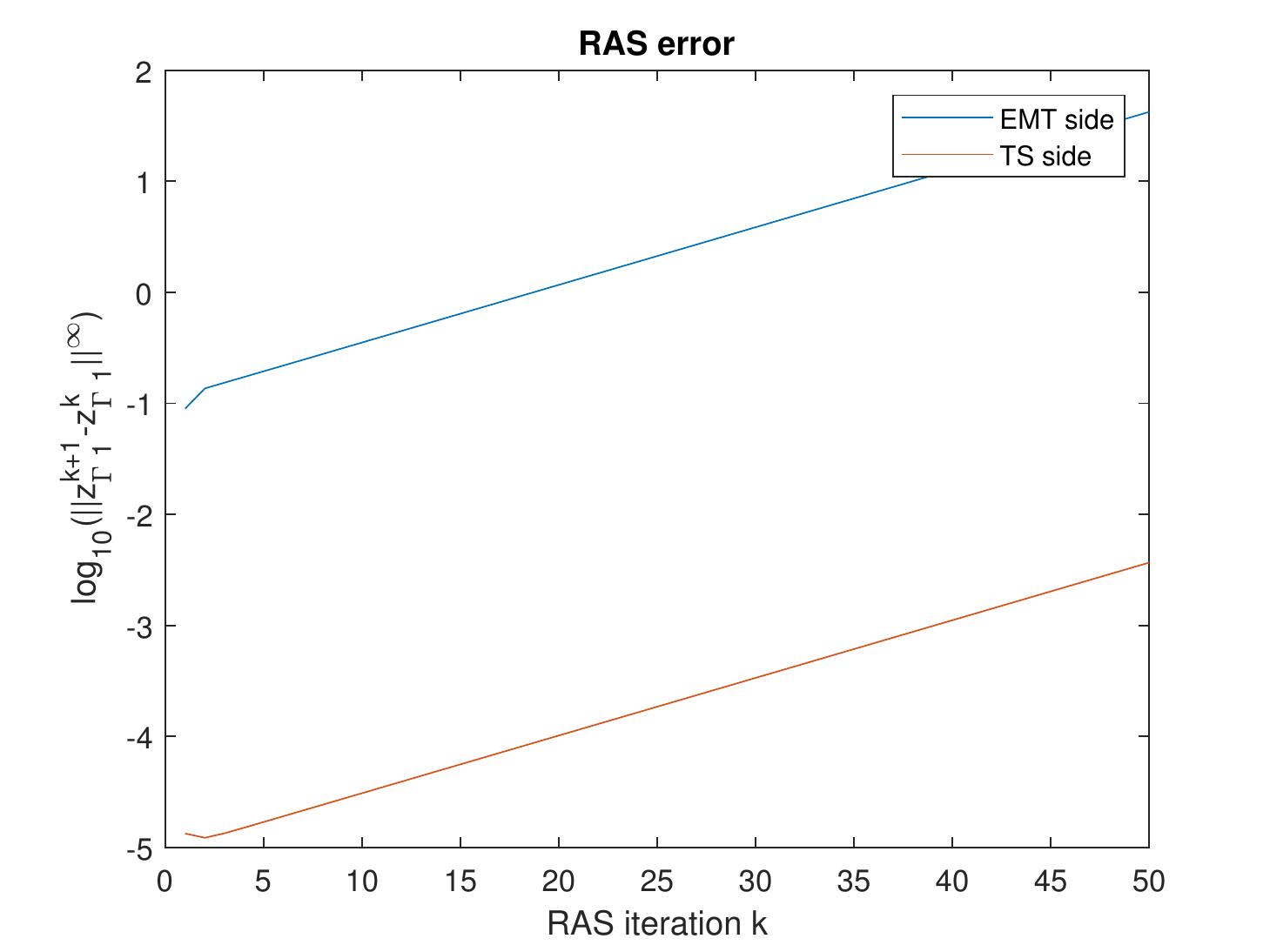}
\end{minipage}
\hfill
\begin{minipage}{5.2cm}
\includegraphics[scale=0.4]{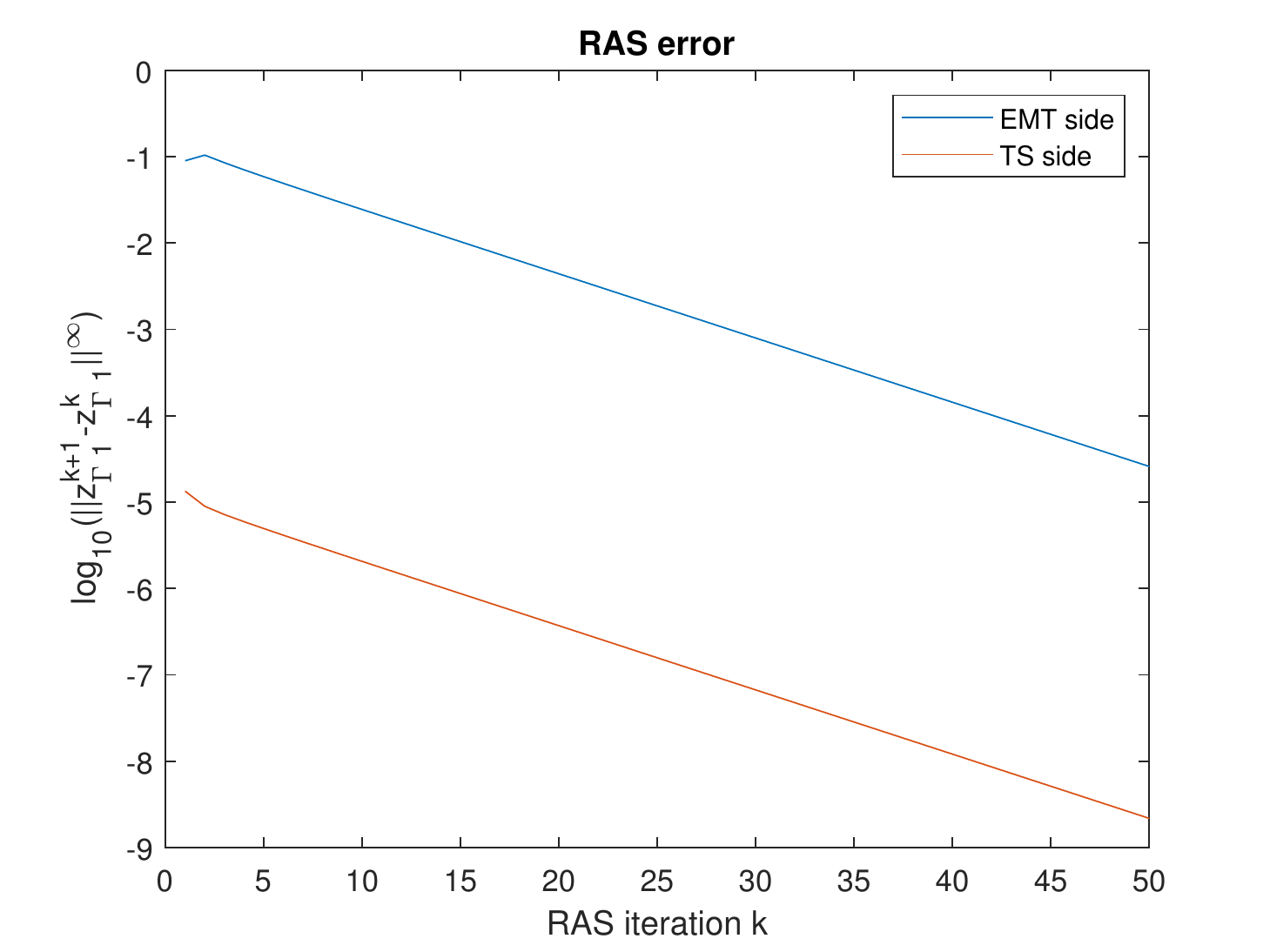}
\end{minipage}
\end{minipage} 
\caption{Convergence of the heterogeneous EMT-TS RAS error between two consecutive iterates with respect to the $\alpha$ parameter in the smoothing and  with voltage source in the TS subsystem. With, $R_1=R_2=7$,$L_1=0.3$, $L_2=0.7$ ,$C_1=C_2=1.10^{-6}$. With top left $\alpha=0$, top right $\alpha=25$ , bottom left  $\alpha=50$ and bottom right  $\alpha=75$}
\label{ErrLissageTS}
\end{figure}

As the increase in $\alpha$ has a beneficial effect on convergence in the case where the voltage source is in the EMT part, we investigated whether this will also be the case for the case where the voltage source is in the TS part which was convergent without smoothing. Figure \ref{ErrLissageTS} shows the error between two heterogeneous EMT-TS RAS iterates by varying $\alpha$ when the source is in the TS part.  It exhibits that the convergence  deteriorates when $\alpha$ becomes larger than $0\%$ until $\alpha$ reaches a threshold (about $50\%$). Once this threshold is exceeded, the error decreases with the increase of $\alpha$. We notice that the convergence improves faster than it has deteriorated. 

In conclusion, smoothing has a beneficial impact on the convergence of the heterogeneous EMT-TS RAS method as we found $\alpha$ values that make the method convergent regardless of the source voltage location. Nevertheless, the convergence behavior is not a monotonic function with respect to $\alpha$.

\subsubsection{Effect of the reduction of $\Delta_{ts}$ and moving history on the convergence}
Finally, let us study the effect of the reduction of the time steps $\Delta_{ts}$ with respect to the period $T$ and the moving history on the convergence.
\begin{figure}[h!]
\centering
\includegraphics[scale=0.4]{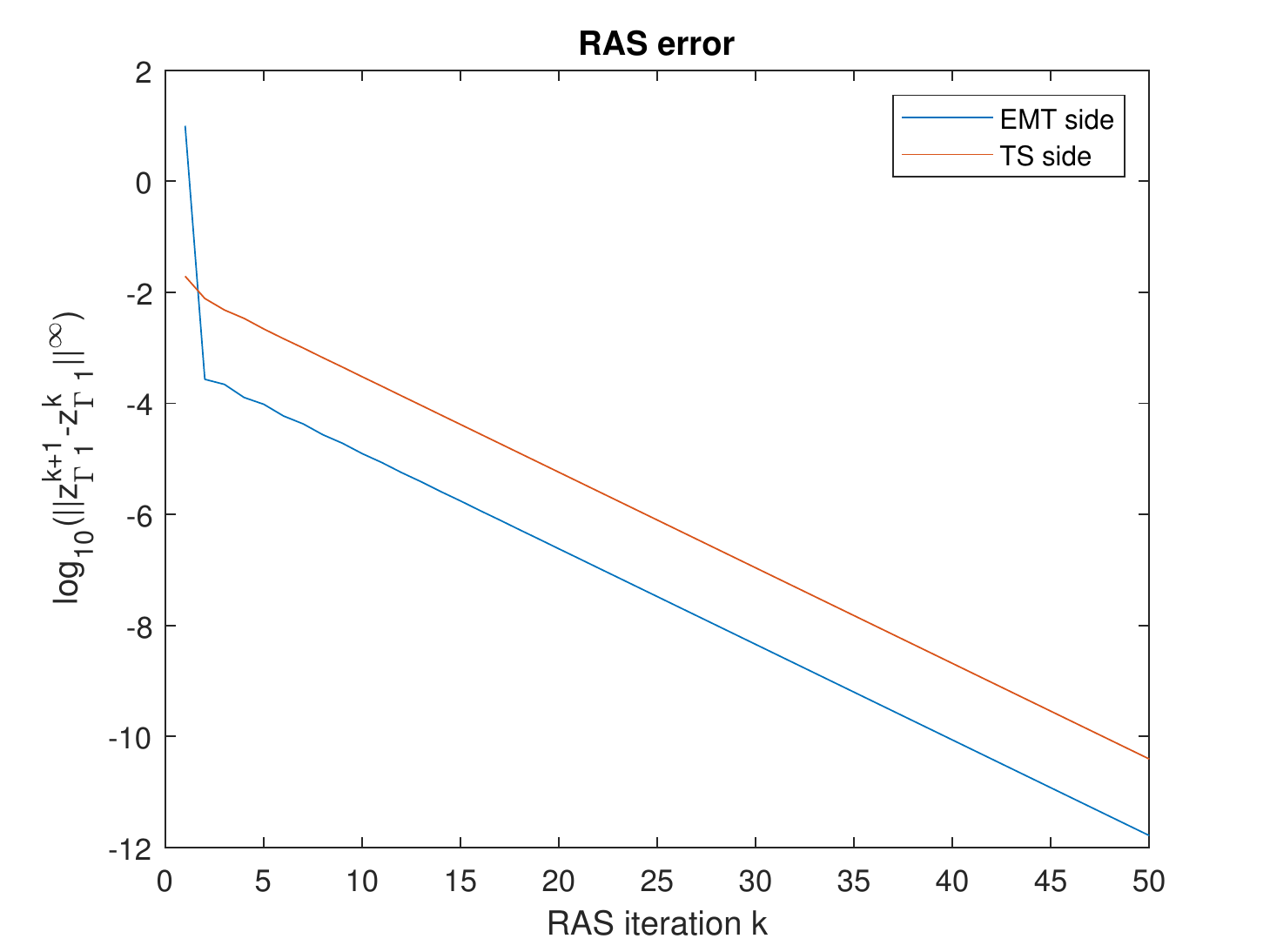}
\includegraphics[scale=0.4]{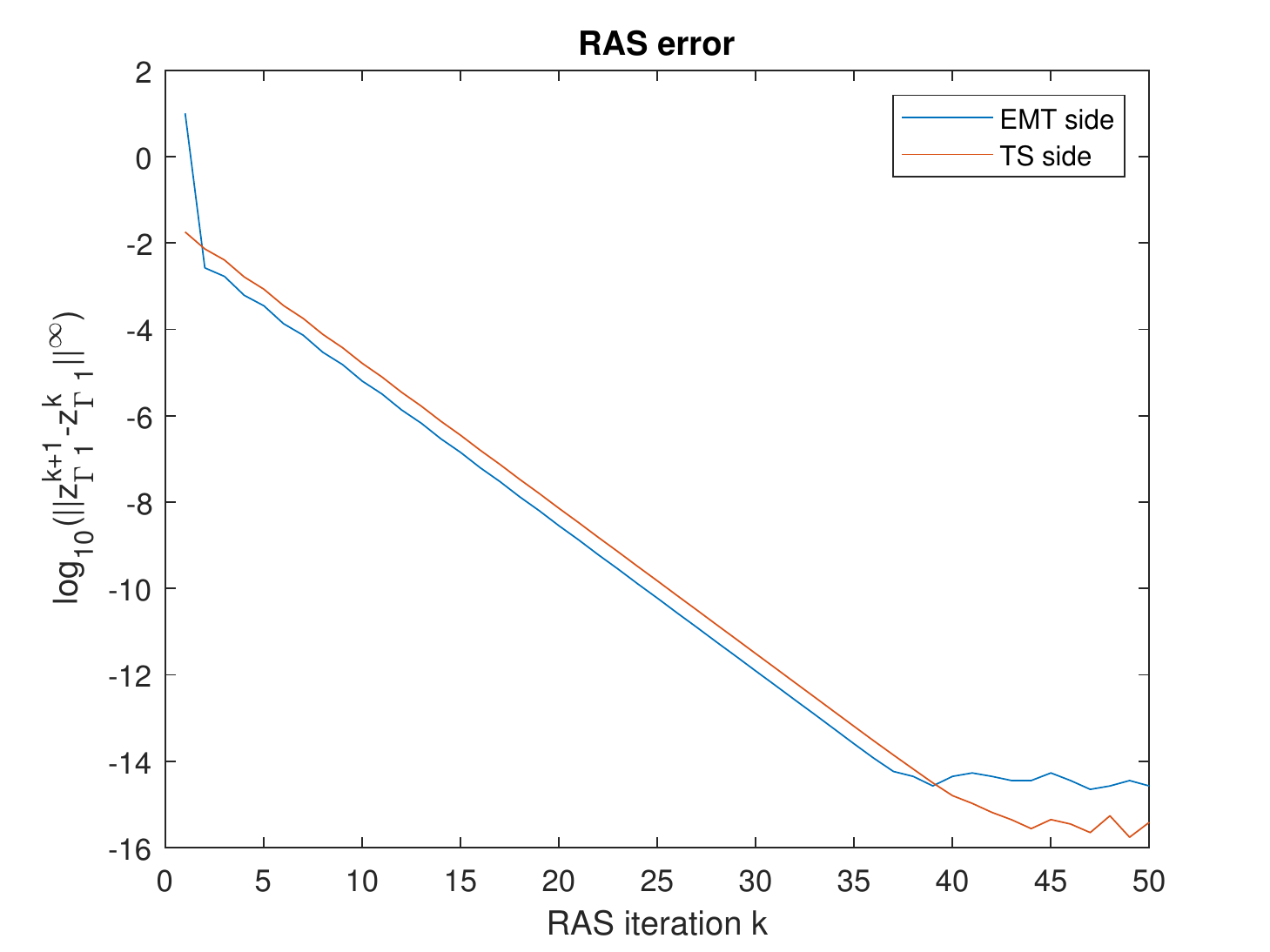}
\caption{Convergence of the heterogeneous EMT-TS RAS error between two consecutive iterates with respect to the $\Delta t_{ts}=10^{-2}$ less than a the period $T$ and rolling history,  with the source voltage in EMT and  $\alpha=0\%$ (left), $\alpha=25\%$ (right), $R_1=R_2=7$,$L_1=0.3$, $L_2=0.7$ ,$C_1=C_2=1.10^{-6}$.}
\label{FourierEMT}
\end{figure}

\begin{figure}[h!]
\centering
\includegraphics[scale=0.4]{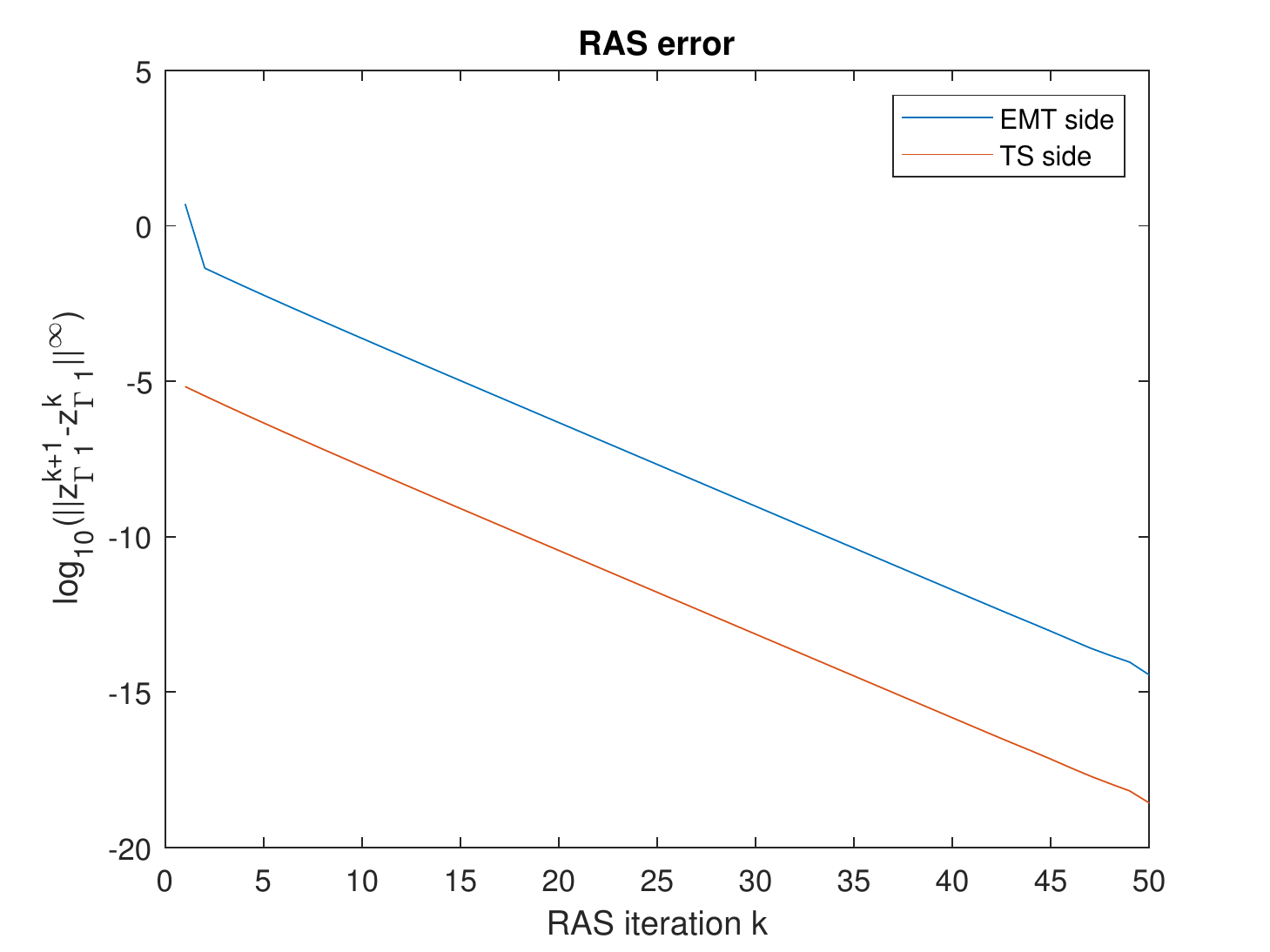}
\includegraphics[scale=0.4]{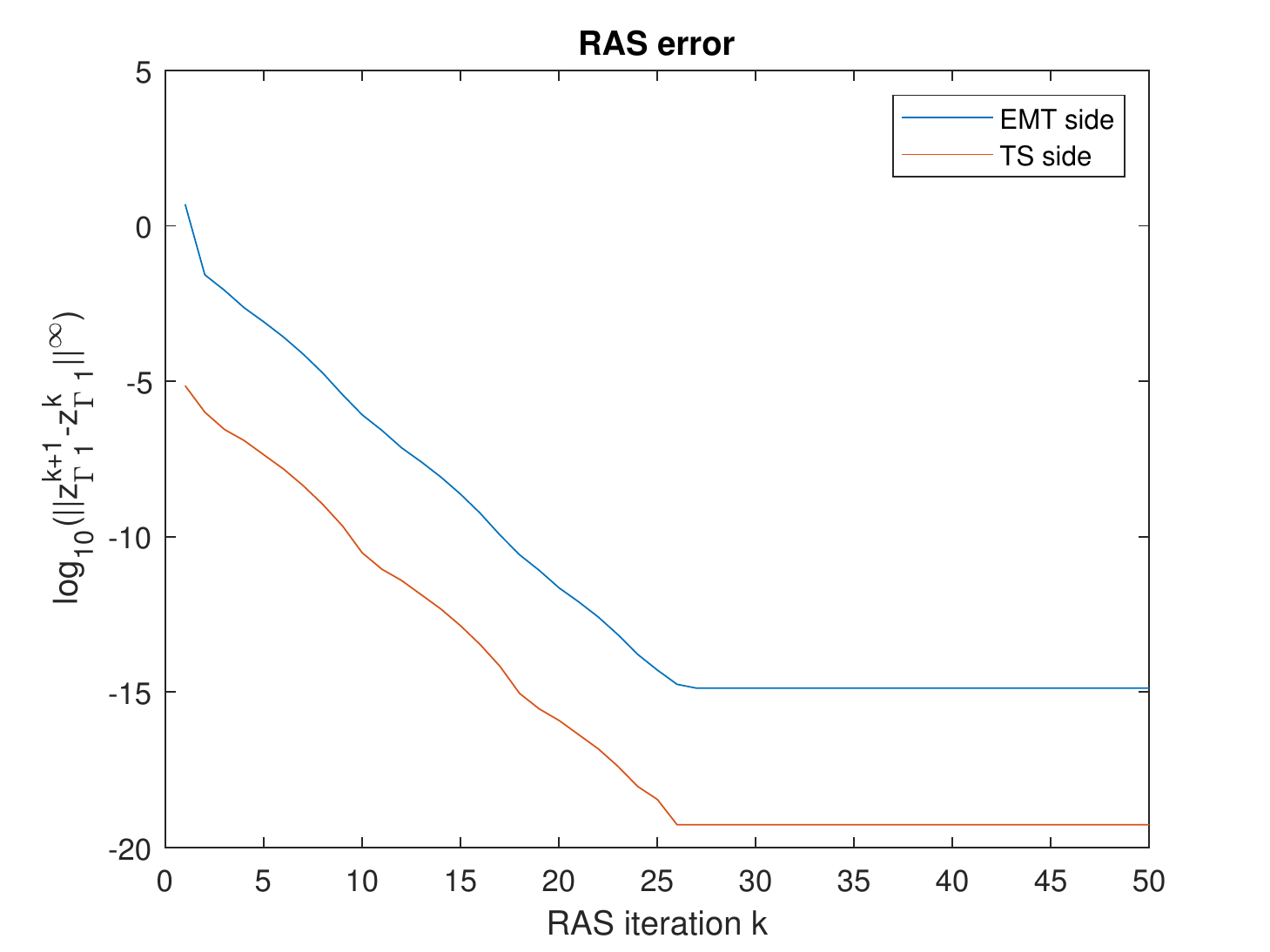}
\caption{Convergence of the heterogeneous EMT-TS RAS error between two consecutive iterates with respect to the $\Delta_{ts=\Delta t_{ts}=10^{-2}}$ less than a the period $T$ and rolling history,  with the source voltage in TS and  $\alpha=0\%$ (left), $\alpha=25\%$ (right), $R_1=R_2=7$,$L_1=0.3$, $L_2=0.7$ ,$C_1=C_2=1.10^{-6}$. }
\label{FourierTS}
\end{figure} 
Figure \ref{FourierEMT}  shows the error between two consecutive heterogeneous EMT-TS RAS iterates  for $\Delta_{ts}={1}/({2T})$ , $\alpha=\left\{0\%, 25\% \right\}$  and the voltage source in EMT. It exhibits that the decrease of $\Delta_{ts}={1}/({2T})$   and moving history lead to a convergent method in comparison of Figure \ref{ErrLissageEMT} (top left) with $\Delta_{ts}={1}/({T})$ which was divergent. The decrease of $\Delta_{ts}$ as also a beneficial impact on the convergence associated to $\alpha=25\%$ with an error nearby $10^{-12}$ for $50$ iterates for $\Delta_{ts}=10^{-2}$ instead of $10^{-6}$ for $\Delta_{ts}=2.~ 10^{-2}$ (Figure \ref{FourierEMT} (right) in comparison of Figure \ref{ErrLissageEMT} (top right)).

 Figure \ref{FourierTS} shows the error between two consecutive heterogeneous EMT-TS RAS iterates  for $\Delta_{ts}=\frac{1}{2T}$ , $\alpha=\left\{0\%, 25\% \right\}$  and the voltage source in TS. It also exhibits the beneficial effect of the decrease of $\Delta_{ts}=\frac{1}{2T}$ and moving history as the method converge for $\alpha=25\%$ that was not the case for $\Delta_{ts}=\frac{1}{T}$ in Figure  \ref{ErrLissageTS} (top right). The convergence is also greatly improved for the $\alpha=0\%$ case with an error of $10^{-15}$ in $40$ iterates  versus and error of $10^{-6}$ for $\Delta_{ts}=\frac{1}{T}$ in Figure \ref{ErrLissageTS} (top left).

In conclusion, decreasing the $\Delta t_{ts}$ time step associated with the moving history has greatly improved the convergence of the heterogeneous EMT-TS RAS, and made it less dependent on the source voltage location. Nevertheless, the decrease in $\Delta t_{ts}$ also decreases the ratio of $\Delta t_{ts}$ to $\Delta t_{emt}$ and thus the potential performance of the heterogeneous EMT-TS RAS.  
 
\subsubsection{Largest eigenvalue of the heterogeneous EMT-TS RAS error operator}
In order to confirm the convergence behavior of the heterogeneous EMT-TS RAS, we numerically calculated the error operator $P$ of the method and its largest eigenvalue.
 \begin{table}[h!]
 \centering
\begin{tabular}{|c|c|c|c|c|}
\hline
\multicolumn{5}{|c|}{\bf Source in the EMT subsystem}\\
\hline
\diagbox{$\alpha$}{$\Delta t_{ts}$}&\footnotesize$2. 10^{-2}$&\footnotesize$1.5~10^{-2}$&\footnotesize$1.~ 10^{-2}$&\footnotesize$2.~ 10^{-3}$\\
\hline
\footnotesize$ 0 \% $&\scriptsize$\lambda(P)=-1.3882$ &\scriptsize$\lambda(P)=-0.735 \pm 0.015i$&\scriptsize$\lambda(P)=-0.7068$&\scriptsize$\lambda(P)=-0.2537$\\
\hline
\footnotesize$ 25 \% $&\scriptsize$\lambda(P)=-0.884 \pm 0.407i$&\scriptsize$\lambda(P)=-0.7231$&\scriptsize$\lambda(P)=-0.6050 $&\scriptsize$\lambda(P)=-0.2538$\\
\hline
\footnotesize$ 50 \% $&\scriptsize$\lambda(P)=-0.876 \pm 0.427i$&\scriptsize$\lambda(P)=-0.6447$&\scriptsize$\lambda(P)=-0.4592$&\scriptsize$\lambda(P)=-0.2523$\\
\hline
\footnotesize$ 75 \% $&\scriptsize$\lambda(P)=-0.812 \pm 0.197i$&\scriptsize$\lambda(P)=-0.4938$&\scriptsize$\lambda(P)=0.3633$&\scriptsize$\lambda(P)=-0.2503$\\
\hline
\multicolumn{5}{|c|}{\bf Source in the TS subsystem}\\
\hline
\diagbox{$\alpha$}{$\Delta t_{ts}$}&\footnotesize$2.~ 10^{-2}$&\footnotesize$1.5 ~10^{-2}$&\footnotesize $1. ~10^{-2}$&\footnotesize $2.~ 10^{-3}$\\
\hline
\footnotesize$ 0 \% $&\scriptsize$\lambda(P)=-1.010 \pm 0.290i$ &\scriptsize$\lambda(P)=-1.0953$&\scriptsize$\lambda(P)=-0.5383$&\scriptsize$\lambda(P)=-0.1716$\\
\hline
\footnotesize$ 25 \% $&\scriptsize$\lambda(P)=-1.1430$&\scriptsize$\lambda(P)=-1.0417$&\scriptsize$\lambda(P)=-0.4478 $&\scriptsize$\lambda(P)=-0.1413$\\
\hline
\footnotesize$ 50 \% $&\scriptsize$\lambda(P)=-1.1789$&\scriptsize$\lambda(P)=-0.8259$&\scriptsize$\lambda(P)=-0.282 \pm0.053i$&\scriptsize$\lambda(P)=-0.1033$\\
\hline
\footnotesize$ 75 \% $&\scriptsize$\lambda(P)=-0.8930$&\scriptsize$\lambda(P)=-0.503 \pm 0.043i$&\scriptsize$\lambda(P)=-0.209 \pm 0.083i$&\scriptsize$\lambda(P)=-0.0831$\\
\hline
\end{tabular}
\caption{Largest eigenvalue of the error operator $P_{\Gamma_{ts}}$, with, $R_1=R_2=7$,$L_1=L_2=0.07$ ,$C_1=C_2=1.10^{-6}$ and $\Delta t_{emt} =2.10^{-4}$. At top with the source in the EMT subsystem and at bottom in the TS subsystem \label{EigenValuesHetero}}
\end{table}

Table \ref{EigenValuesHetero} gives the combined impact of smoothing and moving history on the convergence of the method by showing the value of the largest eigenvalue of the error operator $P$ for source voltage in EMT and TS. It exhibits that except when we are in the particular case where the source of tension is in TS and that $\Delta_{ts}= 1/T$, the increase of $\alpha$ always improves the convergence. The decrease $\Delta_{ts}$ possible thanks to the moving history always improves the convergence too. The combination of the two makes it possible to obtain fairly rapid convergences but still to much from an operational point of view. This is why we need to accelerate its convergence.

\subsubsection{Aitken's acceleration of the heterogeneous EMT-TS RAS convergence}
For all the studies of convergences of the previous sub-section, we see that the convergence/divergence is purely linear, we will therefore use the method of acceleration of convergence of Aitken to obtain the solution after $n_\Gamma+1$ iterations. 
 \begin{figure}[h!]
 \centering
\includegraphics[scale=0.55]{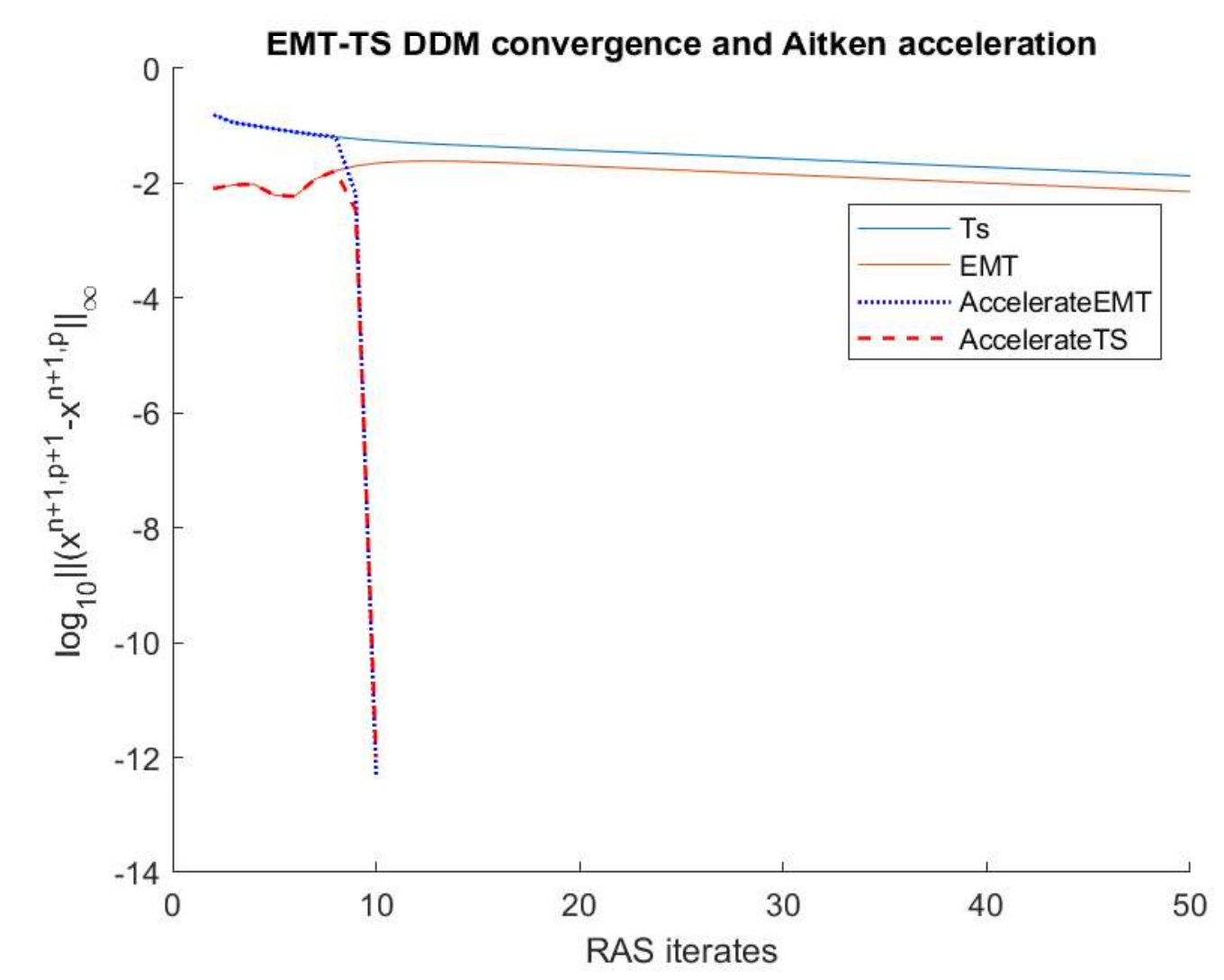}
\caption{ Heterogeneous EMT-TS  RAS convergence error for each subdomain, for the time step $t = 0.02$ and its Aitken's acceleration applied on the TS partition interface  with $\Delta t_{ts}=2.10^{-2}$ and $\Delta t_{emt}=2.10^{-4}$) and with parameters $L_1=0.07, C1=1.10^{-6} R1=7, L2=0.07, C2=1.10^{-6}, R2=7, Zs=0.000001$, with the voltage source in the TS part}
\label{shourick_contrib_Fig3c}
\end{figure}
Figure \ref{shourick_contrib_Fig3c} gives the $log_{10}$ of the error between two consecutive RAS iterations at time $t=0.02$. It shows a linear convergence behavior and can therefore be accelerated by the Aitken acceleration of the convergence technique after $9$ iterates needed to numerically construct the error operator $P_\Gamma$. In this case the method is convergent and it is shown that after the local resolutions, the method accelerated has converged on each subdomain. 

\begin{figure}[h!]
\centering
\includegraphics[scale=0.6]{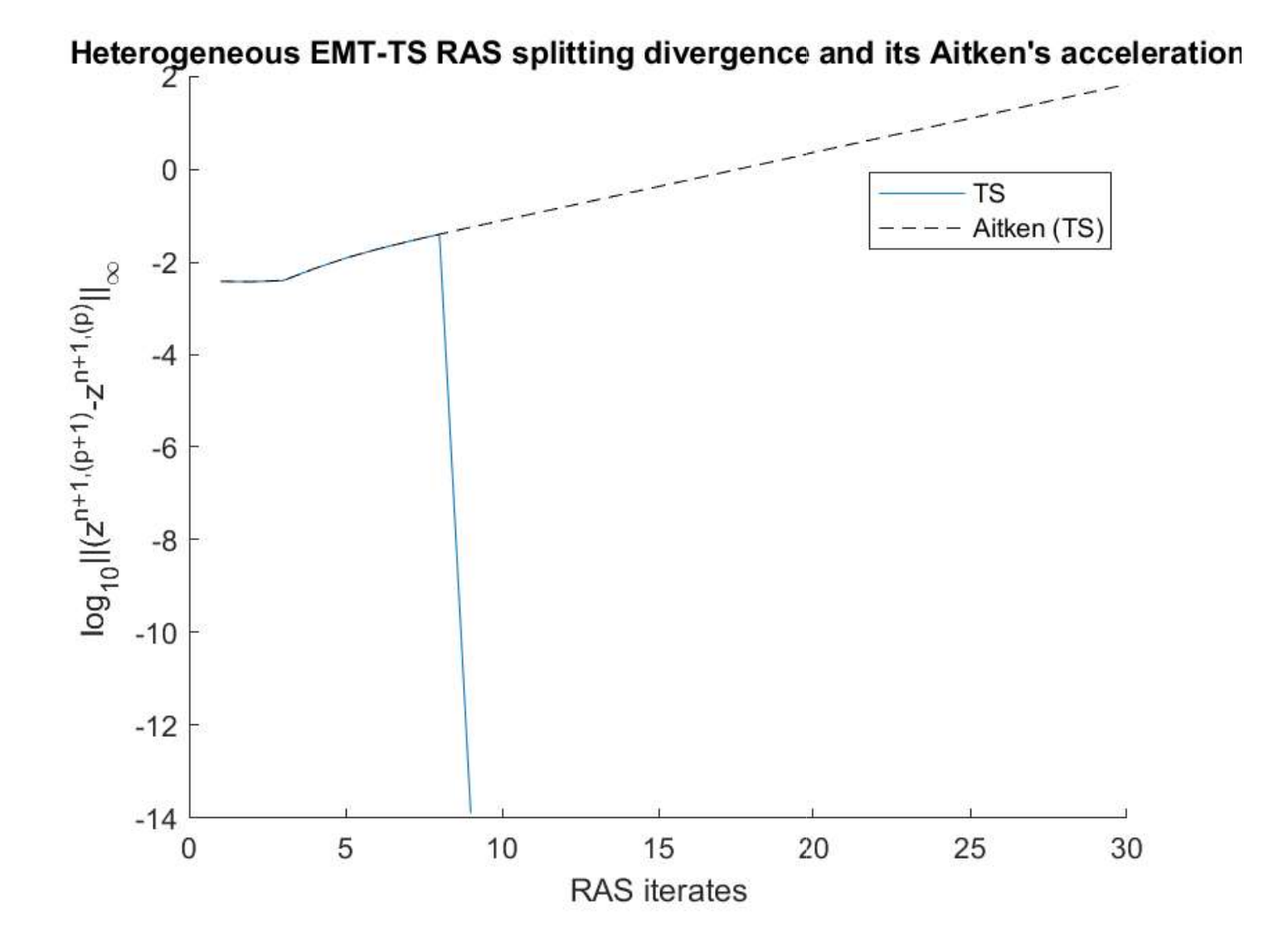}
\hspace*{-2.5cm} \includegraphics[scale=0.37]{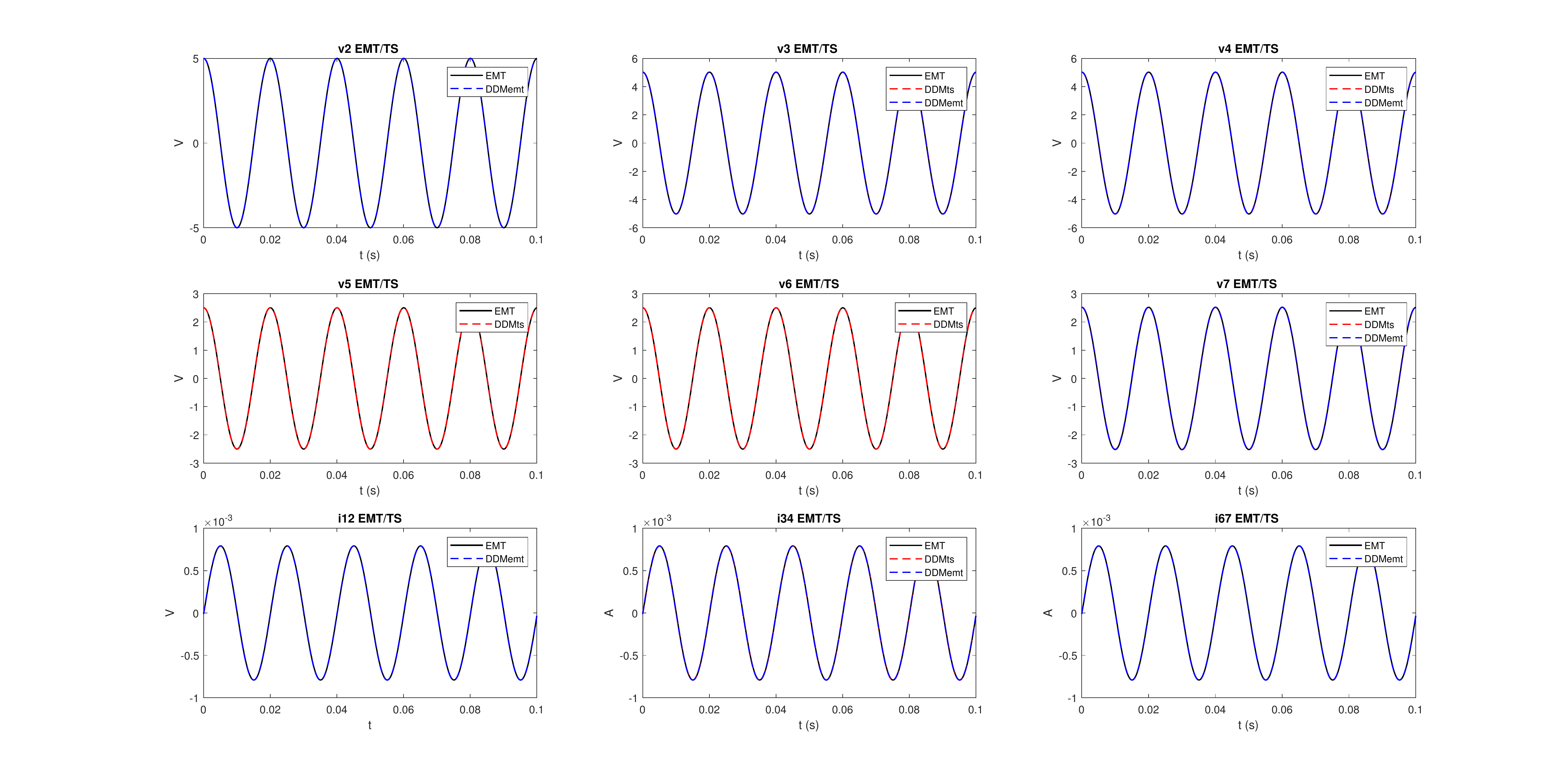}
\caption{ Heterogeneous EMT-TS  RAS convergence error for the TS boundary for the time step $t = 0.02$ (top) and solution comparison with the monolithic EMT case on a time intervalle (bottom)  and its Aitken's acceleration applied on the TS partition interface  with $\Delta t_{ts}=2.10^{-3}$ and $\Delta t_{emt}=2.10^{-5}$) and with parameters $L_1=0.07, C1=1.10^{-6}, R1=7, L2=0.07, C2=1.10^{-6}, R2=7, Zs=0.000001$, with the voltage source in the EMT part }
\label{shourick_contrib_Fig3b}
\end{figure}

Figure \ref{shourick_contrib_Fig3b} (top) shows the $log_{10}$ of the error between two successive heterogeneous EMT-TS RAS iterates with and without the use of the Aitken acceleration method with the voltage source in the EMT part with $\Delta t_{ts}=2.10^{-3}$ and $\Delta t_{emt}=2.10^{-5}$). It exhibits that the method diverges but as the divergence is purely linear the Aitken's acceleration of the convergence technique successes to retrieve the true solution after $9$ iterates as the acceleration is performed on the interface solution of the TS part.

Figure \ref{shourick_contrib_Fig3b} (bottom) shows the solution after using the Aitken acceleration technique in a divergent case. The result is compared to the monolithic EMT case showing that the true solution is obtained.

 \subsection{Qualitative results on the heterogeneous EMT-TS RAS}

 The numerical tests that follow focus on the qualitative advantage of the EMT-TS model over the TS model. The problem is that of the circuit in figure \ref{RLCCut}.
 \subsubsection{The advantage of the EMT part}
 
The purpose of having a part of the circuit simulated in the EMT modelling is to capture events that are not visible with the dynamic phasor modelling. In order to verify that the EMT-TS modeling can capture this advantage, we will create a disturbance that lasts less than a TS time step and therefore cannot be seen by the TS modeling.

\begin{figure}[h!]
  \centering
  \begin{minipage}{14cm}
  \begin{minipage}{6.8cm}
\includegraphics[scale=0.5]{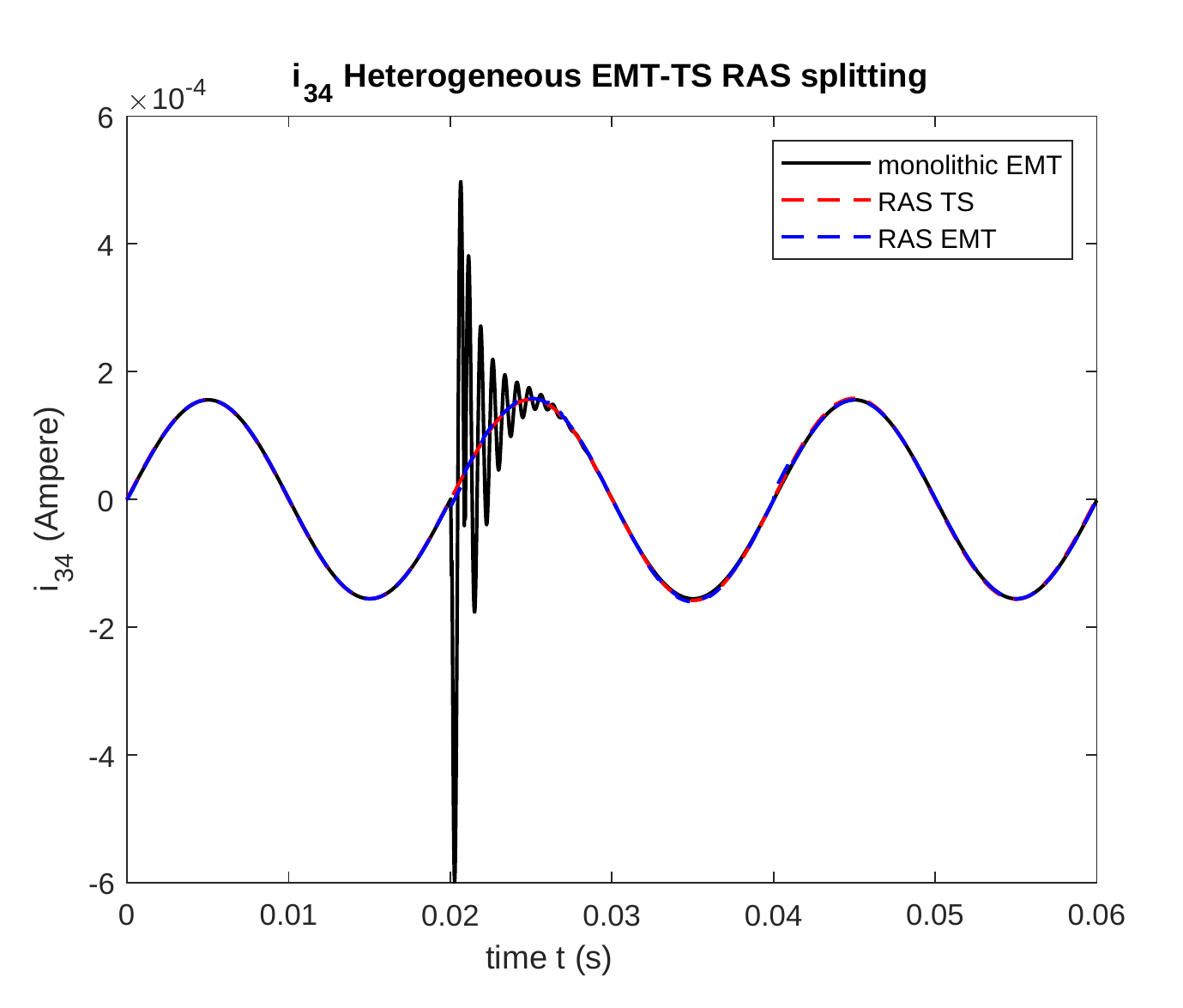}
\end{minipage}
\hfill
  \begin{minipage}{6.8cm}
\includegraphics[scale=0.5]{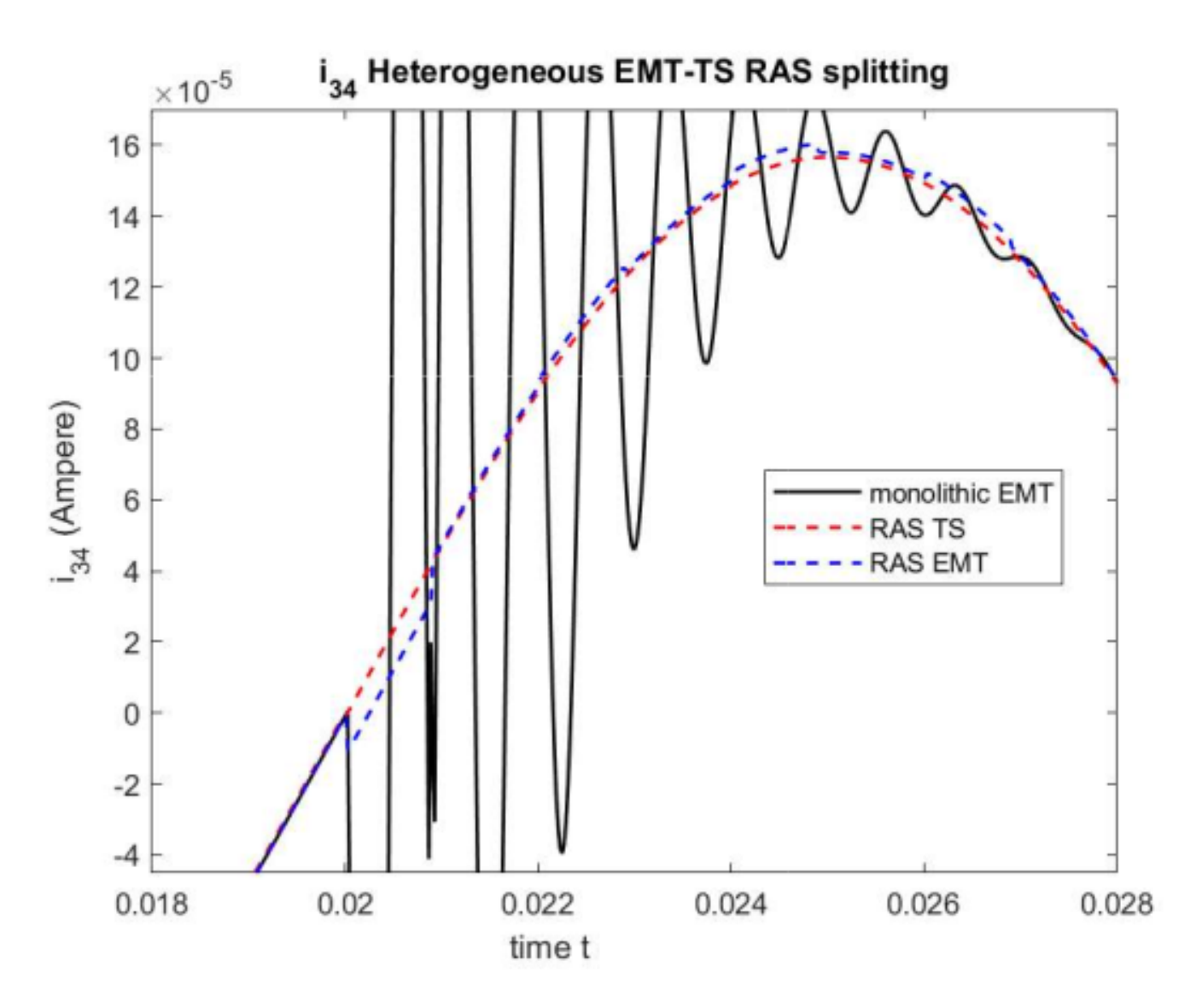}
\end{minipage}
\end{minipage}
   \begin{minipage}{14cm}
  \begin{minipage}{6.8cm}
\includegraphics[scale=0.5]{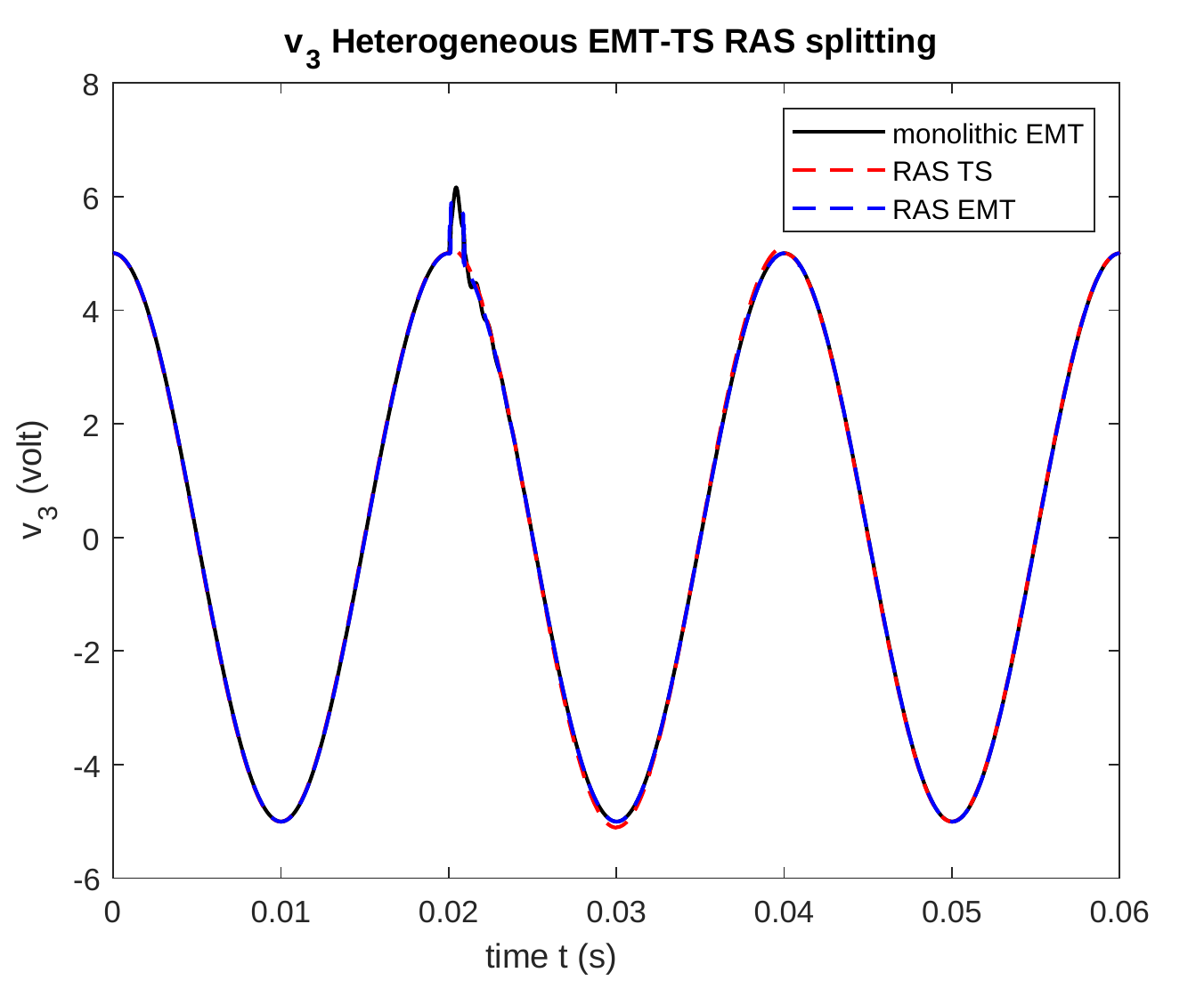}
\end{minipage}
\hfill
  \begin{minipage}{6.8cm}
\includegraphics[scale=0.5]{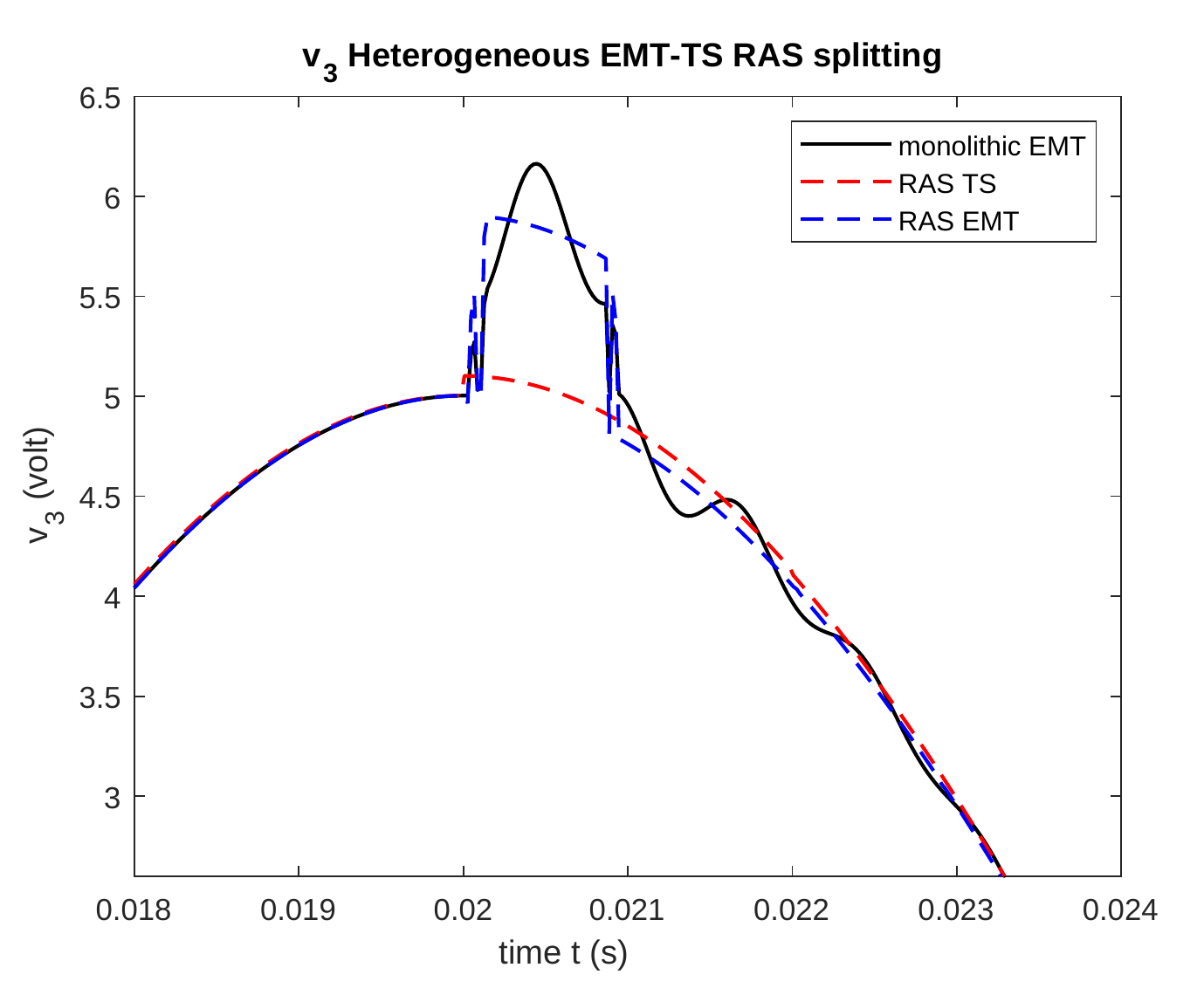}
\end{minipage}
\end{minipage}
\caption{ Comparison of the behavior, with respect to time, of the variables $i_{34}$ (top left) and $v_3$ (bottom left) (the figures on the right are their zoom on the disturbances) computed using the heterogeneous EMT-TS  RAS splitting with the Aitken's technique for accelerating convergence ($\Delta t_{ts}=2.10^{-3}$ and $\Delta t_{emt}=2.10^{-5}$), the reference is the monolithic EMT.  An amplitude perturbation on the voltage source starting  at $t=0.02$s and ending at $t=0.021$s, therefor lasting less than one $\Delta t_{ts}$ is applied. Parameters are $L_1=0.07, C_1=1.10^{-5}, R_1=7, L_2=0.07, C_2=1.10^{-7}, R_2=7, Zs=0.000001$.}\label{shourick_contrib_Fig3}
\end{figure}

Figure \ref{shourick_contrib_Fig3} compares the EMT monolithic reference values for the variables $v_3$ and $i_{34}$ with the EMT-TS heterogeneous RAS splitting where a perturbation on the source voltage that starts at $t=0.02$s and ends at $t=0.021$s is applied. It exhibits that the heterogeneous EMT-TS RAS succeeds in capturing part of the perturbation on the $v_3$. It shows a good agreement between the monolithic and the  heterogeneous EMT-TS RAS for the variable $v_3$. The variable $i_{34}$ in the EMT DDM part captures certain oscillations due to the perturbation. These results show that  heterogeneous EMT-TS RAS  can capture disturbances that last less than one TS time step and therefore would not have been captured by a monolithic TS model. 

 \subsubsection{The impact of the cutting and passing of information}
 \begin{figure}[h!]
 \hspace{-1cm}
   \begin{minipage}{13.5cm}
  \begin{minipage}{6.5cm}
\includegraphics[scale=0.55]{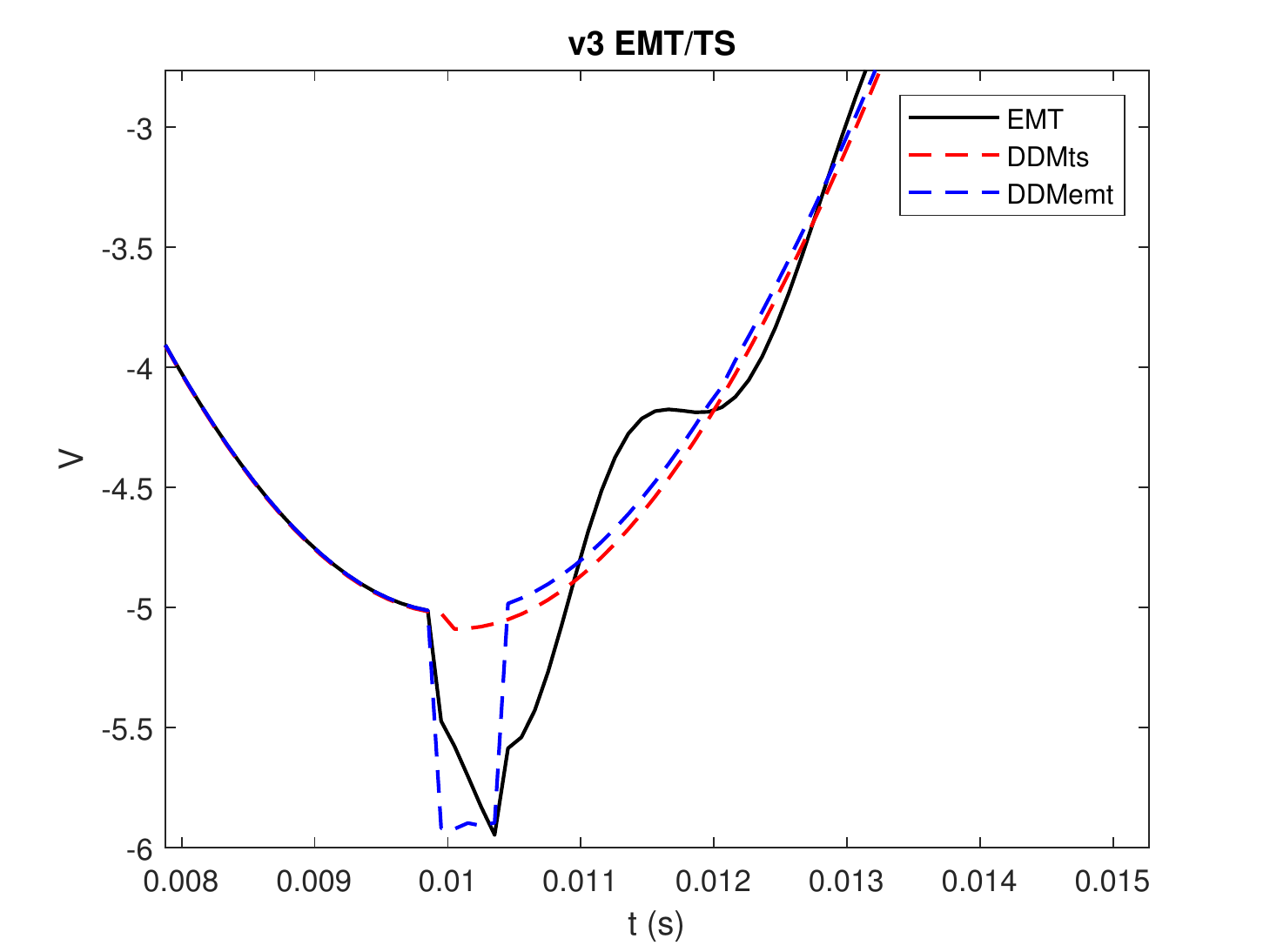}
\end{minipage}
\hfill
  \begin{minipage}{6.cm}
\includegraphics[scale=0.42]{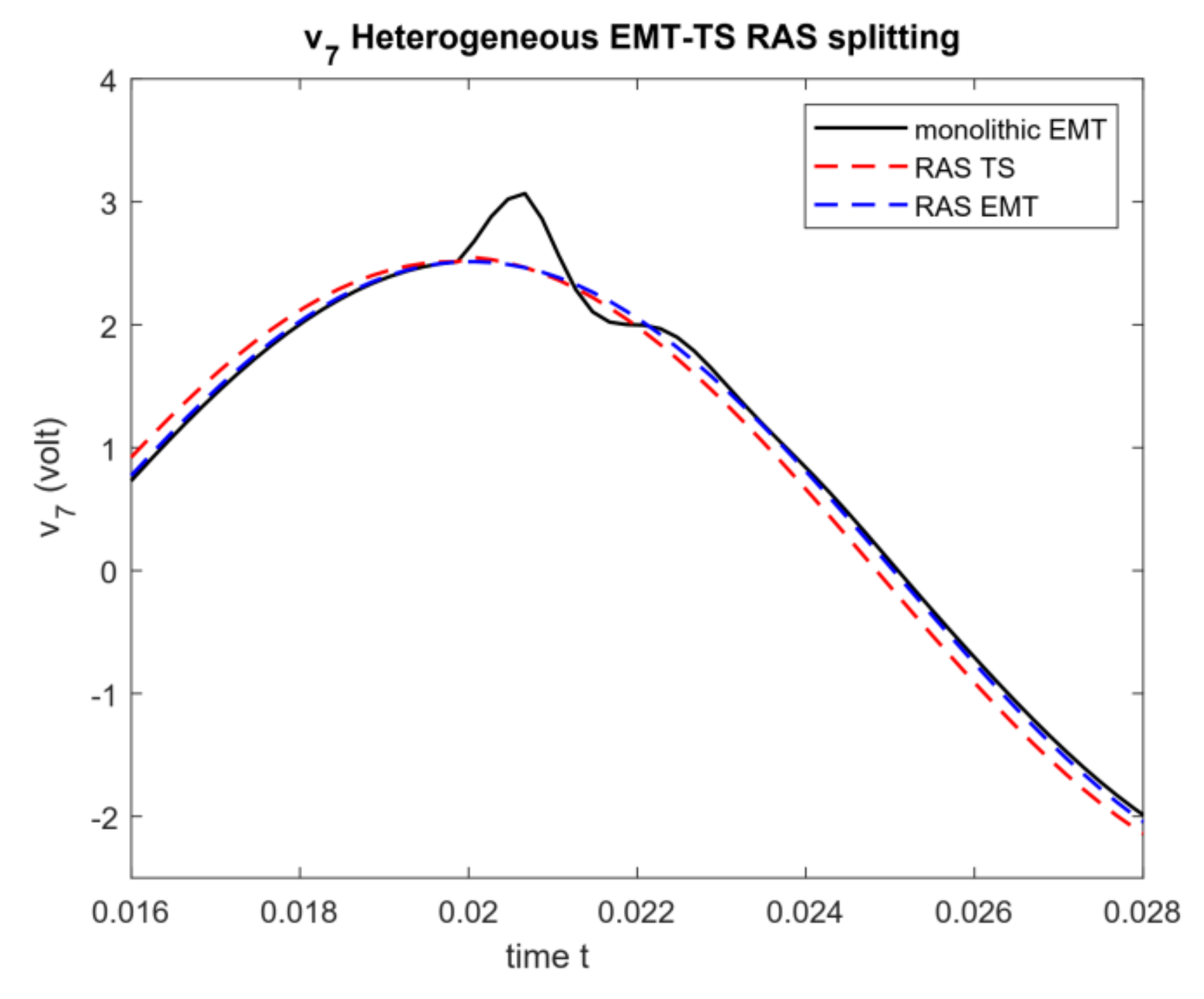}
\end{minipage}
\end{minipage}
\caption{Comparison of the behavior, with respect to time, of the variables $v_3$ (left) and $v_7$ (right) computed using the heterogeneous EMT-TS  RAS splitting with the Aitken's technique for accelerating convergence ($\Delta t_{ts}=2.10^{-3}$ and $\Delta t_{emt}=2.10^{-5}$), the reference is the monolithic EMT.  An amplitude perturbation on the voltage source. Parameters are $L_1=0.07, C_1=1.10^{-5}, R_1=7, L_2=0.07, C_2=1.10^{-7}, R_2=7, Zs=0.000001$.}
\label{DecoupageImpactFig}
\end{figure}
For the figure\ref{DecoupageImpactFig} a short event at the voltage source is created. Two potentials equidistant from the voltage source are compared $v_3$ and $v_7$, it is noted that the solutions found by the EMT subdomain during the co-simulation capture the disturbance much better and are closer to the monolithic solution for the potential $v_3$ compared to the potential $v_7$. Note that $v_3$ is close to the transition from EMT to TS, while $v_7$ is close to a boundary where information is transmitted only from TS to EMT. This gives us an indication of the loss of information implied by the TS part.

\section{Conclusion \label{shourick_contrib_Sec6}}
An iterative  co-simulation algorithm based on a restricted additive Schwarz heterogeneous DDM was used to co-simulate an RLC electrical circuit where part of the domain is modeled with EMT modeling and the other part with TS modeling. It integrates a non-intrusive Aitken's technique to accelerate convergence, possible thanks to  the pure linear convergence or divergence property of the RAS of the  heterogeneous DDM TS-EMT, with or without overlap.  Future developments at SuperGrid-Institute, which have already begun, will consist of building a software platform based on the co-simulation algorithm. The non-intrusive character of this co-simulation algorithm allows it to be implemented in a master-slave approach where it plays the role of a master in charge of orchestrating the call of EMT models and TS slaves launched as MPI (Message Passing Interface) processes. The slaves can be either Functional Mock-up Units (so that models can be developed with dedicated tools such as Modelica electrical libraries) or C++ code implementing the DAE function.

\section*{Acknowledgements}
This work was supported by a grant overseen by the French National Research Agency (ANR) as part of the “Investissements d’Avenir” Program ANE-ITE-002-01.

\bibliographystyle{plain}
\bibliography{biblio1}

\end{document}